\documentclass[11pt]{article}

\usepackage{bm}
\usepackage{amsfonts}
\usepackage{amssymb,amsmath,amsthm, amsfonts}

\usepackage{hyperref}
\usepackage[numbers,sort&compress]{natbib}
\usepackage[title, titletoc]{appendix}

\def\sqr#1#2{{\vcenter{\vbox{\hrule height.#2pt
              \hbox{\vrule width.#2pt height#1pt \kern#1pt \vrule width.#2pt}
          \hrule height.#2pt}}}}

\def\sqr#1#2{{\vcenter{\vbox{\hrule height.#2pt
              \hbox{\vrule width.#2pt height#1pt \kern#1pt \vrule width.#2pt}
              \hrule height.#2pt}}}}
\def\3n{\negthinspace \negthinspace \negthinspace }
\def\2n{\negthinspace \negthinspace }
\def\1n{\negthinspace }

\def\={\buildrel \triangle \over =}

\def\limsup{\mathop{\overline{\rm lim}}}

\def\esssup{\mathop{\rm esssup}}

\def\max{\mathop{\rm max}}
\def\min{\mathop{\rm min}}
\def\exp{\mathop{\rm exp}}
\def\sup{\mathop{\rm sup}}
\def\inf{\mathop{\rm inf}}

\def\({\Big (}
\def\){\Big )}
\def\[{\Big[}
\def\]{\Big]}
\def\be{\begin{equation}}

\def\ee{\end{equation}}

\def\square#1{\vbox{\hrule\hbox{\vrule height#1%
     \kern#1\vrule}\hrule}}
\def\rectangle#1#2{\vbox{\hrule\hbox{\vrule height#1%
     \kern#2\vrule}\hrule}}

\font\tenbb=msbm10 \font\sevenbb=msbm7 \font\fivebb=msbm5

\newfam\bbfam
\scriptscriptfont\bbfam=\fivebb \textfont\bbfam=\tenbb
\scriptfont\bbfam=\sevenbb

\newtheorem{lemma}{Lemma}[section]
\newtheorem{remark}{Remark}[section]
\newtheorem{example}{Example}[section]
\newtheorem{theorem}{Theorem}[section]

\newtheorem{definition}{Definition}[section]
\newtheorem{proposition}{Proposition}[section]
\newtheorem{assumption}{Assumption}[section]

\makeatletter
   
   \@addtoreset{equation}{section}
\makeatother

\begin{document}

\title{ Mean-field BSDEs with jumps under quadratic-exponential growth \thanks{ Juan Li is supported by the NSF of P.R. China (NOs. W2511002, 12031009), the NSF of Shandong Province (NO. ZR2023ZD35), National Key R and D Program of China (NO. 2018YFA0703900), and NSFC-RS (No. 11661130148; NA150344).}}

\author{Weimin Jiang$^{1}$,\,\,Juan Li$^{1,2}$\thanks{Corresponding authors at: School of Mathematics and Statistics, Shandong Uinversity, Weihai, Weihai 264290, P.~R.~China.},\,\,Yan Shen$^{1}$\footnotemark[2] \\
{$^1$\small School of Mathematics and Statistics, Shandong University, Weihai, Weihai 264209, P.~R.~China.}\\
{$^2$\small Research Center for Mathematics and Interdisciplinary Sciences, Shandong University,}\\
 {\small Qingdao 266237, P. R. China.} \\
{\footnotesize{\it Email: weiminjiang@mail.sdu.edu.cn,\ juanli@sdu.edu.cn,\ yanshen12@mail.sdu.edu.cn.}}
}
\date{November 9, 2026}
\maketitle

\textbf{Abstract}.
This paper studies the well-posedness of mean-field backward stochastic differential equations (BSDEs, for short) with Poisson jumps, where the generator has quadratic growth with respect to $Z$ and exponential growth with respect to $U$, where $(Y,Z,U)$ denotes the solution of this BSDE. First, we prove the existence and the uniqueness of local solutions to mean-field backward stochastic differential equations with jumps (BSDEJs, for short) by using a fixed-point argument. Then, based on these local results,  we discuss the solvability of global solutions in the following two cases: in the first one, the generator has linear growth with respect to $Y$ and the corresponding  marginal distribution $\mu_1$; and in the second one, the generator has subquadratic growth with respect to the marginal distribution $\mu_2$ and linear growth with respect to the marginal distribution $\mu_3$, corresponding to $Z$ and $U$, respectively. Finally, by employing a Girsanov transformation we establish a comparison theorem for mean-field BSDEJs under quadratic-exponential growth conditions.


\section{Introduction}
 The origin of linear BSDEs can be traced back to their seminal paper by Bismut \cite{B73} in 1973, cerning his work on the stochastic maximum principle.
 Later on,  Pardoux and Peng  \cite{PP90} pioneered the introduction of nonlinear BSDEs of the following form:
 \begin{equation}\label{eq1.1}
\begin{aligned}
 Y_t=\eta+\int_t^T g(s,Y_s,Z_s)ds-\int_t^T Z_sdB_s,\ \ t\in[0,T],
 \end{aligned}
\end{equation}
 where the terminal value $\eta:\Omega\rightarrow \mathbb{R}^n$ is an $\mathcal{F}_T$-measurable random variable, and the generator $g:\Omega\times[0,T]\times\mathbb{R}^n\times \mathbb{R}^{n\times d}\rightarrow \mathbb{R}^n$ is $\mathbb{F}$-progessively measurable for every fixed $(y,z)\in \mathbb{R}^n\times \mathbb{R}^{n\times d}$; $\mathbb{F}$ is the Brownian filtration generated by $B$.  They established the existence and the uniqueness  for solutions to equation $\eqref{eq1.1}$ when the generator satisfies a Lipschitz condition. This work has laid the foundational framework for the theory of BSDEs and has led to many important applications, including stochastic control  (Peng \cite{P92}), mathematical finance (El Karoui, Peng and Quenez \cite{KPQ97}), and partial differential equations (Pardoux and Peng \cite{PP92}).   Since then, many authors have worked on relaxing the Lipschitz condition and extending the well-posedness theory of BSDEs. For example, Lepeltier and San Martin \cite{LM97} established the existence of adapted solutions for BSDEs with continuous generators of linear growth, while Fan and Jiang \cite{FJ13} investigated multidimensional BSDEs under weak monotonicity and general growth assumptions.

As a generalization of BSDEs, those  with quadratic growth have attracted considerable attention because of their wide range of applications in nonlinear financial problems. Such equations arise ina natural in many classical models, including risk-sensitive control (Kohlmann and Tang \cite{KT03}) and utility maximization (Hu, Imkeller, and M$\mathrm{\ddot{u}}$ller \cite{HIM05}). The study of quadratic BSDEs was initiated with the work of Kobylanski \cite{K00}. In 2000, she considered one-dimensional quadratic BSDEs with bounded terminal values and proved the existence and the uniqueness of solutions. Since then, many important results have been obtained. So Briand and Hu \cite{BH06} extended the existence result to the case of unbounded terminal values. Later, Briand and Hu \cite{BH08} established uniqueness for generators that are convex in $z$ and terminal values that are unbounded. They also derived a nonlinear Feynman-Kac formula connected with PDEs. In 2020, Fan, Hu and Tang \cite{FHT20} further extended the uniqueness result to non-convex generators. For multidimensional quadratic BSDEs, Hu and Tang \cite{HT16} studied equations with diagonally quadratic generators and investigated both local and global solvability. In recent years, quadratic growth conditions have also been considered in several extensions of BSDEs, including anticipated BSDEs \cite{HLW21} and backward doubly stochastic differential equations (BDSDEs) \cite{HWX24}.

The classical BSDE theory has been well developed in the continuous diffusion setting. However, it is not suitable for modeling discontinuous risks caused by sudden events such as defaults, natural disasters, and market crashes. To capture these jump risks, Poisson random measures are introduced, leading to the more general class of BSDEJs. Compared with classical BSDEs, BSDEJs contain an additional jump component $U$, which is used to describe and hedge discontinuous risks. As a result, BSDEJs have become an important tool in credit risk pricing, insurance mathematics, and portfolio optimization in markets with jumps (Lim and Quenez \cite{LQ11}). The study of BSDEJs began with the work of Tang and Li \cite{TL94}, who established the existence and the uniqueness of solutions under Lipschitz assumptions. Later, Barles, Buckdahn and Pardoux \cite{BBP97} investigated the well-posedness of BSDEJs and their links with integro-PDEs. Situ \cite{S97} further extended these results to non-Lipschitz generators. Compared with the continuous case, the theory of quadratic BSDEJs is still less developed. Early studies mainly arose from utility maximization and indifference pricing problems in jump models (Becherer \cite{B06}). Since around 2010, more attention has been paid to the theoretical analysis of quadratic BSDEJs. For example, Kazi-Tani, Possamai and Zhou \cite{KPZ15} studied BSDEJs whose generators are uniformly Lipschitz in $Y$, have quadratic growth in $Y$, and exponential growth in $U$, and they proved the existence and the uniqueness of solutions. Fujii and Takahashi \cite{FT18} further established well-posedness and Malliavin differentiability under a general quadratic-exponential framework.

The study of mean-field stochastic differential equations (McKean-Vlasov SDEs) has its origin in statistical physics. Its early development can be traced back to the work of Kac \cite{K56} on the propagation of chaos, where mean field limits were used to describe systems with a large number of interacting particles. Later, Lasry and Lions \cite{LL07} introduced the theory of mean field games from the perspective of game theory. At about the same time, Huang, Caines and Malham$\mathrm{\acute{e}}$ \cite{HMC06} independently developed a related framework based on the Nash certainty equivalence principle in engineering control. Motivated by these developments, Buckdahn, Djehiche, Li and Peng \cite{BDLP09} introduced mean-field BSDEs and studied a class of mean-field problems by probabilistic methods. Soon afterwards, Buckdahn, Li and Peng \cite{BLP09} established the existence and the uniqueness of solutions, proved a comparison theorem, and provided a probabilistic interpretation of the associated nonlocal PDEs under Lipschitz assumptions on the generator. Since then, the theory of mean-field SDEs and BSDEs has developed rapidly. We refer to Buckdahn, Li, Peng and Rainer \cite{BLPR17}, Li \cite{L18}, and Li, Liang and Zhang \cite{LLZ18} for some representative results. More recently, Jiang, Li and Wei \cite{JLW24} established well-posedness results for multidimensional mean-field BSDEs with diagonally quadratic generators. Hao, Hu, Tang and Wen \cite{HHTW25} studied one-dimensional mean-field BSDEs with quadratic growth and bounded terminal conditions, and investigated the related nonlocal PDEs.

Nasroallah and Ouknine  \cite{NO25} investigated the well-posedness of  BSDEJs with central value reflection under quadratic-exponential growth conditions, where the generator depends on the law of $Y$. Becherer and Hesse \cite{BH26} considered a competitive investment and hedging problem with the common noise contains random measure, and characterized the mean-field equilibrium.

Motivated by the above works, we investigate the following mean-field BSDEJs with quadratic-exponential growth:
\begin{equation}\label{eq1.2}
\begin{aligned}
Y_t=\eta+\int_t^T g\big(s,Y_s,Z_s,U_s^{\rho}(\cdot), \mathbb{P}_{(Y_s,Z_s,U_s^{\rho}(\cdot))}\big)ds
-\int_t^T Z_sdB_s
-\int_t^T\!\int_E U_s(x)\tilde{N}(ds,dx),\  t\in [0, T],
\end{aligned}
\end{equation}
where the terminal value $\eta:\Omega\rightarrow \mathbb{R}$ is an $\mathcal{F}_T$-measurable random variable, the generator $g:\Omega\times[0,T]\times\mathbb{R}\times \mathbb{R}^{d}\times\mathbb{R}\times \mathcal{P}_2(\mathbb{R}^{1+d+1})\rightarrow \mathbb{R}$ is predictable for every fixed $(y,z,h,\mu)\in \mathbb{R}^{1+d+1}\times \mathcal{P}_2(\mathbb{R}^{1+d+1})$ and $U_s^{\rho}(\cdot):=\int_EU_s(x)\rho_s(x)\lambda(dx)$.
This paper studies mean-field BSDEJs with bounded terminal conditions and a generator of quadratic-exponential growth. We establish the local and the global well-posedness of these equations and prove a comparison theorem. First, we prove the existence and the uniqueness of local solutions by a fixed-point argument (Theorem \ref{th3.1}). Based on this local result, we then establish global well-posedness of equation \eqref{eq1.2} in two different settings. In the first setting, the generator $g(t,y,z,u^{\rho_t}(\cdot),\mu)$ has a quadratic growth in $z$, an exponential growth in $u$, and a linear growth in both $y$ and the first marginal distribution $\mu_1$ (Theorem \ref{th4.1}). In the second one, the generator has a quadratic growth in $z$, an exponential growth in $u$, a sub-quadratic growth in the second marginal distribution $\mu_2$, and a linear growth in the third marginal distribution $\mu_3$, with sufficiently small growth coefficients associated with $\mu_2$ and $\mu_3$ (Theorem \ref{th4.2}). We also investigate a class of mean-field BSDEJs whose generator has the additive form $f(t,z,u,\nu):=f_1(t,z,u)+f_2(t,\nu),$ for $(t,z,u,\nu)\in[0,T]\times\mathbb{R}^d\times \mathcal{L}^2(\lambda)\times\mathcal{P}_2(\mathbb{R}^{1+d}\times \mathcal{L}^2(\lambda))$, and we establish a corresponding well-posedness result (Theorem \ref{th4.3}). Finally, by means of a Girsanov transformation  we prove a comparison theorem for mean-field BSDEJs with a quadratic-exponential growth generator (Theorem \ref{th5.1}).

The main contributions of this paper are threefold. First, we establish the existence and the uniqueness of local solutions to mean-field BSDEJs under quadratic-exponential growth, where the generator depends on the joint distributions of $Y,\ Z$ and $U^\rho$. Second, based on this local result, we prove global well-posedness in two cases. In the first case, the generator has linear growth in $Y$ and its marginal distribution. In the second one,  the generator has subquadratic growth in the  marginal distribution of $Z$ and linear growth in the  marginal distribution of $U^\rho$,  with sufficiently small growth coefficients. We also investigate a class of generators with an additive form. In this case, the generator  depends on the law of $Y$ and the law of the $\mathcal{L}^2(\lambda)$-value jump integrand $U$, and the quadratic growth coefficient associated with the law of $Z$ need not be sufficiently small. Third, by using a Girsanov transformation we establish a comparison theorem for mean-field BSDEJs whose generator depends only on law of $Y$ and which satisfy the corresponding monotonicity condition.

The rest of this paper is organized as follows.  Section 2 introduces the function spaces, notation, and  properties of BMO martingales.  Section 3 establishes  local well-posedness, and  Section 4 studies  global solvability  in the stated growth assumptions and the additive setting.  Section 5 establishes a comparison theorem for bounded solutions.  Appendix A gives the proof of  the auxiliary existence result used in Proposition 3.1.

\section{Preliminaries}
Let $T>0$ be an arbitrarily  finite time horizon. Let $(\Omega ,\mathcal{F}, \mathbb{P})$ be a complete probability space on which is defined a $d$-dimensional Brownian motion $B=\{ B_t\}_{ t\in[0,T]}$, and an independent Poisson random measure $N$ on $\mathbb{R}^{+}\times E.\ E=\mathbb{R} \backslash \{0\}$ is equipped with the $\sigma$-field $\mathcal{B}(E)$. The compensator $\nu(dt, dx) = \lambda(dx)dt$ of $N$ is such that $\big\{ \tilde{N}([0, t]\times K) = (N-\nu)([0, t]\times K) \big\}_{t\geq 0}$ is a martingale for all $K\in \mathcal{B}(E)$ satisfying $\lambda(K) < +\infty$, where $\lambda$ is a given $\sigma$-finite L$\acute{\mathrm{e}}$vy measure on $(E, \mathcal{B}(E))$, i.e., a measure on $(E, \mathcal{B}(E))$ such that $\int _{E} (1\wedge |x|^{2})\lambda(dx)< +\infty$. We write $\textbf{1}_A(\cdot)$ for the indicator of set $A$. 
The notation $\delta_0$ denotes the Dirac measure at the origin. We assume that there is a sub-$\sigma$-field $\mathcal{G}^0\subset\mathcal{F}$ such that\\
\indent(i) the Brownian motion $B$ and the Poisson random measure $N$ are independent of $\mathcal{G}^0$;\\
\indent(ii)\ $\mathcal{G}^{0}$ is ``rich enough", i.e., $\mathcal{P}_2(\mathbb{R}^k)=\big\{\mathbb{P}_\varsigma,\ \varsigma\in L^2_{\mathcal{G}^0}(\Omega,\mathbb{R}^k)\big\}$, $k\geq1$, where $\mathbb{P}_\varsigma:=\mathbb{P}\circ{[\varsigma]^{-1}}$ denotes the law of random variable\ $\varsigma$ under probability\ $\mathbb{P}$; \\
\indent(iii)\ $\mathcal{G}^{0}\supset \mathcal{N}_{\mathbb{P}}$, where $\mathcal{N}_{\mathbb{P}}$ is the set of all $\mathbb{P}$-null subsets of $\mathcal{F}$.\\
By $\mathbb{F} = \{ \mathcal{F}_t ,\ 0\leq t\leq T  \}$ we denote the filtration generated by the Brownian motion $B$ and the Poisson
random measure $N$, augmented by $\mathcal{G}^{0}$, i.e.,
$$\mathcal{F}^{0}_{t}:=\sigma \big\{B_s,\ N([0, s]\times K)|\ s\leq t,\ K\in \mathcal{B}(E)\big\},$$
$$\mathcal{F}_{t}:=\mathcal{F}^{0}_{t+}\vee \mathcal{G}^{0}\Big(=\Big(\bigcap \limits_{s>t} \mathcal{F}^{0}_{s}\Big)\vee \mathcal{G}^{0}\Big),\ \ \ t\in [0,T].$$
Let us introduce the following spaces which are needed in what follows. For $p\geq 1$, $0\leq s\leq t\leq T$,\\
\noindent$\bullet$\ $L^p_{\mathcal{F}_{t}}(\Omega;\mathbb{R}) :=\big\{\phi:\Omega \rightarrow \mathbb{R}\ |\  \phi\ \text{is an}\ {\mathcal{F}_{t}}\text{-measurable r.v. with}\ \|\phi \|^{p} _{L^p(\Omega)}:=\mathbb{E}\big[|\phi|^p\big]<+\infty\big\}$;

\noindent$\bullet$\ $L^\infty_{\mathcal{F}_{t}}(\Omega;\mathbb{R}) :=\big\{\phi: \Omega \rightarrow \mathbb{R}\ |\ \phi\ \text{is an}\ {\mathcal{F}_{t}}\text{-measurable r.v. with} \ \|\phi \| _{\infty}:=\esssup _{\omega\in \Omega}|\phi(\omega)|<+\infty\big\}$;

\noindent$\bullet$\ $\mathbb{S}^p(s,t;\mathbb{R}) :=\big\{\phi: \Omega \times [s,t]\rightarrow \mathbb{R}\ |\ \phi\ \text{is an}\ \mathbb{F}\text{-adapted}\ \text{c}\grave{\mathrm{a}}\text{dl}\grave{\mathrm{a}}\text{g}\ \text{process such that}\ \| \phi \|^{p} _{\mathbb{S}^p(s,t)}:=\mathbb{E}\big[\sup_{r\in [s,t]}|\phi_r|^p\big]<+\infty\big\}$;

\noindent$\bullet$\ $\mathbb{S}^\infty(s,t;\mathbb{R}) :=\big\{\phi: \Omega \times [s,t]\rightarrow \mathbb{R}\ |\ \phi\ \text{is an}\ \mathbb{F}\text{-adapted}\ \text{c}\grave{\mathrm{a}}\text{dl}\grave{\mathrm{a}}\text{g}\ \text{process with}\ \| \phi\| _{\mathbb{S}^\infty(s,t)}:=\|\sup_{r\in [s,t]}|\phi_r|\|_{\infty}<+\infty\big\}$;

\noindent$\bullet$\ $\mathbb{H}^p(s,t;\mathbb{R}^d):=\big\{\phi: \Omega \times [s,t]\rightarrow \mathbb{R}^d\ |\ \phi \ \text{is an}\ \mathbb{F}\text{-progressively measurable process with}\ \| \phi \|^p _{\mathbb{H}^p(s,t)}:= \mathbb{E}\big[\big(\int_s^t\!|\phi_r|^2dr\big)^{\frac{p}{2}}\big]<+\infty\big\}$;

\noindent$\bullet$\ $\mathcal{L}^p(\lambda):=\big\{U:\ E\rightarrow \mathbb{R}\ |\ U\ \text{is }\ \mathcal{B}(E)\text{-measurable with}\ \| U \|^p _{\mathcal{L}^p(\lambda)}:= \int_E\!|U(x)|^p\lambda(dx)<+\infty\big\}$;

\noindent$\bullet$\ $\mathcal{L}^\infty(\lambda):=\big\{U: E\rightarrow \mathbb{R}\ |\ U\ \text{is }\  \mathcal{B}(E)\text{-measurable  with}\ \| U \| _{\mathcal{L}^\infty(\lambda)}:= \esssup_{x\in E}|U(x)|<+\infty \big\};$

\noindent$\bullet$\ $\mathbb{J}^p(s,t;\mathbb{R}):=\big\{U: \Omega \times [s,t]\times E \rightarrow \mathbb{R}\ |\ U\  \text{is }\ \mathcal{P}^{0}\otimes\mathcal{B}(E)\text{-measurable  process with}\ \| U \|^p _{\mathbb{J}^p(s,t)}:= \mathbb{E}\big[\big(\int_s^t\int_E\!|U_r(x)|^2\lambda(dx)dr\big)^{\frac{p}{2}}\big]<+\infty\big\}$;

\noindent$\bullet$\ $\mathbb{J}^\infty(s,t;\mathbb{R}):=\big\{U: \Omega \times [s,t]\times E \rightarrow \mathbb{R}\ |\ U\ \text{is }\ \mathcal{P}^{0}\otimes\mathcal{B}(E)\text{-measurable process with}\ \| U \| _{\mathbb{J}^\infty(s,t)}:=\|\esssup_{(\omega,r,x)\in  \Omega\times[s,t]\times E}|U_r(\omega,x)|<+\infty\big\}$;

\noindent$\bullet$\ $\mathbb{H}^2_{BMO}(s,t;\mathbb{R}^d):=\big\{\phi: \Omega \times [s,t]\rightarrow \mathbb{R}^d\ |\ \phi\ \text{is an}\ \mathbb{F}\text{-progressively measurable process with}\\ \| \phi \|^2 _{\mathbb{H}^2_{BMO}(s,t)}:= \sup _{\tau\in \mathcal{T}_{s}^{t}}\big\|\mathbb{E}_{\tau}\big[\int_\tau^t|\phi_r|^2dr\big]\big\| _{\infty}<+\infty\big\}$;

\noindent$\bullet$\ $\mathbb{J}^2_{B}(s,t;\mathbb{R}):=\big\{\phi: \Omega \times [s,t]\times E \rightarrow \mathbb{R}\ |\ \phi\ \text{is }\ \mathcal{P}^{0}\otimes\mathcal{B}(E)\text{-measurable process with}\  \| \phi \|^2 _{\mathbb{J}^2_{B}(s,t)}:= \sup _{\tau\in \mathcal{T}_{s}^{t}}\big\|\mathbb{E}_{\tau}\big[\int_\tau^t\int_E|\phi_r(x)|^2\lambda(dx)dr\big]\big\| _{\infty}<+\infty\big\}$;

\noindent$\bullet$\ $\mathbb{J}^2_{BMO}(s,t;\mathbb{R}):=\big\{\phi: \Omega \times [s,t]\times E\rightarrow \mathbb{R}\ |\ \phi\ \text{is }\ \mathcal{P}^{0}\otimes\mathcal{B}(E)\text{-measurable process with}\ \| \phi \|^2 _{\mathbb{J}^2_{BMO}(s,t)}:= \sup _{\tau\in \mathcal{T}_{s}^{t}}\big\|\mathbb{E}_{\tau}\big[\int_\tau^t\int_E|\phi_r(x)|^2\lambda(dx)dr\big]+(\Delta M_{\tau})^{2}\big\| _{\infty}<+\infty\big\}$,

\noindent where\ $\mathcal{P}^{0}$ denotes the $\sigma$-field of $\mathbb{F}$-predictable subsets on $\Omega\times [0, T ]$, $ \mathcal{T}_{s}^{t}$ is the set of all $\mathbb{F}$-stopping times $\tau$ with values in $[s,t]$, $\mathbb{E}_{\tau}$ is the conditional expectation with respect to
 $\sigma$-field $\mathcal{F}_{\tau}$ and $\Delta M_{\tau}:=\int_E\phi_\tau(x)N(\{\tau\},dx).$ We denote by $\mathbb{E}^{\mathbb{Q}}$ and $\mathbb{E}^{\mathbb{Q}}_{\tau}$ the expectation and the conditional expectation under the probability measure $\mathbb{Q}$, respectively.

 Moreover, let $M=(M_t)_{t\in[0,T]}$  be a uniformly integrable martingale with $M_0 = 0$, and define $\mathcal{E}(M)=(\mathcal{E}(M)_t)_{t\in[0,T]}:=\big(e^{(M_t-\frac{1}{2}\langle M^c\rangle_t)}\prod_{0\leq s \leq t}(1+\Delta M_s)e^{-\Delta M_s}\big)_{t\in[0,T]}$ as the Dol$\acute{\mathrm{e}}$ans-Dade exponential of $M$, where $M^c$  represents the continuous part of $M$. For $p\geq 1,$  we set
$$\|M\|_{BMO_p(\mathbb{P})}:=\sup \limits_{\tau\in \mathcal{T}_0^T}\Big\|\mathbb{E}_{\tau}\big[\big|M_T-M_{\tau-}\textbf{1}_{\tau>0}\big|^p\big]^{\frac{1}{p}}\Big\| _{\infty}.$$
The class $\big\{M : \|M\|_{BMO_p(\mathbb{P})}< +\infty\big\}$ is denoted by $BMO_p(\mathbb{P})$. Observe that $BMO_p(\mathbb{P})$ is a Banach space equipped with the norm $\|\cdot\|_{BMO_p(\mathbb{P})}$. In the sequel for simplicity, we write $BMO(\mathbb{P})$ for the space of $BMO_2(\mathbb{P})$. We recall that the process $t \rightarrow \int_0^t Z_s dB_s$ (denoted by $Z \cdot B$) belongs to $BMO(\mathbb{P})$ if and only if $Z \in \mathbb{H}^2_{BMO}(0,T;\mathbb{R}^d)$; likewise, the process $t \rightarrow \int_0^t \int_E U_s(x) \tilde{N}(ds,dx)$ (denoted by $U \cdot \tilde{N}$) belongs to $BMO(\mathbb{P})$ if and only if $U \in \mathbb{J}^2_{BMO}(0,T;\mathbb{R})$.  Moreover, in this case,
\begin{equation*}
\begin{split}
\|Z\cdot B\|^2_{BMO(\mathbb{P})}=\| Z \|^2 _{\mathbb{H}^2_{BMO}(0,T)},\ \ \ \|U\cdot \tilde{N}\|^2_{BMO(\mathbb{P})}=\| U \|^2 _{\mathbb{J}^2_{BMO}(0,T)}.
\end{split}
\end{equation*}
 Let $ \mathcal{P}_{2}(\mathbb{R}^d)$ denote the set of all probability measures $\mu$ on $(\mathbb{R}^d, \mathcal{B}(\mathbb{R}^d))$ with finite
 2-order moment, that is, $\int_{\mathbb{R}^d}|x|^{2}\mu(dx)< +\infty$.\ The space $\mathcal{P}_2(\mathbb{R}^d)$ is endowed with the $2$-Wasserstein metric: For $ \mu,\ \nu\in\mathcal{P}_2(\mathbb{R}^d)$,
\begin{equation*}
\begin{split}
W_2(\mu,\nu):=\inf\Big\{&\Big(  \int_{\mathbb{R}^d\times\mathbb{R}^d}|x-y|^2\rho(dxdy)\Big)^{\frac{1}{2}},\ \ \rho\in\mathcal{P}_2(\mathbb{R}^{2d})\ \text{with marginals}\ \mu\ \text{and}\ \nu    \Big\}.
\end{split}
\end{equation*}
For our discussion which follows, we introduce the following definitions concerning the solutions of  equation \eqref{eq1.2}.
\begin{definition}
\normalfont
 A triple of processes $(Y, Z, U)\in \mathbb{S}^{2}(0,T;\mathbb{R})\times \mathbb{H}^{2}(0,T;\mathbb{R}^d)\times \mathbb{J}^{2}(0,T;\mathbb{R}) $ is called an adapted
solution of  mean-field BSDEJs \eqref{eq1.2}, if $\mathbb{P}$-almost surely, it satisfies \eqref{eq1.2}. Furthermore, it is called a bounded  adapted solution if it also belongs to $\mathbb{S}^{\infty}(0,T;\mathbb{R})\times \mathbb{H}^{2}_{BMO}(0,T;\mathbb{R}^d)\times \mathbb{J}^{2}_{BMO}(0,T;\mathbb{R}).$
\end{definition}
We present three technical lemmas that will be used repeatedly in the subsequent analysis.  The proofs of the results below are similar to those of Lemma 2.2 in \cite{FT19} and Corollary 1 in \cite{M09}.
\begin{lemma}\label{le2.1}
If $(Y, Z, U)\in \mathbb{S}^{\infty}(0,T;\mathbb{R})\times \mathbb{H}^{2}_{BMO}(0,T;\mathbb{R}^d)\times \mathbb{J}^{2}_{BMO}(0,T;\mathbb{R})$ is a solution of the mean-field BSDEJs \eqref{eq1.2}, then we have

\noindent$\mathrm{(i)}\ \big(\|U\|^{2}_{\mathbb{J}_{B}^2(0,T)}\vee \|U\|^{2}_{\mathbb{J}^\infty(0,T)}\big)\leq \|U\|^{2}_{\mathbb{J}^2_{BMO}(0,T)}\leq \|U\|^{2}_{\mathbb{J}_{B}^2(0,T)}+ \|U\|^{2}_{\mathbb{J}^\infty(0,T)};$\\
\noindent $\mathrm{(ii)}$ We have \ $\|U_s\|_{\mathcal{L}^\infty(\lambda)}\leq 2\|Y\|_{\mathbb{S}^\infty(0,T)},\ d\mathbb{P}\times ds\mbox{-a.e.}$ In particular, $U$ has a bounded predictable version, which we use throughout.
\end{lemma}
The following result, known as  Kazamaki's criterion (see \cite[Remark 3.1]{K94}), establishes a connection between the martingale property of the stochastic exponential $\mathcal{E}(M)$ and a BMO condition on the martingale $M$.
\begin{lemma}\label{le2.2}
Let $0 < \delta < +\infty$ and M be a BMO-martingale such that
 $\Delta M_{t}\geq -1+\delta, \ \mathbb{P}\mbox{-a.s.}, \  t\in[0, T],$
then $\mathcal{E}(M)$ is a uniformly integrable martingale.
\end{lemma}
The following result characterizes the equivalence of BMO norms between different probability measures under Girsanov transformations.
\begin{lemma}\label{le2.3}
Let $\kappa>0$ and $\delta\in(0,1)$. Suppose that $K:=\Lambda\cdot B+\Theta \cdot \tilde{N}$ satisfies $\|K\|_{BMO(\mathbb{P})}\le\kappa$ and $ |\Theta_s(x)|\le 1-\delta.$  Then, for every $M:=Z\cdot B+U\cdot\tilde N\in BMO(\mathbb{P}),$  there exist constants $c_1,\ c_2>0$, depending only on $\kappa$ and $\delta$, such that
$$c_{1}\|M\|_{BMO(\mathbb{P})}\leq \|\widehat{M}\|_{BMO(\mathbb{Q})}\leq c_{2}\|M\|_{BMO(\mathbb{P})},$$
where $d\mathbb{Q}:=\mathcal{E}(K)_Td\mathbb{P}$ and
$$\widehat{M}_t:=M_t-\int_0^t Z_s\Lambda_sds-\int_0^t\int_EU_s(x)\Theta_s(x)\lambda(dx)ds,\ \  \ t\in[0,T].$$
\end{lemma}
\begin{proof}
Since $K\in BMO(\mathbb{P})$ and $\Delta K\ge-1+\delta$, the process $L:=\mathcal{E}(K)$ is a positive uniformly integrable martingale. Thus, measure
$\mathbb{Q}$ is equivalent to $\mathbb{P}$. By Kazamaki's isomorphism theorem \cite[Theorem]{K79}, we know that there exist constants $a,\ b>0$, depending only on $\kappa$ and $\delta$, such that
\begin{equation}\label{eq2.1}
\begin{aligned}
a\|X\|_{BMO(\mathbb{P})}\le\|\Phi_K(X)\|_{BMO(\mathbb{Q})}\le b\|X\|_{BMO(\mathbb{P})},
\end{aligned}
\end{equation}
where
$$\Phi_K(X)_t:=X_t-\int_0^t L_s^{-1} d[X,L]_s,\ \  \ t\in[0,T].$$
Define $\widetilde{M}:=Z\cdot B+((1+\Theta)U)\cdot \tilde{N}.$ Since $\delta\le 1+\Theta\le 2-\delta$, we obtain $\widetilde{M}\in BMO(\mathbb{P})$ and
\begin{equation}\label{eq2.2}
\begin{aligned}
\delta\|M\|_{BMO(\mathbb{P})} \le \|\widetilde{M}\|_{BMO(\mathbb{P})} \le (2-\delta)\|M\|_{BMO(\mathbb{P})}.
\end{aligned}
\end{equation}
By direct calculation, we obtain
$$
\int_0^t L_s^{-1}d[\widetilde{M},L]_s=\int_0^t Z_s\Lambda_sds+\int_0^t\int_E U_s(x)\Theta_s(x)N(ds,dx),\ \ \ t\in[0,T].
$$
Therefore, we have
\begin{equation*}
\begin{aligned}
\Phi_K(\widetilde{M})_t=M_t-\int_0^t Z_s \Lambda_sds-\int_0^t\int_E U_s(x)\Theta_s(x)\lambda(dx)ds=\widehat{M}_t,\ \ \ t\in[0,T].
\end{aligned}
\end{equation*}
Combining equations \eqref{eq2.1} and \eqref{eq2.2}, we obtain
$$
a\delta\|M\|_{BMO(\mathbb{P})} \le \|\widehat M\|_{BMO(\mathbb{Q})} \le b(2-\delta)\|M\|_{BMO(\mathbb{P})}.
$$
The proof is complete.
\end{proof}

\section{The existence and the uniqueness of the local solution}
In this section, our goal is to prove the existence and the uniqueness of local solutions for  mean-field BSDEJs \eqref{eq1.2} under quadratic-exponential growth,
we  rewrite the equation \eqref{eq1.2} as follows:
\begin{equation}\label{eq3.1}
\begin{aligned}
Y_t=\eta+\int_t^T\! g\big(s,Y_s,Z_s,U_s^{\rho_s}(\cdot), \mathbb{P}_{(Y_s,Z_s,U_s^{\rho_s}(\cdot))}\big)ds
-\int_t^T\! Z_sdB_s
-\int_t^T \!\int_E \!U_s(x)\tilde{N}(ds,dx),\ t\in [0, T],
\end{aligned}
\end{equation}
where $U_s^{\rho_s}(\cdot):=\int_EU_s(x)\rho_s(x)\lambda(dx)$.

In what follows, we  assume  that $\theta: \Omega \times[0,T] \rightarrow \mathbb{R}^{+}$ is $\mathbb{F}$-progressively measurable, $\phi$ and $ \phi_0:[0,\infty)\rightarrow[0,\infty)$ are two nondecreasing continuous functions,  $\beta_0,\ \beta_1,\ \gamma,\ \gamma_0$ and $K_0$ are positive constants and $\alpha\in [0,1)$. For simplicity, we denote the marginal distributions of $\mu\in \mathcal{P}_{2}(\mathbb{R}^{1+d+1})$ by $\mu _{1}:=\mu(\cdot\times \mathbb{R}^{d+1}),\ \mu _{2}:=\mu( \mathbb{R}\times\cdot\times\mathbb{R})$ and $ \mu _{3}:=\mu( \mathbb{R}^{1+d}\times\cdot)$. For $u,\ \rho\in \mathcal{L}^{2}(\lambda)$, we introduce the notation  $u^{\rho}(\cdot):=\int_Eu(x)\rho(x)\lambda(dx)$ and the nonnegative functional $[u]_{\gamma}:=\int_E\frac{e^{\gamma u(x)}-\gamma u(x)-1}{\gamma}\lambda(dx)$.

 Let $\eta: \Omega \rightarrow \mathbb{R}$  be $\mathcal{F}_T$-measurable,  $g:=g(\omega,t,y,z,h,\mu):\Omega \times[0,T]\times \mathbb{R}\times\mathbb{R}^d \times \mathbb{R}\times \mathcal{P}_{2}(\mathbb{R}^{1+d+1})\rightarrow \mathbb{R}$  be $\mathcal{P}^0$-measurable, for every fixed $(y,z,h,\mu)\in \mathbb{R}\times\mathbb{R}^d \times \mathbb{R}\times \mathcal{P}_{2}(\mathbb{R}^{1+d+1}),$ and $ \rho:\Omega \times[0,T]\times E\rightarrow \mathbb{R}$ be $\mathcal{P}^0\otimes\mathcal{B}(E)$-measurable. Furthermore, we make the following assumptions.
\begin{assumption}\label{as3.1}
\normalfont
 $\mathrm{(i)}$
  For all $t\in [0,T],\ y\in \mathbb{R},\ z\in \mathbb{R}^d,\ u\in \mathcal{L}^{2}(\lambda),$ and $ \mu\in \mathcal{P}_{2}(\mathbb{R}^{1+d+1})$,
\begin{equation*}
\begin{aligned}
-[-u]_{\gamma}-f(t,y,z,\mu)\leq g\big(t,y,z,u^{\rho_t}(\cdot),\mu\big)\leq [u]_{\gamma}+f(t,y,z,\mu),\ \ d\mathbb{P}\times dt\mbox{-a.e.},
\end{aligned}
\end{equation*}
where $$f(t,y,z,\mu):=\theta _{t}+\phi( |y|)+\frac{\gamma}{2}|z|^{2}+\phi_0\big(W_2(\mu _{1},\delta _{0})\big)+\gamma _{0}W_2(\mu _{2},\delta _{0})^{1+\alpha}+\beta _{0}W_2(\mu _{3},\delta _{0}).$$
$\mathrm{(ii)}$
For all $t\in [0,T],\ y,\ \bar{y}\in \mathbb{R},\ z,\ \bar{z}\in \mathbb{R}^d,\ u\in \mathcal{L}^{2}(\lambda), $ and $ \mu,\  \bar{\mu}\in \mathcal{P}_{2}(\mathbb{R}^{1+d+1})$,
\begin{equation*}
\begin{aligned}
&\big|g\big(t,y,z,u^{\rho_t}(\cdot),\mu\big)-g\big(t,\bar{y},\bar{z},u^{\rho_t}(\cdot),\bar{\mu}\big)\big|\\
\leq\ & \Big[\big(1+|z|+|\bar{z}|+W_2(\mu _{2},\delta _{0})+W_2(\bar{\mu} _{2},\delta _{0})\big)\big(|y-\bar{y}|+|z-\bar{z}|+W_2(\mu _{1},\bar{\mu}_{1})\big)\\
& +\big(1+|z|^\alpha+|\bar{z}|^\alpha+W_2(\mu _{2},\delta _{0})^{\alpha}+W_2(\bar{\mu} _{2},\delta _{0})^{\alpha}\big)W_2(\mu _{2},\bar{\mu} _{2})+W_2(\mu _{3},\bar{\mu} _{3})\Big]\\
&\times \phi\big(|y|\vee|\bar{y}|\vee W_2(\mu _{1},\delta _{0})\vee W_2(\bar{\mu} _{1},\delta _{0})\big),\ \ \ d\mathbb{P}\times dt\mbox{-a.e.}
\end{aligned}
\end{equation*}
$\mathrm{(iii)}$\
For all $t\in [0,T],$ and $ x\in E$,
$$ |\rho_t(x)|\leq K_0(1\wedge |x|),\ \ d\mathbb{P}\times dt\mbox{-a.e.}$$
$\mathrm{(iv)}$ For each $M>0$, there exist constants $K_M\geq0$ and  $\delta_M\in(0,1)$ such that, $d\mathbb{P}\times dt\mbox{-a.e.}$, for all $t\in [0,T],\ y\in \mathbb{R},\ z\in \mathbb{R}^d,\ u,\ \bar{u}\in \mathcal{L}^{2}(\lambda)\cap \mathcal{L}^\infty(\lambda)$ with
$\|u\|_{\mathcal{L}^\infty(\lambda)}\vee\|\bar u\|_{\mathcal{L}^\infty(\lambda)}\le M,$ and $ \mu\in \mathcal{P}_{2}(\mathbb{R}^{1+d+1})$,
\begin{equation*}
\begin{aligned}
\big|g\big(t,y,z,u^{\rho_t}(\cdot),\mu\big)-g\big(t,y,z,\bar{u}^{\rho_t}(\cdot),\mu\big)\big|\leq K_M|u^{\rho_t}(\cdot)-\bar{u}^{\rho_t}(\cdot)|.
\end{aligned}
\end{equation*}
Moreover, for all $t\in[0,T]$, and  $x\in E$,
 $$K_M|\rho_t(x)|\leq (1-\delta_M)(1\wedge |x|),\ \  \ d\mathbb{P}\times dt\mbox{-a.e.}$$
 $\mathrm{(v)}$ There exist  constants $K_{1}> 0$ and $K_{2}> 0$ such that
 $$\| \eta\|_{\infty}\leq K_{1}\ \ \ \text{and}\ \ \  \Big\| \int _{0}^{T}\theta _{t} dt\Big\|_{\infty}\leq K_{2}.$$
\end{assumption}
In what follows, we define
$$\Lambda_\lambda:=\int_E(1\wedge|x|^2)\lambda(dx),\ \ \ \ \  C_\rho:=K_0\Lambda_\lambda^{\frac{1}{2}}.$$
\begin{remark}\label{re3.1}
\normalfont
 Using H$\ddot{\mathrm{o}}$lder's inequality and Assumption \ref{as3.1}-(iii), for all $U\in \mathbb{J}^{2}(0,T;\mathbb{R})$ and $t\in [0,T]$, we have
\begin{equation*}
\begin{aligned}
 W_2\big(\mathbb{P}_{U_t^{\rho_t}(\cdot)},\delta_0\big) = \Big(\mathbb{E}\Big[\Big|\int_EU_t(x)\rho_t(x)\lambda(dx)\Big|^2\Big]\Big)^\frac{1}{2}  \leq & \Big(\mathbb{E}\Big[\int_E|U_t(x)|^2\lambda(dx)\int_E|\rho_t(x)|^2\lambda(dx)\Big]\Big)^\frac{1}{2}\\
 \leq &C_\rho\|\|U_t\|_{\mathcal{L}^{2}(\lambda)}\|_{L^2(\Omega)}.
 \end{aligned}
\end{equation*}
Similarly, for all $U,\ \bar{U}\in \mathbb{J}^{2}(0,T;\mathbb{R})$ and $t\in [0,T]$, we obtain
$$W_2\big(\mathbb{P}_{U_t^{\rho_t}(\cdot)},\mathbb{P}_{\bar{U}_t^{\rho_t}(\cdot)}\big)\leq C_\rho\|\|U_t-\bar{U}_t\|_{\mathcal{L}^{2}(\lambda)}\|_{L^2(\Omega)}.$$
\end{remark}

Before proving the existence and the uniqueness of local solutions to equation \eqref{eq3.1}, we first present the following proposition.
\begin{proposition} \label{prop3.1}
Given a triple of processes $(P,Q,H)\in \mathbb{S}^{\infty}(0,T;\mathbb{R})\times\mathbb{H}^{2}_{BMO}(0,T;\mathbb{R}^d) \times \mathbb{J}^{2}_{BMO}(0,T;\mathbb{R})$, assume that the generator $f:=f(\omega,t, z,u):\Omega \times[0,T] \times\mathbb{R}^d \times \mathcal{L}^{2}(\lambda)\rightarrow \mathbb{R}$ is  $\mathcal{P}^0$-measurable for each fixed $(z,u)\in\mathbb{R}^d \times \mathcal{L}^{2}(\lambda)$. $\eta$ and $f$ satisfy the following conditions:\\
$\mathrm{(1)}$  For all $t\in [0,T],\ z\in \mathbb{R}^d,$ and $ u\in \mathcal{L}^{2}(\lambda)$,
\begin{equation*}
\begin{aligned}
-\bar{f}(t,z)-[-u]_{\gamma}
\leq\ &f(t,z,u)
\leq \bar{f}(t,z)+[u]_{\gamma},\ \ d\mathbb{P}\times dt\mbox{-a.e.},
\end{aligned}
\end{equation*}
where $$\bar{f}(t,z):=\theta _{t}+\frac{\gamma}{2}|z|^{2}+\phi(|P_{t}|)+\phi_0\big(\|P_{t}\|_{L^{2}(\Omega)}\big)+\gamma _{0}\|Q_{t}\|_{L^{2}(\Omega)}^{1+\alpha}+\beta _{1}\|\|H_{t}\|_{\mathcal{L}^{2}(\lambda)}\|_{L^{2}(\Omega)}.$$
$\mathrm{(2)}$\
For each $M>0$,\ there exists a constant $C_M\ge 0$  such that, $d\mathbb{P}\times dt\mbox{-a.e.}$,  for all $t\in [0,T],\ z,\ \bar{z}\in \mathbb{R}^d,$ and $ u,\ \bar{u}\in \mathcal{L}^{2}(\lambda)\cap \mathcal{L}^{\infty}(\lambda)$ with $\|u\|_{\mathcal L^\infty(\lambda)}\vee\|\bar u\|_{\mathcal L^\infty(\lambda)}\le M$,
\begin{equation*}
\begin{aligned}
|f(t,z,u)-f(t,\bar{z},\bar{u})|\leq \phi\big( |P_{t}|\vee\|P_{t}\|_{L^{2}(\Omega)}\big)\big(1+|z|+|\bar{z}|+2\|Q_{t}\|_{L^{2}(\Omega)}\big)|z-\bar{z}|+C_M\|u-\bar{u}\|_{\mathcal{L}^{2}(\lambda)}.
\end{aligned}
\end{equation*}
$\mathrm{(3)}$\ For each $M>0$,  for all $t\in [0,T],\ z\in \mathbb{R}^d$ and $ u,\ \bar{u}\in \mathcal{L}^{2}(\lambda)\cap \mathcal{L}^{\infty}(\lambda)$ with $\|u\|_{\mathcal L^\infty(\lambda)}\vee\|\bar u\|_{\mathcal L^\infty(\lambda)}\le M$, there exists a $\mathcal{P}^0\otimes \mathcal{B}(E)$-measurable process $\Pi^{z,u,\bar{u}}$ satisfying
\begin{equation*}
\begin{aligned}
f(t,z,u)-f(t,z,\bar{u})
\leq\ &\int_E\Pi_t^{z,u,\bar{u}}(x)(u(x)-\bar{u}(x))\lambda(dx),\ \ \ d\mathbb{P}\times dt\mbox{-a.e.},
\end{aligned}
\end{equation*}
and $(-1+\delta_M)(1\wedge |x|)\leq\Pi_t^{z,u,\bar{u}}(x)\leq D_M(1\wedge |x|)$. Here, $D_M>0$ and $\delta_M>0$ depend only on $M$.\\
 $\mathrm{(4)}$ There exist constants $K_{1}> 0$ and $K_{2}> 0$ such that
 $$\| \eta\|_{\infty}\leq K_{1}\ \ \ \text{and}\ \ \ \Big\| \int _{0}^{T}\theta _{t} dt\Big\|_{\infty}\leq K_{2}.$$
Then the following BSDEJs
\begin{equation}\label{eq3.2}
\begin{aligned}
Y_t=\eta+\int_t^T f(s,Z_s,U_s)ds-\int_t^T Z_sdB_s-\int_t^T \int_E U_s(x)\tilde{N}(ds,dx),\ t\in [0,T]
\end{aligned}
\end{equation}
admits a unique  solution $(Y,Z,U)\in \mathbb{S}^{\infty}(0,T;\mathbb{R})\times \mathbb{H}^{2}_{BMO}(0,T;\mathbb{R}^d)\times \mathbb{J}^{2}_{BMO}(0,T;\mathbb{R})$. Furthermore, for all $t \in[0, T]$ and stopping time $\tau\in \mathcal{T}_t^T$, we have
\begin{equation}\label{eq3.3}
\begin{aligned}
|Y_{t}|\leq\ &\| \eta\|_{\infty}+\Big\| \int _{0}^{T}\theta _{t} dt\Big\|_{\infty}+(\phi+\phi _{0})\big(\|P\|_{\mathbb{S}^{\infty}(t,T)}\big)(T-t)
 +\gamma _{0}\|Q\|_{\mathbb{H}^{2}_{BMO}(t,T)}^{1+\alpha}(T-t)^{\frac{1-\alpha}{2}}\\
 &+\beta _{1}\|H\|_{\mathbb{J}^{2}_{BMO}(t,T)}(T-t)^{\frac{1}{2}},
\end{aligned}
\end{equation}
and
\begin{equation}\label{eq3.4}
\begin{aligned}
&\mathbb{E}_{\tau}\Big[\int _{\tau} ^{T}|Z_{s}|^{2}ds+\int _{\tau} ^{T}\int _{E}|U_{s}(x)|^{2}\lambda(dx)ds\Big]\\
\leq\ &\frac{1}{\gamma ^{2}}e^{2\gamma \| \eta\|_{\infty}}+\frac{2}{\gamma} e^{2\gamma \| Y\|_{\mathbb{S}^{\infty}(t,T)}}\cdot\bigg[\Big\| \int _{0}^{T}\theta _{t} dt\Big\|_{\infty}+(\phi+\phi _{0})\big(\|P\|_{\mathbb{S}^{\infty}(t,T)}\big)(T-t)\\
&  +\gamma _{0}\|Q\|_{\mathbb{H}^{2}_{BMO}(t,T)}^{1+\alpha}(T-t)^{\frac{1-\alpha}{2}}+\beta _{1}\|H\|_{\mathbb{J}^{2}_{BMO}(t,T)}(T-t)^{\frac{1}{2}}\bigg].
\end{aligned}
\end{equation}
\end{proposition}
\begin{proof}
For all $t\in [0,T]$, we define
 $$A_{t}:= \theta _{t}+\phi\big( |P_{t}|\big)+\phi_0\big(\|P_{t}\|_{L^{2}(\Omega)}\big)+\gamma _{0}\|Q_{t}\|_{L^{2}(\Omega)}^{1+\alpha}+\beta _{1}\|\|H_{t}\|_{\mathcal{L}^{2}(\lambda)}\|_{L^{2}(\Omega)}.$$
Then for all $(t,z,u)\in [0,T]\times \mathbb{R}^d\times \mathcal{L}^2(\lambda),$ the Assumption  (1) can be rewritten as
\begin{equation}\label{eq3.5}
\begin{aligned}
-A_{t}-\frac{\gamma}{2}|z|^{2}-[-u]_{\gamma}\leq f(t,z,u)\leq A_{t}+\frac{\gamma}{2}|z|^{2}+[u]_{\gamma},
\end{aligned}
\end{equation}
and for all $t\in[0,T]$,  we have
\begin{equation}\label{eq3.6}
\begin{aligned}
\int _{t}^{T}A_{s}ds\leq\ &\Big\| \int _{0}^{T}\theta _{t} dt\Big\|_{\infty}+(\phi+\phi _{0})\big(\|P\|_{\mathbb{S}^{\infty}(t,T)}\big)(T-t)
  +\gamma _{0}\|Q\|_{\mathbb{H}^{2}_{BMO}(t,T)}^{1+\alpha}(T-t)^{\frac{1-\alpha}{2}}\\
  &+\beta _{1}\|H\|_{\mathbb{J}^{2}_{BMO}(t,T)}(T-t)^{\frac{1}{2}}.
\end{aligned}
\end{equation}
Indeed, inequality \eqref{eq3.5} clearly holds. We now prove that inequality \eqref{eq3.6} holds. Using H$\ddot{\mathrm{o}}$lder's inequality, we obtain for $t\in[0,T]$,
\begin{equation*}
\begin{aligned}
\int _{t}^{T}\big(\mathbb{E}[|Q_{s}|^{2}]\big)^{\frac{1+\alpha}{2}}ds\leq \Big(\int _{t}^{T}\mathbb{E}[|Q_{s}|^{2}]ds\Big)^{\frac{1+\alpha}{2}}(T-t)^{\frac{1-\alpha}{2}}\leq\|Q\|_{\mathbb{H}^{2}_{BMO}(t,T)}^{1+\alpha}(T-t)^{\frac{1-\alpha}{2}},
\end{aligned}
\end{equation*}
and
\begin{equation*}
\begin{aligned}
\int _{t}^{T}\Big(\mathbb{E}\Big[\int _{E}|H_{s}(x)|^{2}\lambda(dx)\Big]\Big)^{\frac{1}{2}}ds\leq \Big(\int _{t}^{T}\mathbb{E}\Big[\int _{E}|H_{s}(x)|^{2}\lambda(dx)\Big]ds\Big)^{\frac{1}{2}}(T-t)^{\frac{1}{2}}\leq\|H\|_{\mathbb{J}^{2}_{BMO}(t,T)}(T-t)^{\frac{1}{2}},
\end{aligned}
\end{equation*}
therefore, inequality \eqref{eq3.6} holds.

Applying Lemma \ref{leA.1} with $A$ defined as above and $d_s:=\|Q_s\|_{L^2(\Omega)}$, we obtain equation \eqref{eq3.2} admits a unique solution $(Y,Z,U)\in\mathbb{S}^{\infty}(0,T;\mathbb{R})\times \mathbb{H}^{2}_{BMO}(0,T;\mathbb{R}^d)\times \mathbb{J}^{2}_{BMO}(0,T;\mathbb{R}).$ It remains to prove the estimates \eqref{eq3.3} and \eqref{eq3.4}.

First, under inequality \eqref{eq3.5}, following the same calculation as in  \cite[Lemma 3.2]{FT18} yields the following inequality:
$$ |Y_{t}|\leq\ \frac{1}{\gamma}\ln \mathbb{E}_t\Big[\exp\Big(\gamma |\eta|+\gamma\int_t^TA_sds\Big)\Big],\ \ \ \ t\in[0,T].$$
Then using Assumption (4) and estimate \eqref{eq3.6},  we obtain
 \begin{equation*}
\begin{aligned}
|Y_{t}|\leq\ & \| \eta\|_{\infty}+\Big\| \int _{0}^{T}\theta _{t} dt\Big\|_{\infty}+(\phi+\phi _{0})\big(\|P\|_{\mathbb{S}^{\infty}(t,T)}\big)(T-t)
+\gamma _{0}\|Q\|_{\mathbb{H}^{2}_{BMO}(t,T)}^{1+\alpha}(T-t)^{\frac{1-\alpha}{2}}\\
&+\beta _{1}\|H\|_{\mathbb{J}^{2}_{BMO}(t,T)}(T-t)^{\frac{1}{2}},\ \ \ \ \ t\in[0,T].
\end{aligned}
\end{equation*}
We  now prove that estimate \eqref{eq3.4}. Applying It$\mathrm{\hat{o}}$'s formula to $e^{2\gamma Y_{\tau}}$ with $\tau\in \mathcal{T}_t^T$ and using  inequality \eqref{eq3.5}, we get
\begin{equation*}
\begin{aligned}
e^{2\gamma Y_{\tau}}=\ &e^{2\gamma \eta}+\int _{\tau}^{T}\!2\gamma e^{2\gamma Y_{s}} f(s,Z_s,U_s)ds-\int _{\tau}^{T}\!2\gamma e^{2\gamma Y_{s}} Z_{s}dB_{s}-\int _{\tau }^{T}\!2\gamma ^{2}e^{2\gamma Y_{s}}|Z_{s}|^{2}d{s}\\
&-\int _{\tau}^{T}\int _{E}e^{2\gamma Y_{s-}}\big(e^{2\gamma U_{s}(x)}-1\big)\tilde{N}(ds,dx)
-\int _{\tau}^{T}\int _{E}e^{2\gamma Y_{s}}\big(e^{2\gamma U_{s}(x)}-1-2\gamma U_{s}(x)\big)\lambda(dx)ds\\
\leq\ & e^{2\gamma \eta}+\int _{\tau}^{T}2\gamma e^{2\gamma Y_{s}} \big(A_{s}+[U_{s}]_\gamma+\frac{\gamma }{2} |Z_{s}|^{2}\big)ds-\int _{\tau}^{T}2\gamma e^{2\gamma Y_{s}} Z_{s}dB_{s}-\int _{\tau}^{T}2\gamma ^{2}e^{2\gamma Y_{s}}|Z_{s}|^{2}d{s}\\
&-\int _{\tau}^{T}\int _{E}e^{2\gamma Y_{s-}}\big(e^{2\gamma U_{s}(x)}-1\big)\tilde{N}(ds,dx)
-\int _{\tau}^{T}\int _{E}e^{2\gamma Y_{s}}\big(e^{2\gamma U_{s}(x)}-1-2\gamma U_{s}(x)\big)\lambda(dx)ds.\\
\end{aligned}
\end{equation*}
Taking the conditional expectation $\mathbb{E}_{\tau}[\cdot]$ on both sides of the above inequality and using
$$2\gamma[U_{s}]_\gamma-\int _{E}\big(e^{2\gamma U_{s}(x)}-1-2\gamma U_{s}(x)\big)\lambda(dx)=-\int _{E}\big(e^{\gamma U_{s}(x)}-1\big)^{2}\lambda(dx),$$
we obtain
\begin{equation}\label{eq3.7}
\begin{aligned}
&\mathbb{E}_{\tau}\Big[\int_{\tau}^{T}\!\Big(\gamma ^{2}e^{2\gamma Y_{s}}|Z_{s}|^{2}+\int _{E}e^{2\gamma Y_{s}}\big(e^{\gamma U_{s}(x)}-1\big)^{2}\lambda(dx)\Big)ds\Big]
\leq \mathbb{E}_{\tau}\Big[e^{2\gamma \eta}+\int _{\tau}^{T}2\gamma e^{2\gamma Y_{s}} A_{s}ds\Big].
\end{aligned}
\end{equation}
Next, applying It$\mathrm{\hat{o}}$'s formula to $e^{-2\gamma Y_{\tau}}$ with $\tau\in \mathcal{T}_t^T$ and using  inequality \eqref{eq3.5}, we have
\begin{equation*}
\begin{aligned}
e^{-2\gamma Y_{\tau}}\!\leq\ & e^{-2\gamma \eta}+\!\int _{\tau}^{T}\!\!\!2\gamma e^{-2\gamma Y_{s}} \big(A_{s}+[-U_{s}]_\gamma+\frac{\gamma }{2} |Z_{s}|^{2}\big)ds+\!\int _{\tau}^{T}\!\!\!2\gamma e^{-2\gamma Y_{s}} Z_{s}dB_{s}-\!\int _{\tau}^{T}\!\!\!2\gamma ^{2}e^{-2\gamma Y_{s}}|Z_{s}|^{2}d{s}\\
&-\!\int _{\tau}^{T}\!\!\int _{E}\!\!e^{-2\gamma Y_{s-}}\big(e^{-2\gamma U_{s}(x)}-1\big)\tilde{N}(ds,dx)
-\!\int _{\tau}^{T}\!\!\int _{E}\!\!e^{-2\gamma Y_{s}}\big(e^{-2\gamma U_{s}(x)}-1+2\gamma U_{s}(x)\big)\lambda(dx)ds.\\
\end{aligned}
\end{equation*}
Taking the conditional expectation $\mathbb{E}_{\tau}[\cdot]$ on both sides of the above inequality and using
$$2\gamma[-U_{s}]_\gamma-\int _{E}\big(e^{-2\gamma U_{s}(x)}-1+2\gamma U_{s}(x)\big)\lambda(dx)=-\int _{E}\big(1-e^{-\gamma U_{s}(x)}\big)^{2}\lambda(dx),$$
we get
\begin{equation}\label{eq3.8}
\begin{aligned}
\mathbb{E}_{\tau}\Big[\int _{\tau} ^{T}\!\Big(\gamma ^{2}e^{-2\gamma Y_{s}}|Z_{s}|^{2}+\int _{E}\!e^{-2\gamma Y_{s}}\big(1-e^{-\gamma U_{s}(x)}\big)^{2}\lambda(dx)\Big)ds\Big]
 \leq \mathbb{E}_{\tau}\Big[e^{-2\gamma \eta}+\int _{\tau}^{T}\!2\gamma e^{-2\gamma Y_{s}} A_{s}ds\Big].
\end{aligned}
\end{equation}
Noticing the fact that $e^{m}+e^{-m}\geq 2$ and $a(e^{m}-1)^{2}+\frac{1}{a}(1-e^{-m})^{2}\geq 2(e^{m}-1)(1-e^{-m})=2(e^{m}+e^{-m}-2)\geq 2m^{2}$, we obtain\\
\begin{equation*}
\begin{aligned}
&\mathbb{E}_{\tau}\Big[\int _{\tau} ^{T}\Big(2\gamma ^{2}|Z_{s}|^{2}+\int _{E}2\gamma ^{2}|U_{s}(x)|^{2}\lambda(dx)\Big)ds\Big]\\
 \leq\ & \mathbb{E}_{\tau}\Big[e^{2\gamma \eta}+e^{-2\gamma \eta}+\int _{\tau}^{T}2\gamma \big(e^{2\gamma Y_{s}}+e^{-2\gamma Y_{s}}\big) A_{s}ds\Big]\\
\leq\ & 2\bigg\{e^{2\gamma \| \eta\|_{\infty}}+2\gamma e^{2\gamma \| Y\|_{\mathbb{S}^{\infty}(t,T)}}\cdot\bigg[\Big\| \int _{0}^{T}\theta _{t} dt\Big\|_{\infty}+(\phi+\phi _{0})\big(\|P\|_{\mathbb{S}^{\infty}(t,T)}\big)(T-t)\\
&\ \  +\gamma _{0}\|Q\|_{\mathbb{H}^{2}_{BMO}(t,T)}^{1+\alpha}(T-t)^{\frac{1-\alpha}{2}}+\beta _{1}\|H\|_{\mathbb{J}^{2}_{BMO}(t,T)}(T-t)^{\frac{1}{2}}\bigg]\bigg\},
\end{aligned}
\end{equation*}
which implies that
\begin{equation*}
\begin{aligned}
&\mathbb{E}_{\tau}\Big[\int _{\tau} ^{T}|Z_{s}|^{2}ds+\int _{\tau} ^{T}\int _{E}|U_{s}(x)|^{2}\lambda(dx)ds\Big]\\
\leq\ & \frac{1}{\gamma ^{2}}e^{2\gamma \| \eta\|_{\infty}}+\frac{2}{\gamma} e^{2\gamma \| Y\|_{\mathbb{S}^{\infty}(t,T)}}\cdot\bigg[\Big\| \int _{0}^{T}\theta _{t} dt\Big\|_{\infty}+(\phi+\phi _{0})\big(\|P\|_{\mathbb{S}^{\infty}(t,T)}\big)(T-t)\\
&  +\gamma _{0}\|Q\|_{\mathbb{H}^{2}_{BMO}(t,T)}^{1+\alpha}(T-t)^{\frac{1-\alpha}{2}}+\beta _{1}\|H\|_{\mathbb{J}^{2}_{BMO}(t,T)}(T-t)^{\frac{1}{2}}\bigg].
\end{aligned}
\end{equation*}
Hence,  estimates \eqref{eq3.3} and \eqref{eq3.4} hold true. This proof is complete.
\end{proof}

We now establish the local existence and the uniqueness of mean-field BSDEJs \eqref{eq3.1} in a complete closed subset of $\mathbb{S}^{\infty}(T-\varepsilon,T;\mathbb{R})\times \mathbb{H}^{2}_{BMO}(T-\varepsilon,T;\mathbb{R}^d)\times \mathbb{J}^{2}_{BMO}(T-\varepsilon,T;\mathbb{R})$ where $ \varepsilon\in(0, T)$. Furthermore, for positive constants $L_{1},\ L_{2}$ and $L_{3}$, we introduce the following  space
\begin{equation*}
\begin{aligned}
&\mathcal{B}_{\varepsilon}(L_{1}, L_{2}, L_{3}):=\Big\{(Y,Z,U)\in \mathbb{S}^{\infty}(T-\varepsilon,T;\mathbb{R})\times \mathbb{H}^{2}_{BMO}(T-\varepsilon,T;\mathbb{R}^d)\times \mathbb{J}^{2}_{BMO}(T-\varepsilon,T;\mathbb{R}):\\
&\ \ \ \ \ \ \ \ \ \ \ \ \ \ \ \ \ \ \ \ \ \ \ \ \ \ \ \ \  \ \|Y\|_{\mathbb{S}^{\infty}(T-\varepsilon,T)}\leq L_{1},\ \ \|Z\|^2_{\mathbb{H}^{2}_{BMO}(T-\varepsilon,T)}\leq  L_{2},\ \
  \|U\|^2_{\mathbb{J}^{2}_{BMO}(T-\varepsilon,T)}\leq L_{3}\Big\}
\end{aligned}
\end{equation*}
 with the metric induced by the product norm
$$\big\|(Y,Z,U)\big\|^2_{\mathcal{B}_{\varepsilon}(L_{1}, L_{2}, L_{3})}:=\|Y\|^2_{\mathbb{S}^{\infty}(T-\varepsilon,T)}+ \|Z\|^2_{\mathbb{H}^{2}_{BMO}(T-\varepsilon,T)}+
  \|U\|^2_{\mathbb{J}^{2}_{BMO}(T-\varepsilon,T)}.$$
\begin{theorem}\label{th3.1}
Under Assumption \ref{as3.1}, there exists a positive constant $\varepsilon$, depending only on $T$ and the datas in Assumption \ref{as3.1}, such that on the interval $[T - \varepsilon, T]$, mean-field BSDEJs \eqref{eq3.1} admits a unique local  solution $(Y, Z,U) \in \mathcal{B}_{\varepsilon}(L_{1}, L_{2}, L_{3})$ with
\begin{equation}\label{eq3.9}
\begin{aligned}
L_{1}=2(K_{1}+K_{2}),\ \ L_{2}=\frac{2}{\gamma ^2}e^{2\gamma K_1}+\frac{4K_2}{\gamma}e^{2\gamma L_1},\ \ L_{3}=\frac{2}{\gamma ^2}e^{2\gamma K_1}+\frac{4K_2}{\gamma}e^{2\gamma L_1}+8L^2_1.
\end{aligned}
\end{equation}
\end{theorem}
\begin{proof}
In the sequel, for all $U\in \mathbb{J}^{2}(0,T;\mathbb{R})$, we denote $U_t^{\rho}(\cdot):=U_t^{\rho_t}(\cdot),\ t\in[0,T]$.
Given a triple of processes $(P,Q,H)\in \mathbb{S}^{\infty}(0,T;\mathbb{R})\times \mathbb{H}^{2}_{BMO}(0,T;\mathbb{R}^d)\times \mathbb{J}^{2}_{BMO}(0,T;\mathbb{R})$, we consider the following BSDEJs:
\begin{equation}\label{eq3.10}
\begin{aligned}
Y_t=\eta+\int_t^T\!g\big(s,P_s,Z_s,U_{s}^{\rho}(\cdot),\mathbb{P}_{(P_s,Q_s,H_{s}^{\rho}(\cdot))}\big)ds-\int_t^T\!Z_sdB_s
-\int_t^T\!\int_E\!U_s(x)\tilde{N}(ds,dx),  \ t\in[0, T].
\end{aligned}
\end{equation}
We define the generator  $f^{P,Q,H}$:
\begin{equation*}
\begin{aligned}
f^{P,Q,H}(t,z,u)&:=\ g\big(t,P_t,z,u^{\rho_t}(\cdot),\mathbb{P}_{(P_t,Q_t,H_{t}^{\rho}(\cdot))}\big),\ \ (t,z,u)\in [0,T]\times \mathbb{R}^d \times \mathcal{L}^2(\lambda).
\end{aligned}
\end{equation*}
Under Assumption \ref{as3.1}, we can directly  verify that the generator $f^{P,Q,H}$ satisfies conditions (1)-(4) in Proposition \ref{prop3.1}.

Indeed, conditions (1), (2) and (4) are clearly satisfied; we only need to verify that condition (3) holds. For all $t\in [0,T],\ z\in \mathbb{R}^d,\ u,\ \bar{u}\in \mathcal{L}^{2}(\lambda)\cap \mathcal{L}^{\infty}(\lambda)$ satisfying $\|u\|_{\mathcal L^\infty(\lambda)}\vee\|\bar u\|_{\mathcal L^\infty(\lambda)}\le M$, on the one hand, if $u^{\rho_t}(\cdot)\geq \bar{u}^{\rho_t}(\cdot)$, using  Assumption \ref{as3.1}-(iv), we have
\begin{equation*}
\begin{aligned}
f^{P,Q,H}(t,z,u)-f^{P,Q,H}(t,z,\bar{u})\leq\ &\big|g\big(t,P_t,z,u^{\rho_t}(\cdot),\mathbb{P}_{(P_t,Q_t,H_{t}^{\rho}(\cdot))}\big)
-g\big(t,P_t,z,\bar{u}^{\rho_t}(\cdot),\mathbb{P}_{(P_t,Q_t,H_{t}^{\rho}(\cdot))}\big)\big|\\
\leq\ &K_M(u^{\rho_t}(\cdot)- \bar{u}^{\rho_t}(\cdot))=K_M\int_E(u(x)-\bar{u}(x))\rho_t(x)\lambda(dx)
\end{aligned}
\end{equation*}
and $K_M|\rho_t(x)|\leq(1-\delta_M)(1\wedge |x|),$ for some $\delta_M>0$.\\
On the other hand, if $u^{\rho_t}(\cdot)< \bar{u}^{\rho_t}(\cdot)$,  we obtain
\begin{equation*}
\begin{aligned}
f^{P,Q,H}(t,z,u)-f^{P,Q,H}(t,z,\bar{u})\leq\ &\big|g\big(t,P_t,z,u^{\rho_t}(\cdot),\mathbb{P}_{(P_t,Q_t,H_{t}^{\rho}(\cdot))}\big)
-g\big(t,P_t,z,\bar{u}^{\rho_t}(\cdot),\mathbb{P}_{(P_t,Q_t,H_{t}^{\rho}(\cdot))}\big)\big|\\
\leq\ &K_M(\bar{u}^{\rho_t}(\cdot)-u^{\rho_t}(\cdot))=-K_M\int_E(u(x)-\bar{u}(x))\rho_t(x)\lambda(dx).
\end{aligned}
\end{equation*}
Therefore, setting $\Pi_t^{z,u,\bar{u}}(x):=K_M\rho_t(x)\textbf{1}_{\{u^{\rho_t}(\cdot)\geq \bar{u}^{\rho_t}(\cdot)\}}-K_M\rho_t(x)\textbf{1}_{\{u^{\rho_t}(\cdot)< \bar{u}^{\rho_t}(\cdot)\}}$ yields
\begin{equation*}
\begin{aligned}
f^{P,Q,H}(t,z,u)-f^{P,Q,H}(t,z,\bar{u})\leq &\int_E\Pi_t^{z,u,\bar{u}}(x)(u(x)-\bar{u}(x))\lambda(dx)
\end{aligned}
\end{equation*}
with $|\Pi_t^{z,u,\bar{u}}(x)|\leq (1-\delta_M)(1\wedge |x|)$,\ for some $\delta_M>0$.
This implies that condition (3) holds true.

Then according to Proposition \ref{prop3.1},  BSDEJs \eqref{eq3.10} admits a unique solution $(Y,Z,U)\in \mathbb{S}^{\infty}(0,T;\mathbb{R})\times \mathbb{H}^{2}_{BMO}(0,T;\mathbb{R}^d)\times \mathbb{J}^{2}_{BMO}(0,T;\mathbb{R})$.  Consequently, we can define a mapping $\Theta$ on $\mathbb{S}^{\infty}(0,T;\mathbb{R})\times \mathbb{H}^{2}_{BMO}(0,T;\mathbb{R}^d)\times \mathbb{J}^{2}_{BMO}(0,T;\mathbb{R})$ as follows:
$$\Theta(P,Q,H): =(Y,Z,U),\ (P,Q,H)\in \mathbb{S}^{\infty}(0,T;\mathbb{R})\times \mathbb{H}^{2}_{BMO}(0,T;\mathbb{R}^d)\times \mathbb{J}^{2}_{BMO}(0,T;\mathbb{R}).$$
Next, we seek to identify an appropriate subset $\mathcal{B}_{\varepsilon}(L_1, L_2, L_3)$ on which the mapping $\Theta$ is both stable and contractive. For this purpose, we combine estimates \eqref{eq3.3}-\eqref{eq3.4} with Assumption \ref{as3.1}-(v), we obtain for all $t\in[0,T]$,
\begin{equation}\label{eq3.11}
\begin{aligned}
\|Y\|_{\mathbb{S}^{\infty}(t,T)}\leq\ & K_1+K_2+(\phi+\phi _{0})\big(\|P\|_{\mathbb{S}^{\infty}(t,T)}\big)(T-t)
 +\gamma _{0}\|Q\|_{\mathbb{H}^{2}_{BMO}(t,T)}^{1+\alpha}(T-t)^{\frac{1-\alpha}{2}}\\
 &+\beta _{1}\|H\|_{\mathbb{J}^{2}_{BMO}(t,T)}(T-t)^{\frac{1}{2}},
\end{aligned}
\end{equation}
\begin{equation}\label{eq3.12}
\begin{aligned}
\|Z\|^2_{\mathbb{H}^{2}_{BMO}(t,T)} \leq\ &\frac{1}{\gamma ^{2}}e^{2\gamma K_1}+\frac{2}{\gamma} e^{2\gamma \| Y\|_{\mathbb{S}^{\infty}(t,T)}}\cdot\Big[K_2+(\phi+\phi _{0})\big(\|P\|_{\mathbb{S}^{\infty}(t,T)}\big)(T-t)\\
&+\gamma _{0}\|Q\|_{\mathbb{H}^{2}_{BMO}(t,T)}^{1+\alpha}(T-t)^{\frac{1-\alpha}{2}}+\beta _{1}\|H\|_{\mathbb{J}^{2}_{BMO}(t,T)}(T-t)^{\frac{1}{2}}\Big],
\end{aligned}
\end{equation}
and by Lemma \ref{le2.1} we have
\begin{equation}\label{eq3.13}
\begin{aligned}
\|U\|^2_{\mathbb{J}^{2}_{BMO}(t,T)}\leq\ & \|U\|^2_{\mathbb{J}^{2}_{B}(t,T)} +4\|Y\|^2_{\mathbb{S}^{\infty}(t,T)}\\
\leq\ & \frac{1}{\gamma ^{2}}e^{2\gamma K_1}+\frac{2}{\gamma}e^{2\gamma \| Y\|_{\mathbb{S}^{\infty}(t,T)}}\cdot\Big[K_2+(\phi+\phi _{0})\big(\|P\|_{\mathbb{S}^{\infty}(t,T)}\big)(T-t)\\
& +\gamma _{0}\|Q\|_{\mathbb{H}^{2}_{BMO}(t,T)}^{1+\alpha}(T-t)^{\frac{1-\alpha}{2}}+\beta _{1}\|H\|_{\mathbb{J}^{2}_{BMO}(t,T)}(T-t)^{\frac{1}{2}}\Big]+4\|Y\|^2_{\mathbb{S}^{\infty}(t,T)}.
\end{aligned}
\end{equation}
Now, we define
$$L_{1}:=2(K_{1}+K_{2}),\ \ L_{2}:=\frac{2}{\gamma ^2}e^{2\gamma K_1}+\frac{4K_2}{\gamma}e^{2\gamma L_1},\ \ L_{3}:=\frac{2}{\gamma ^2}e^{2\gamma K_1}+\frac{4K_2}{\gamma}e^{2\gamma L_1}+8L^2_1.$$
Let $\varepsilon_{1},\ \varepsilon_{2},$ and $\varepsilon_{3}$ be the unique solutions to the following three equations, respectively:
\begin{equation}\label{eq3.14}
\left\{
\begin{aligned}
&\big[\phi(L_1)+\phi_0(L_1)\big]x+\gamma _0L^{\frac{1+\alpha}{2}}_2x^{\frac{1-\alpha}{2}}+\beta _1L^{\frac{1}{2}}_3x^{\frac{1}{2}}=\frac{L_1}{2},\\
&\big[\phi(L_1)+\phi_0(L_1)\big]y+\gamma _0L^{\frac{1+\alpha}{2}}_2y^{\frac{1-\alpha}{2}}+\beta _1L^{\frac{1}{2}}_3y^{\frac{1}{2}}=\frac{\gamma L_2}{4}e^{-2\gamma L_1},\\
&\big[\phi(L_1)+\phi_0(L_1)\big]z+\gamma _0L^{\frac{1+\alpha}{2}}_2z^{\frac{1-\alpha}{2}}+\beta _1L^{\frac{1}{2}}_3z^{\frac{1}{2}}=\frac{\gamma L_3}{4}e^{-2\gamma L_1}.
\end{aligned}
\right.
\end{equation}
 Hence, combining  \eqref{eq3.11}-\eqref{eq3.14}, we see that if
$$\|P\|_{\mathbb{S}^{\infty}(T-\varepsilon,T)}\leq L_{1},\ \ \ \|Q\|^2_{\mathbb{H}^{2}_{BMO}(T-\varepsilon,T)}\leq  L_{2},\ \ \
  \|H\|^2_{\mathbb{J}^{2}_{BMO}(T-\varepsilon,T)}\leq L_{3},$$
then we obtain
$$\|Y\|_{\mathbb{S}^{\infty}(T-\varepsilon,T)}\leq L_{1},\ \ \|Z\|^2_{\mathbb{H}^{2}_{BMO}(T-\varepsilon,T)}\leq  L_{2},\ \
  \|U\|^2_{\mathbb{J}^{2}_{BMO}(T-\varepsilon,T)}\leq L_{3},\ \ \text{for all} \  \varepsilon \in(0,\varepsilon^*],$$
where $\varepsilon ^*=\min(\varepsilon_{1},\ \varepsilon_{2},\ \varepsilon_{3})> 0$. This implies that,
$$\Theta(P,Q,H)\in \mathcal{B}_{\varepsilon}(L_{1}, L_{2}, L_{3}),\ \ \ \forall (P,Q,H)\in \mathcal{B}_{\varepsilon}(L_{1}, L_{2}, L_{3}).$$
Therefore, the mapping $\Theta$ is stable in $\mathcal{B}_{\varepsilon}(L_{1}, L_{2}, L_{3})$.

We now prove that the mapping $\Theta$ is contractive on $\mathcal{B}_{\varepsilon}(L_{1}, L_{2}, L_{3})$. For a fixed $\varepsilon \in (0, \varepsilon^*)$ and any two triples of processes $(P,Q,H),\ (\bar{P},\bar{Q},\bar{H}) \in \mathcal{B}_{\varepsilon}(L_{1}, L_{2}, L_{3})$, we define
$$(Y,Z,U):=\Theta(P,Q,H),\ \ \ (\bar{Y},\bar{Z},\bar{U}):=\Theta(\bar{P},\bar{Q},\bar{H}).$$
In addition,  for all $t\in[T-\varepsilon,T]$, we define
$$\Delta Y_t:=Y_t-\bar{Y}_t, \ \Delta Z_t:=Z_t-\bar{Z}_t, \ \Delta U_t:=U_t-\bar{U}_t,\ \Delta P_t:=P_t-\bar{P}_t, \ \Delta Q_t:=Q_t-\bar{Q}_t, \ \Delta H_t:=H_t-\bar{H}_t.$$
Then, for all $t\in[T-\varepsilon,T]$, we have
$$\Delta Y_t=\int_t^T(I_{1,s}+I_{2,s}+I_{3,s})ds-\int_t^T \Delta Z_sdB_s-\int_t^T \int_E \Delta U_s(x)\tilde{N}(ds,dx),$$
where
\begin{equation*}
\begin{aligned}
I_{1,s}:=\ &g\big(s,P_s,Z_s,U_{s}^{\rho}(\cdot),\mathbb{P}_{(P_s,Q_s,H_{s}^{\rho}(\cdot))}\big)
-g\big(s,P_s,\bar{Z}_s,U_{s}^{\rho}(\cdot),\mathbb{P}_{(P_s,Q_s,H_{s}^{\rho}(\cdot))}\big),\\
I_{2,s}:=\ &g\big(s,P_s,\bar{Z}_s,U_{s}^{\rho}(\cdot),\mathbb{P}_{(P_s,Q_s,H_{s}^{\rho}(\cdot))}\big)
-g\big(s,P_s,\bar{Z}_s,\bar{U}_{s}^{\rho}(\cdot),\mathbb{P}_{(P_s,Q_s,H_{s}^{\rho}(\cdot))}\big),\\
I_{3,s}:=\ &g\big(s,P_s,\bar{Z}_s,\bar{U}_{s}^{\rho}(\cdot),\mathbb{P}_{(P_s,Q_s,H_{s}^{\rho}(\cdot))}\big)
-g\big(s,\bar{P}_s,\bar{Z}_s,\bar{U}_{s}^{\rho}(\cdot),\mathbb{P}_{(\bar{P}_s,\bar{Q}_s,\bar{H}_{s}^{\rho}(\cdot))}\big).
\end{aligned}
\end{equation*}
Note that the terms $I_{1,s}$ and $ I_{2,s}$ can be expressed as $I_{1,s} =\Lambda _s\Delta Z_s$ and $ I_{2,s} =\Gamma _s\int_E\Delta U_s(x)\rho_s(x)\lambda(dx)$, respectively, where
\begin{equation*}
\begin{aligned}
\Lambda _s=\frac{I_{1,s}}{|\Delta Z_s|^2}{\Delta Z_s^\top\textbf{1}}_{\{\Delta Z_s\neq0\}},\ \ \ \
\Gamma _s=\frac{I_{2,s}}{\int_E\Delta U_s(x)\rho_s(x)\lambda(dx)}{\textbf{1}}_{\{\int_E\Delta U_s(x)\rho_s(x)\lambda(dx)\neq0\}}.
\end{aligned}
\end{equation*}
By Lemma \ref{le2.1} we get $\|U_s\|_{\mathcal{L}^\infty(\lambda)}\vee \|\bar{U}_s\|_{\mathcal{L}^\infty(\lambda)}\leq 2\big(\|Y\|_{\mathbb{S}^{\infty}(T - \varepsilon,T)}\vee \|\bar{Y}\|_{\mathbb{S}^{\infty}(T - \varepsilon,T)}\big)\leq 2L_1,\ \forall s\in[T - \varepsilon,T]$. Hence, Assumption 3.1-(ii) and Assumption 3.1-(iv) with $M=2L_1$ imply
$$|\Lambda _s|\leq \phi(L_1)\big(1+|Z_s|+|\bar{Z_s}|+2\|Q_{s}\|_{L^{2}(\Omega)}\big)\ \ \text{and}\ \
|\Gamma _s\rho_s(x)|\leq (1-\delta_{2L_1})(1\wedge |x|).$$
 Since  $Z\cdot B$ and  $ \bar{Z}\cdot B$ are BMO-martingales, the above estimates imply that  $\Lambda\textbf{1}_{[T-\varepsilon,T]} \cdot B+ \Gamma \rho\textbf{1}_{[T-\varepsilon,T]}\cdot\tilde{N}$ is also a BMO-martingale and its BMO norm is bounded uniformly  for all $(P,Q,H),\ (\bar{P},\bar{Q},\bar{H}) \in \mathcal{B}_{\varepsilon}(L_{1}, L_{2}, L_{3})$.   Hence the constants in Lemma \ref{le2.3} can be chosen independently of these  process. Furthermore, by Lemma  \ref{le2.2}, the Dol$\mathrm{\acute{e}}$ans-Dade exponential
$$(\mathcal{E}_t)_{t\in[0,T]}:=\mathcal{E}\big(\Lambda \textbf{1}_{[T-\varepsilon,T]}\cdot B+ \Gamma \rho\textbf{1}_{[T-\varepsilon,T]}\cdot\tilde{N}\big)$$
is a  uniformly integrable martingale. Consequently, $d\mathbb{Q} = \mathcal{E}_T d\mathbb{P}$ defines a probability measure equivalent to $\mathbb{P}$. Under $\mathbb{Q}$, the processes
\begin{equation}\label{eq3.15}
\begin{aligned}
dB^{\mathbb{Q}}_s:=dB_s-\Lambda _s\textbf{1}_{[T-\varepsilon,T]}(s)ds,\ \ \tilde{N}^{\mathbb{Q}}(ds,dx):=\tilde{N}(ds,dx)-\Gamma _s\rho_s(x)\textbf{1}_{[T-\varepsilon,T]}(s)\lambda(dx)ds
\end{aligned}
\end{equation}
are respectively a Brownian motion and a compensated random measure with compensator $\big(1+\Gamma _s\rho_s(x)\textbf{1}_{[T-\varepsilon,T]}(s)\big)\lambda(dx)ds$. Then we have
\begin{equation}\label{eq3.16}
\begin{aligned}
\Delta Y_t=\int_t^TI_{3,s}ds-\int_t^T \Delta Z_sdB^{\mathbb{Q}}_s-\int_t^T \int_E \Delta U_s(x)\tilde{N}^{\mathbb{Q}}(ds,dx),\ \ \ t\in[T-\varepsilon,T].
\end{aligned}
\end{equation}
Applying It$\mathrm{\hat{o}}$'s formula to $|\Delta Y_t|^{2}$, we obtain
\begin{equation*}
\begin{aligned}
|\Delta Y_t|^{2}=&\int_t^T2\Delta Y_sI_{3,s}ds-\int_t^T 2\Delta Y_s\Delta Z_sdB^{\mathbb{Q}}_s-\int_t^T |\Delta Z_s|^2ds\\
&-\int_t^T \int_E\big(|\Delta Y_{s}+\Delta U_s(x)|^{2}-|\Delta Y_{s}|^{2}-2\Delta Y_s\Delta U_s(x)\big)(1+\Gamma _s\rho_s(x))\lambda(dx)ds\\
&-\int_t^T \int_E \big(|\Delta Y_{s-}+\Delta U_s(x)|^{2}-|\Delta Y_{s-}|^{2}\big)\tilde{N}^{\mathbb{Q}}(ds,dx).
\end{aligned}
\end{equation*}
Taking conditional expectations $\mathbb{E}^{\mathbb{Q}}_{\tau}[\cdot]$ on both sides and applying Young's inequality, we have for all $\tau \in \mathcal{T}_t^T$ with $t\in[T-\varepsilon,T],$
\begin{equation}\label{eq3.17}
\begin{aligned}
&|\Delta Y_\tau|^{2}+\mathbb{E}^{\mathbb{Q}}_{\tau}\Big[\int_{\tau}^T |\Delta Z_s|^2ds\Big]+\mathbb{E}^{\mathbb{Q}}_{\tau}\Big[\int_{\tau}^T \int_E|\Delta U_s(x)|^{2}(1+\Gamma _s\rho_s(x))\lambda(dx)ds\Big]\\
=\ &2\mathbb{E}^{\mathbb{Q}}_{\tau}\Big[\int_{\tau}^T\Delta Y_sI_{3,s}ds\Big]\\
\leq\ &2\mathbb{E}^{\mathbb{Q}}_{\tau}\Big[\|\Delta Y\|_{\mathbb{S}^{\infty}(t,T)}\int_{\tau}^T|I_{3,s}|ds\Big]
 \leq \frac{1}{2}\|\Delta Y\|^2_{\mathbb{S}^{\infty}(t,T)}+2\mathbb{E}^{\mathbb{Q}}_{\tau}\Big[\Big(\int_{\tau}^T|I_{3,s}|ds\Big)^2\Big].
\end{aligned}
\end{equation}
By the Assumption \ref{as3.1}-(ii) and Remark \ref{re3.1}, we have
\begin{equation*}
\begin{aligned}
|I_{3,s}|\leq\ & \phi\left( |P_s|\vee|\bar{P}_s|\vee W_2(\mathbb{P}_{P_s},\delta_0)\vee  W_2(\mathbb{P}_{\bar{P}_s},\delta_0)\right)\\
&\times\Big[ \big(1+2|\bar{Z}_s|+W_2(\mathbb{P}_{Q_s},\delta_0)+  W_2(\mathbb{P}_{\bar{Q}_s},\delta_0)\big)\cdot\big(|\Delta P_s|+ W_2(\mathbb{P}_{P_s},\mathbb{P}_{\bar{P}_s})\big)\\
&+\big(1+2|\bar{Z}_s|^\alpha+W_2(\mathbb{P}_{Q_s},\delta_0)^\alpha+W_2(\mathbb{P}_{\bar{Q}_s},\delta_0)^\alpha\big)W_2(\mathbb{P}_{Q_s},\mathbb{P}_{\bar{Q}_s})
+W_2\big(\mathbb{P}_{H_s^{\rho}(\cdot)},\mathbb{P}_{\bar{H}_s^{\rho}(\cdot)}\big)\Big]\\
\leq\ & \phi( L_1)\Big[ \big(1+2|\bar{Z}_s|+\|Q_{s}\|_{L^{2}(\Omega)}+\|\bar{Q}_{s}\|_{L^{2}(\Omega)}\big)\cdot2\|\Delta P\|_{\mathbb{S}^{\infty}(s,T)}\ \ \ \ \ \ \ \ \ \ \ \  \ \ \ \ \ \ \ \ \ \ \\
&+\big(1+2|\bar{Z}_s|^\alpha+\|Q_{s}\|^\alpha_{L^{2}(\Omega)}+\|\bar{Q}_{s}\|^\alpha_{L^{2}(\Omega)}\big)\|\Delta Q_{s}\|_{L^{2}(\Omega)} +C_\rho\|\|\Delta H_{s}\|_{\mathcal{L}^{2}(\lambda)}\|_{L^{2}(\Omega)}\Big],
\end{aligned}
\end{equation*}
then we can rewrite \eqref{eq3.17} as
\begin{equation}\label{eq3.18}
\begin{aligned}
&|\Delta Y_\tau|^{2}+\mathbb{E}^{\mathbb{Q}}_{\tau}\Big[\int_{\tau}^T |\Delta Z_s|^2ds\Big]+\mathbb{E}^{\mathbb{Q}}_{\tau}\Big[\int_{\tau}^T \int_E|\Delta U_s(x)|^{2}(1+\Gamma _s\rho_s(x))\lambda(dx)ds\Big]\\
\leq\ & \frac{1}{2}\|\Delta Y\|^2_{\mathbb{S}^{\infty}(t,T)}
+6\phi( L_1)^2\mathbb{E}^{\mathbb{Q}}_{\tau}\Big[\Big(\int_{\tau}^T\!\!2\big(1+2|\bar{Z}_s|+\|Q_{s}\|_{L^{2}(\Omega)}+\|\bar{Q}_{s}\|_{L^{2}(\Omega)}\big)
\|\Delta P\|_{\mathbb{S}^{\infty}(s,T)}ds\Big)^2\Big]\\
&+6\phi(L_1)^2\mathbb{E}^{\mathbb{Q}}_{\tau}\Big[\Big(\int_{\tau}^T \big(1+2|\bar{Z}_s|^\alpha+\|Q_{s}\|^\alpha_{L^{2}(\Omega)}+\|\bar{Q}_{s}\|^\alpha_{L^{2}(\Omega)}\big)\|\Delta Q_{s}\|_{L^{2}(\Omega)}ds\Big)^2\Big]\\
&+6\phi(L_1)^2C_\rho^2\mathbb{E}^{\mathbb{Q}}_{\tau}\Big[\Big(\int_{\tau}^T \|\|\Delta H_{s}\|_{\mathcal{L}^{2}(\lambda)}\|_{L^{2}(\Omega)}ds\Big)^2\Big]\\
=:\ & \frac{1}{2}\|\Delta Y\|^2_{\mathbb{S}^{\infty}(t,T)}+6\phi( L_1)^2\mathrm{II}+6\phi(L_1)^2\mathrm{III}+6\phi(L_1)^2C_\rho^2\mathrm{IV}.
\end{aligned}
\end{equation}
We now handle the second term in the right hand side of the last inequality, using H$\ddot{\mathrm{o}}$lder's inequality together with Lemma \ref{le2.3}, we obtain that  for $\tau \in \mathcal{T}_t^T$ with $t\in[T-\varepsilon,T]$,
\begin{equation*}
\begin{aligned}
\mathrm{II}\leq\ & 4\mathbb{E}^{\mathbb{Q}}_{\tau}\Big[\int_{\tau}^T\big(1+2|\bar{Z}_s|+\|Q_{s}\|_{L^{2}(\Omega)}+\|\bar{Q}_{s}\|_{L^{2}(\Omega)}\big)^2ds\Big]\|\Delta P\|^2_{\mathbb{S}^{\infty}(t,T)}(T-t)\\
\leq\ & 16\varepsilon \mathbb{E}^{\mathbb{Q}}_{\tau}\Big[\int_{\tau}^T\big(1+4|\bar{Z}_s|^2+\|Q_{s}\|^2_{L^{2}(\Omega)}+\|\bar{Q}_{s}\|^2_{L^{2}(\Omega)}\big)ds\Big]\|\Delta P\|^2_{\mathbb{S}^{\infty}(t,T)}\\
\leq\ & 16\varepsilon\big(T+4\| \bar{Z}\|^{2}_{\mathbb{H}^{2}_{BMO(\mathbb{Q})}(t,T)}+\| Q\|^{2}_{\mathbb{H}^{2}_{BMO}(t,T)}+\| \bar{Q}\|^{2}_{\mathbb{H}^{2}_{BMO}(t,T)}\big)\|\Delta P\|^2_{\mathbb{S}^{\infty}(t,T)}\\
\leq\ &16\varepsilon\big[T+2\big(2c^2_2+1\big)L_2\big]\|\Delta P\|^2_{\mathbb{S}^{\infty}(t,T)}.
\end{aligned}
\end{equation*}
For the third term, using Holder's inequality, we see that for $\tau \in \mathcal{T}_t^T$ with $t\in[T-\varepsilon,T]$,
\begin{equation*}
\begin{aligned}
\mathrm{III}\leq\ & \mathbb{E}^{\mathbb{Q}}_{\tau}\Big[\int_{\tau}^T\!\big(1+2|\bar{Z}_s|^\alpha+\|Q_{s}\|^\alpha_{L^{2}(\Omega)}+\|\bar{Q}_{s}\|^\alpha_{L^{2}(\Omega)}\big)^2ds\cdot \int_{\tau}^T\!\|\Delta Q_{s}\|^2_{L^{2}(\Omega)}ds\Big]\\
\leq\ &4 \|\Delta Q\|^2_{\mathbb{H}^{2}_{BMO}(t,T)}\Big[\varepsilon+ 4\mathbb{E}^{\mathbb{Q}}_{\tau}\Big[\int_{\tau}^T|\bar{Z}_s|^{2\alpha}ds\Big]+\int_{\tau}^T\!\! \big(\|Q_{s}\|^{2\alpha}_{L^{2}(\Omega)}+ \|\bar{Q}_{s}\|^{2\alpha}_{L^{2}(\Omega)}\big)ds\Big].
\end{aligned}
\end{equation*}
Since $\alpha\in [0,1)$, using H$\ddot{\mathrm{o}}$lder's inequality  together with Lemma \ref{le2.3}, we get, for $\tau \in \mathcal{T}_t^T$ with $t\in[T-\varepsilon,T]$,
\begin{equation*}
\begin{aligned}
\mathbb{E}^{\mathbb{Q}}_{\tau}\Big[\int_{\tau}^T|\bar{Z}_s|^{2\alpha}ds\Big]\leq \Big(\mathbb{E}^{\mathbb{Q}}_{\tau}\Big[\int_{\tau}^T |\bar{Z}_{s}|^2 ds\Big]\Big)^\alpha(T-t)^{1-\alpha}
\leq \|\bar{Z}\|^{2\alpha}_{\mathbb{H}^{2}_{BMO(\mathbb{Q})}(t,T)}(T-t)^{1-\alpha} \leq c^{2\alpha}_2L_2^{\alpha}\varepsilon ^{1-\alpha},
\end{aligned}
\end{equation*}
and
\begin{equation*}
\begin{aligned}
\int_{\tau}^T \|Q_{s}\|^{2\alpha}_{L^{2}(\Omega)}ds\leq \Big(\int_{\tau}^T \mathbb{E}[|Q_{s}|^2] ds\Big)^\alpha(T-t)^{1-\alpha}
\leq \| Q\|^{2\alpha}_{\mathbb{H}^{2}_{BMO}(t,T)}(T-t)^{1-\alpha} \leq L_2^{\alpha}\varepsilon ^{1-\alpha}.
\end{aligned}
\end{equation*}
Therefore combining the above estimates we obtain
\begin{equation*}
\begin{aligned}
\mathrm{III}\leq\ & 4 \varepsilon^{1-\alpha}\big[T^\alpha+2\big(2c^{2\alpha}_2+1\big)L_2^{\alpha}\big]\|\Delta Q\|^2_{\mathbb{H}^{2}_{BMO}(t,T)}.
\end{aligned}
\end{equation*}
For the last term, using the following estimate
 $$\int_{\tau}^T \|\|\Delta H_{s}\|_{\mathcal{L}^{2}(\lambda)}\|^2_{L^{2}(\Omega)}ds=\int_{\tau}^T \mathbb{E}\Big[\int _{E}|\Delta H_{s}(x)|^{2}\lambda(dx)\Big] ds\leq  \|\Delta H\|^2_{\mathbb{J}^{2}_{BMO}(t,T)},$$
then we have
\begin{equation*}
\begin{aligned}
\mathrm{IV}\leq \mathbb{E}^{\mathbb{Q}}_{\tau}\Big[\int_{\tau}^T\|\|\Delta H_{s}\|_{\mathcal{L}^{2}(\lambda)}\|^2_{L^{2}(\Omega)}ds\Big](T-t)
\leq \varepsilon \|\Delta H\|^2_{\mathbb{J}^{2}_{BMO}(t,T)}.
\end{aligned}
\end{equation*}
Inserting the results of II-IV into \eqref{eq3.18} yields
\begin{equation}\label{eq3.19}
\begin{aligned}
&|\Delta Y_\tau|^{2}+\mathbb{E}^{\mathbb{Q}}_{\tau}\Big[\int_{\tau}^T |\Delta Z_s|^2ds\Big]+\mathbb{E}^{\mathbb{Q}}_{\tau}\Big[\int_{\tau}^T \int_E|\Delta U_s(x)|^{2}(1+\Gamma _s\rho_s(x))\lambda(dx)ds\Big]\\
\leq\ & \frac{1}{2}\|\Delta Y\|^2_{\mathbb{S}^{\infty}(t,T)}+6\phi(L_1)^2\Big\{16\varepsilon\big[T+2\big(2c^2_2+1\big)L_2\big] \|\Delta P\|^2_{\mathbb{S}^{\infty}(t,T)}\\
&+4\varepsilon^{1-\alpha}\big[T^\alpha+2\big(2c^{2\alpha}_2+1\big)L_2^{\alpha}\big]\|\Delta Q\|^2_{\mathbb{H}^{2}_{BMO}(t,T)}+ \varepsilon C_\rho^2\|\Delta H\|^2_{\mathbb{J}^{2}_{BMO}(t,T)}\Big\}\\
=:\ &\frac{1}{2}\|\Delta Y\|^2_{\mathbb{S}^{\infty}(t,T)}+\Sigma.
\end{aligned}
\end{equation}
Then we obtain
$$\|\Delta Y\|^{2}_{\mathbb{S}^{\infty}(t,T)}\leq \frac{1}{2}\|\Delta Y\|^2_{\mathbb{S}^{\infty}(t,T)}+\Sigma,\ \ \ t\in[T-\varepsilon,T],$$
which implies that $\|\Delta Y\|^{2}_{\mathbb{S}^{\infty}(t,T)}\leq 2\Sigma$ for all $t\in[T-\varepsilon,T]$. Furthermore, for all $t\in[T-\varepsilon,T]$, we have
$$\|\Delta Z\|^2_{\mathbb{H}^{2}_{BMO(\mathbb{Q})}(t,T)}+\|\Delta U\|^2_{\mathbb{J}^{2}_{B(\mathbb{Q})}(t,T)}\leq 4\Sigma.$$
Thanks to Lemma \ref{le2.1} and Lemma \ref{le2.3}, we have, for all $t\in[T-\varepsilon,T]$,
\begin{equation*}
\begin{aligned}
&\|\Delta Y\|^{2}_{\mathbb{S}^{\infty}(t,T)}+\|\Delta Z\|^2_{\mathbb{H}^{2}_{BMO}(t,T)}+\|\Delta U\|^2_{\mathbb{J}^{2}_{BMO}(t,T)}\\
\leq\ &\|\Delta Y\|^{2}_{\mathbb{S}^{\infty}(t,T)}+\frac{1}{c^2_1}\big(\|\Delta Z\|^2_{\mathbb{H}^{2}_{BMO(\mathbb{Q})}(t,T)}+\|\Delta U\|^2_{\mathbb{J}^{2}_{BMO(\mathbb{Q})}(t,T)}\big)\\
\leq\ &\Big(1+\frac{4}{c^2_1}\Big)\|\Delta Y\|^{2}_{\mathbb{S}^{\infty}(t,T)}+\frac{1}{c^2_1}\big(\|\Delta Z\|^2_{\mathbb{H}^{2}_{BMO(\mathbb{Q})}(t,T)}+\|\Delta U\|^2_{\mathbb{J}^{2}_{B(\mathbb{Q})}(t,T)}\big)
\leq  2\Big(1+\frac{6}{c^2_1}\Big) \Sigma.
\end{aligned}
\end{equation*}
In particular, let $t=T-\varepsilon$, we get
\begin{equation}\label{eq3.20}
\begin{aligned}
&\|\Delta Y\|^{2}_{\mathbb{S}^{\infty}(T-\varepsilon,T)}+\|\Delta Z\|^2_{\mathbb{H}^{2}_{BMO}(T-\varepsilon,T)}+\|\Delta U\|^2_{\mathbb{J}^{2}_{BMO}(T-\varepsilon,T)}\\
\leq\ & 12\Big(1+\frac{6}{c^2_1}\Big)\phi(L_1)^2\Big\{16\varepsilon\big[T+2\big(2c^2_2+1\big)L_2\big] \|\Delta P\|^2_{\mathbb{S}^{\infty}(T-\varepsilon,T)}\\
&+4\varepsilon^{1-\alpha}\big[T^\alpha+2\big(2c^{2\alpha}_2+1\big)L_2^{\alpha}\big]\|\Delta Q\|^2_{\mathbb{H}^{2}_{BMO}(T-\varepsilon,T)}+ \varepsilon C_\rho^2\|\Delta H\|^2_{\mathbb{J}^{2}_{BMO}(T-\varepsilon,T)}\Big\}.
\end{aligned}
\end{equation}
If we choose a sufficiently small $\varepsilon> 0$ such that
$$12\Big(1+\frac{6}{c^2_1}\Big)\phi(L_1)^2 \max\Big\{16\varepsilon\big(T+2\big(2c^2_2+1\big)L_2\big) ,\
4\varepsilon^{1-\alpha}\big(T^\alpha+2\big(2c^{2\alpha}_2+1\big)L_2^{\alpha}\big),\  \varepsilon C_\rho^2\Big \}\leq \frac{1}{2},$$
then we have,
$$\big\|(\Delta Y,\Delta Z,\Delta U)\big\|^2_{\mathcal{B}_{\varepsilon}(L_{1}, L_{2}, L_{3})}\leq \frac{1}{2}\big\|(\Delta P,\Delta Q,\Delta H)\big\|^2_{\mathcal{B}_{\varepsilon}(L_{1}, L_{2}, L_{3})}.$$
Therefore, $\Theta$ is a contraction in $\mathcal{B}_{\varepsilon}(L_{1}, L_{2}, L_{3})$. By  the Banach fixed point theorem, mean-field BSDEJs \eqref{eq3.1} admits a unique local  solution $(Y, Z,U)\in\mathcal{B}_{\varepsilon}(L_{1}, L_{2}, L_{3})$ on the interval $[T - \varepsilon, T]$. The proof is complete.
\end{proof}
\begin{example}\label{ex3.1}
\normalfont
Suppose that Assumption \ref{as3.1}-(iii) holds and the terminal  $\eta$ is bounded. For simplicity, let $\rho$ be deterministic and independent of time $t$. For every $(t,y,z,u,\mu)\in [0,T]\times \mathbb{R}^{1+d}\times \mathcal{L}^{2}(\lambda)\times \mathcal{P}_{2}(\mathbb{R}^{1+d+1})$, we consider the following generator:
$$g\big(t,y,z,u^{\rho}(\cdot),\mu\big)=|y|^2+\frac{\gamma}{2}|z|^2+A\frac{u^{\rho}(\cdot)|u^{\rho}(\cdot)|}{1+2K_0|u^{\rho}(\cdot)|}
+W_2(\mu_1,\delta_0)^3+W_2(\mu_2,\delta_0)^{\frac{3}{2}}+W_2(\mu_3,\delta_0),$$
where $\gamma>0,\ 0<A\leq\min\Big\{\frac{1}{2},\inf\limits_{u\in\mathcal{L}^{2}(\lambda),u^{\rho}(\cdot)\neq0}\frac{[u]_{\gamma}(1+2K_0|u^{\rho}(\cdot)|)}{|u^{\rho}(\cdot)|^2}\Big\}$. Then by Theorem \ref{th3.1}, we obtain the mean-field BSDEJs with terminal condition $\eta$ and generator $g$ has a unique local solution. To this  end, we need to prove the generator $g$ satisfies Assumption \ref{as3.1}-(i),(ii),(iv).

Step 1. We prove that $\inf\limits_{u\in \mathcal{L}^{2}(\lambda),u^{\rho}(\cdot)\neq 0}\frac{[u]_{\gamma}(1+2K_0|u^{\rho}(\cdot)|)}{|u^{\rho}(\cdot)|^2}>0$. Define $w(x):=1\wedge |x|,\  M:=\int_E w^2(x)\lambda(dx)<\infty.$ For $u\in \mathcal{L}^2(\lambda)$, denote $h:=|u^\rho(\cdot)| =\left|\int_E u(x)\rho(x)\lambda(dx)\right|>0.$
 We first use the elementary inequality
$$\frac{e^{\gamma r}-\gamma r-1}{\gamma}\ge\frac{\gamma r^2}{2(1+\gamma |r|)},\ \ \  r\in\mathbb R .$$
Hence, we get $[u]_\gamma\ge\frac{\gamma}{2}\int_E \frac{u^2(x)}{1+\gamma |u(x)|}\lambda(dx).$\\
On the other hand, by $|\rho(x)|\le K_0w(x)$, we have $h\le K_0\int_E |u(x)|w(x)\lambda(dx).$ Define $B:=\int_E |u(x)|w(x)\lambda(dx),$  then we get $h\le K_0B$. By the Cauchy-Schwarz inequality and $0\le w\le1$, we obtain
$$B^2\le\int_E\frac{u^2(x)}{1+\gamma |u(x)|}\lambda(dx)\cdot\int_E w^2(x)(1+\gamma |u(x)|)\lambda(dx)\le\int_E\frac{u^2(x)}{1+\gamma |u(x)|}\lambda(dx)(M+\gamma B).$$
It follows that $\int_E\frac{u^2(x)}{1+\gamma |u(x)|}\lambda(dx)\ge\frac{B^2}{M+\gamma B}.$
Since the function $s\mapsto \frac{s^2}{M+\gamma s}$ is increasing on $[0,+\infty)$, and $B\ge \frac{h}{K_0}$, we get
$$[u]_\gamma\ge\frac{\gamma}{2K_0}\frac{h^2}{K_0M+\gamma h}\ \ \text{and}\ \ \frac{[u]_\gamma(1+2K_0h)}{h^2}\ge\frac{\gamma}{2K_0}\frac{1+2K_0h}{K_0M+\gamma h }.$$
Define $f(h):=\frac{\gamma}{2K_0}\frac{1+2K_0h}{K_0M+\gamma h},\  h\ge0.$
Since the function $f$ is monotone on $[0,+\infty)$, it follows that $f(h)\ge\min\left\{f(0),\lim_{h\to+\infty}f(h)\right\}=\min\{\frac{\gamma}{2K_0^2M}, 1\}.$
Therefore,
$$\inf_{\substack{u\in \mathcal{L}^2(\lambda), u^\rho(\cdot)\neq0}}\frac{[u]_\gamma(1+2K_0|u^\rho(\cdot)|)}{|u^\rho(\cdot)|^2}\ge\min\Big\{\frac{\gamma}{2K_0^2M},1\Big\}>0.$$

Step 2. We now verify that Assumption 3.1-(i),(ii), and (iv) hold. Assumption \ref{as3.1}-(ii) clearly holds, and we only need to verify Assumption  \ref{as3.1}-(i),(iv). First, we verify Assumption  \ref{as3.1}-(i) holds.    From the upper and lower bounds of $A$ and the definition of $[u]_\gamma$, it follows that for all $u\in \mathcal{L}^{2}(\lambda)$,
\begin{equation}\label{3.2}
\begin{aligned}
-[-u]_{\gamma}\leq A\frac{u^{\rho}(\cdot)|u^{\rho}(\cdot)|}{1+2K_0|u^{\rho}(\cdot)|}\leq [u]_{\gamma}.
 \end{aligned}
\end{equation}
Indeed, we set $h:=u^{\rho}(\cdot)$. For $h>0$, the upper bound of $A$ and  the definition of $[-u]_\gamma$ yield
$$-[-u]_{\gamma}\leq 0< A\frac{h|h|}{1+2K_0|h|}\leq \frac{[u]_{\gamma}(1+2K_0h)}{h^2}\cdot \frac{h^2}{1+2K_0h}=[u]_{\gamma}.$$
On the other hand, since $|(-u)^{\rho}(\cdot)|=|-u^{\rho}(\cdot)|=|h|$, for $h < 0$, we have
$$-[-u]_{\gamma}=\frac{[-u]_{\gamma}(1-2K_0h)}{h^2}\cdot \frac{-h^2}{1-2K_0h} \leq A\frac{h|h|}{1+2K_0|h|}<0\leq [u]_{\gamma}.$$
Therefore,   we obtain equation \eqref{3.2}.
Furthermore,  this implies that Assumption \ref{as3.1}-(i) holds.

Next, we will prove that Assumption  \ref{as3.1}-(iv) holds. Define $G(h):=A\frac{h|h|}{1+2K_0|h|}$ for $ h\in \mathbb{R}.$ Using $0<A \leq \frac{1}{2}$, we obtain $|G'(h)|\leq \frac{1}{4K_0}$ for all $h\in\mathbb{R}$. In fact, by a direct computation, we have $|G'(h)|=\frac{2A|h|(1+K_0|h|)}{(1+2K_0|h|)^2}.$
Setting $a=K_0|h|\geq 0$,  we have $|G'(h)|=\frac{2A}{K_0}\cdot \frac{a(1+a)}{(1+2a)^2}$. Since the function $\frac{a(1+a)}{(1+2a)^2}$  is non-decreasing  on $[0,+\infty)$ and takes values on $[0, \frac{1}{4}]$, using $0<A \leq \frac{1}{2}$, we conclude $|G'(h)|\leq\frac{2A}{K_0}\cdot \frac{1}{4}\leq \frac{1}{4K_0}.$

From the form of $g$ and $|G'|\leq \frac{1}{4K_0},$ for all $t\in [0,T],\ y\in \mathbb{R},\ z\in \mathbb{R}^d,\ u,\ \bar{u}\in \mathcal{L}^{2}(\lambda)$ and $\mu\in \mathcal{P}_{2}(\mathbb{R}^{1+d+1})$, we have
$$|g(t,y,z,u^\rho(\cdot),\mu)-g(t,y,z,\bar{u}^\rho(\cdot),\mu)|\leq \frac{1}{4K_0}|u^\rho(\cdot)-\bar{u}^\rho(\cdot)|.$$
Together with Assumption \ref{as3.1}-(iii), this shows Assumption \ref{as3.1}-(iv) holds. This completes the proof.
\end{example}

\section{The existence and the uniqueness of the  global solution}
Building on the results established in Section 3, this section is devoted to the global well-posedness of mean-field BSDEJs. In Subsection 4.1,  we consider equation \eqref{eq3.1} with generators of at most quadratic growth in $z$. We first study the case of linear growth in $y$ and $\mu_1$. We then consider subquadratic growth in $\mu_2$ and linear growth in $\mu_3$, with sufficiently small growth coefficients for these two terms. In Subsection 4.2, we consider generators of the form $f(t,z,u,\nu):=f_1(t,z,u)+f_2(t,\nu),$ where $f_2$ is deterministic. Under this additive structure, we establish global existence and uniqueness while allowing quadratic growth in $\nu_2$, without requiring the corresponding growth coefficient to be small.

In the following, we suppose that $\theta: \Omega \times[0,T] \rightarrow \mathbb{R}^{+}$ is $\mathbb{F}$-progressively measurable, $\phi:[0,\infty)\rightarrow[0,\infty)$ is a nondecreasing continuous function, $\beta_0,\ \gamma,\ \gamma_0$ and $K_0$ are positive constants and $\alpha\in [0,1)$.  Let $\eta: \Omega \rightarrow \mathbb{R}$  be $\mathcal{F}_T$-measurable,  $g:=g(\omega,t,y,z,h,\mu):\Omega \times[0,T]\times \mathbb{R}\times\mathbb{R}^d \times \mathbb{R}\times \mathcal{P}_{2}(\mathbb{R}^{1+d+1})\rightarrow \mathbb{R}$  be $\mathcal{P}^0$-measurable for every fixed $(y,z,h,\mu)\in \mathbb{R}\times\mathbb{R}^d \times \mathbb{R}\times \mathcal{P}_{2}(\mathbb{R}^{1+d+1})$, and $ \rho:\Omega \times[0,T]\times E\rightarrow \mathbb{R}$ be $\mathcal{P}^0\otimes\mathcal{B}(E)$-measurable.

\subsection{The case where the generator $g$ has quadratic growth in $z$}
We begin by studying the global well-posedness of equation \eqref{eq3.1} under the assumption that the generator $g$ has linear growth in $y$ and $\mu_1$. The precise assumptions are as follows.
\begin{assumption}\label{as4.1}
\normalfont
$\mathrm{(i)}$\ There exists a constant  $ L\geq 0$  such that,  for all $t\in [0,T],\ y\in \mathbb{R},\ z\in \mathbb{R}^d,\ u\in \mathcal{L}^{2}(\lambda),$ and $ \mu\in \mathcal{P}_{2}(\mathbb{R}^{1+d+1})$,
\begin{equation*}
\begin{aligned}
-[-u]_{\gamma}-f(t,y,z,\mu)
\leq g\big(t,y,z,u^{\rho_t}(\cdot),\mu\big)
\leq [u]_{\gamma}+f(t,y,z,\mu),\ d\mathbb{P}\times dt\mbox{-a.e.},
\end{aligned}
\end{equation*}
where $$f(t,y,z,\mu):=\theta _{t}+L \big(|y|+W_2(\mu _{1},\delta _{0})\big)+\frac{\gamma}{2}|z|^{2}.$$
$\mathrm{(ii)}$\ There exists  a constant $L\geq 0$ such that, $d\mathbb{P}\times dt\mbox{-a.e.}$, for all $t\in [0,T],\ y,\ \bar{y}\in\mathbb{R},\ z,\ \bar{z}\in \mathbb{R}^d,\ u\in \mathcal{L}^{2}(\lambda),$ and $\mu,\ \bar{\mu}\in \mathcal{P}_{2}(\mathbb{R}^{1+d+1})$,
\begin{equation*}
\begin{aligned}
&\big|g\big(t,y,z,u^{\rho_t}(\cdot),\mu\big)-g\big(t,\bar{y},\bar{z},u^{\rho_t}(\cdot),\bar{\mu}\big)\big|\\
\leq\ &
\phi\big(|y|\vee|\bar{y}|\vee W_2(\mu _{1},\delta _{0})\vee W_2(\bar{\mu} _{1},\delta _{0})\big)\big[\big(1+|z|+|\bar{z}|+W_2(\mu _{2},\delta _{0})+ W_2(\bar{\mu} _{2},\delta _{0})\big)|z-\bar{z}|\\
&+W_2(\mu_2 ,\bar{\mu}_2)\big]+L\big(|y-\bar{y}|+W_2(\mu_1 ,\bar{\mu}_1)+W_2(\mu_3 ,\bar{\mu}_3)\big).
\end{aligned}
\end{equation*}
$\mathrm{(iii)}$\ For all $t\in [0,T],$ and $ x\in E$,
$$|\rho_t(x)|\leq K_0(1\wedge |x|),\ \ \ d\mathbb{P}\times dt\mbox{-a.e.}$$
$\mathrm{(iv)}$
For each $M>0$, there exist constants $K_M\geq0$ and  $\delta_M\in(0,1)$ such that, $d\mathbb{P}\times dt\mbox{-a.e.}$, for all $t\in [0,T],\ y\in \mathbb{R},\ z\in \mathbb{R}^d,\ u,\ \bar{u}\in \mathcal{L}^{2}(\lambda)\cap \mathcal{L}^\infty(\lambda)$ with
$\|u\|_{\mathcal{L}^\infty(\lambda)}\vee\|\bar u\|_{\mathcal{L}^\infty(\lambda)}\le M,$ and $ \mu\in \mathcal{P}_{2}(\mathbb{R}^{1+d+1}),$
\begin{equation*}
\begin{aligned}
\big|g\big(t,y,z,u^{\rho_t}(\cdot),\mu\big)-g\big(t,y,z,\bar{u}^{\rho_t}(\cdot),\mu\big)\big|\leq K_M|u^{\rho_t}(\cdot)-\bar{u}^{\rho_t}(\cdot)|.
\end{aligned}
\end{equation*}
Moreover, for all $t\in[0,T]$ and $x\in E$,
 $$K_M|\rho_t(x)|\leq (1-\delta_M)(1\wedge |x|),\ \ \ d\mathbb{P}\times dt\mbox{-a.e.}$$
 $\mathrm{(v)}$\ There exist constants $K_{1}> 0$ and $K_{3}> 0$ such that
 $$\| \eta\|_{\infty}\leq K_{1}\ \ \ \text{and}\ \ \ \Big\| \int _{0}^{T}|\theta _{t}|^2 dt\Big\|_{\infty}\leq K_{3}.$$
 \end{assumption}
 \begin{remark}
 \normalfont
 Using H$\ddot{\mathrm{o}}$lder's inequality, we have
 $$\Big\| \int _{0}^{T}\theta _{t}dt\Big\|_{\infty}\leq \Big\| \sqrt{T}\Big(\int _{0}^{T}|\theta _{t}|^2 dt\Big)^{\frac{1}{2}}\Big\|_{\infty}\leq \sqrt{TK_3}=K_2.$$
 \end{remark}

 \begin{theorem} \label{th4.1}
Under Assumption \ref{as4.1}, mean-field BSDEJs (\ref{eq3.1}) admits a unique global solution $(Y,Z,U)\in \mathbb{S}^{\infty}(0,T;\mathbb{R})\times \mathbb{H}^{2}_{BMO}(0,T;\mathbb{R}^d)\times \mathbb{J}^{2}_{BMO}(0,T;\mathbb{R})$ on the  interval $[0, T]$. Moreover, there exist  positive constants $M_{1},\ M_{2}$ and $ M_{3}$ depending only on $\gamma,\ L,\  K_{1},\  K_{3} $ and $ T$ such that
\begin{equation}\label{4.1}
\begin{aligned}
\|Y\|_{\mathbb{S}^{\infty}(0,T)}\leq M_{1},\ \ \ \|Z\|^2_{\mathbb{H}^{2}_{BMO}(0,T)}\leq  M_{2},\ \  \
  \|U\|^2_{\mathbb{J}^{2}_{BMO}(0,T)}\leq M_{3}.
\end{aligned}
\end{equation}
\end{theorem}
\begin{proof}
We first consider the solvability of  mean-field BSDEJs \eqref{eq3.1} on some interval $[T - \kappa _{\iota}, T]$, where $\kappa _{\iota}$ is a positive constant to be determined later. To this end, we define
$$\bar{C}:=K_1^2+K_3+L^2+2L+2.$$
Furthermore, we let the function $\Theta$ be the unique solution of the following ODE
\begin{equation}\label{eq4.1}
\begin{aligned}
\Theta (t)=\bar{C}+\int _t^T\bar{C}\Theta (s)ds,\ \ \ t\in [0,T].
\end{aligned}
\end{equation}
It is clear that the function $\Theta$ is continuous and decreasing in $t$, and define
$$\iota:=\sup \limits_{t\in [0,T]}\Theta (t)=\Theta (0).$$
Since $\| \eta \|_\infty^2 \le \bar{C} \leq \iota$, Theorem \ref{th3.1} implies that there exists a positive constant $\kappa_\iota$, depending only on $\iota$ and the fixed structural data in Theorem \ref{th3.1}, such that   mean-field BSDEJs \eqref{eq3.1} admits a  local solution $(Y^1, Z^1, U^1)$ on  the interval $[T - \kappa_\iota, T]$. Next, we will show that the following inequality:
$$\|Y^{1}\|_{\mathbb{S}^{\infty}(t,T)}^2\leq \Theta (t),\ \ \ t\in[T - \kappa _{\iota}, T].$$
In fact, for $t\in[T-\kappa _{\iota}, T]$,
\begin{equation}\label{eq4.2}
\begin{aligned}
Y^{1}_t=&\ \eta+\int_t^T\!\! \Big[g\big(s,Y^{1}_s,Z^{1}_s,(U^{1}_s)^{\rho}(\cdot),\mathbb{P}_{(Y^{1}_s,Z^{1}_s,(U^{1}_s)^{\rho}(\cdot))}\big)
-g\big(s,0,Z^{1}_s,(U^{1}_s)^{\rho}(\cdot),\mathbb{P}_{(0,Z^{1}_s,(U^{1}_s)^{\rho}(\cdot))}\big)\\
&+g\big(s,0,Z^{1}_s,(U^{1}_s)^{\rho}(\cdot),\mathbb{P}_{(0,Z^{1}_s,(U^{1}_s)^{\rho}(\cdot))}\big)
-g\big(s,0,0,(U^{1}_s)^{\rho}(\cdot),\mathbb{P}_{(0,Z^{1}_s,(U^{1}_s)^{\rho}(\cdot))}\big)\\
&+g\big(s,0,0,(U^{1}_s)^{\rho}(\cdot),\mathbb{P}_{(0,Z^{1}_s,(U^{1}_s)^{\rho}(\cdot))}\big)
-g\big(s,0,0,0,\mathbb{P}_{(0,Z^{1}_s,(U^{1}_s)^{\rho}(\cdot))}\big)\\
&+g\big(s,0,0,0,\mathbb{P}_{(0,Z^{1}_s,(U^{1}_s)^{\rho}(\cdot))}\big)\Big]ds
-\int_t^T Z^{1}_sdB_s-\int_t^T \int_E U^{1}_s(x)\tilde{N}(ds,dx).
\end{aligned}
\end{equation}
Note that for all $s\in[T-\kappa _{\iota}, T]$, we have
\begin{equation*}
\begin{aligned}
g\big(s,0,Z^{1}_s,(U^{1}_s)^{\rho}(\cdot),\mathbb{P}_{(0,Z^{1}_s,(U^{1}_s)^{\rho}(\cdot))}\big)
-g\big(s,0,0,(U^{1}_s)^{\rho}(\cdot),\mathbb{P}_{(0,Z^{1}_s,(U^{1}_s)^{\rho}(\cdot))}\big)
=\ \tilde{\Lambda} _sZ ^{1}_s,
\end{aligned}
\end{equation*}
\begin{equation*}
\begin{aligned}
g\big(s,0,0,(U^{1}_s)^{\rho}(\cdot),\mathbb{P}_{(0,Z^{1}_s,(U^{1}_s)^{\rho}(\cdot))}\big)
-g\big(s,0,0,0,\mathbb{P}_{(0,Z^{1}_s,(U^{1}_s)^{\rho}(\cdot))}\big)= \tilde{\Gamma}_s\int_EU^{1}_s(x)\rho_s(x)\lambda(dx),
\end{aligned}
\end{equation*}
where
\begin{equation*}
\begin{aligned}
\tilde{\Lambda} _s:=&\frac{g\big(s,0,Z^{1}_s,(U^{1}_s)^{\rho}(\cdot),\mathbb{P}_{(0,Z^{1}_s,(U^{1}_s)^{\rho}(\cdot))}\big)
-g\big(s,0,0,(U^{1}_s)^{\rho}(\cdot),\mathbb{P}_{(0,Z^{1}_s,(U^{1}_s)^{\rho}(\cdot))}\big)}{|Z^{1}_s|^2}
{(Z^{1}_s)^{\top}\textbf{1}}_{\{Z^{1}_s\neq0\}},\\
\tilde{\Gamma} _s:=&\frac{g\big(s,0,0,(U^{1}_s)^{\rho}(\cdot),\mathbb{P}_{(0,Z^{1}_s,(U^{1}_s)^{\rho}(\cdot))}\big)
-g\big(s,0,0,0,\mathbb{P}_{(0,Z^{1}_s,(U^{1}_s)^{\rho}(\cdot))}\big)}{\int_EU^{1}_s(x)\rho_s(x)\lambda(dx)}
{\textbf{1}}_{\{\int_EU^{1}_s(x)\rho_s(x)\lambda(dx)\neq0\}}.
\end{aligned}
\end{equation*}
By Lemma \ref{le2.1} we get $\|U^1_s\|_{\mathcal{L}^\infty(\lambda)}\leq 2\|Y^{1}\|_{\mathbb{S}^{\infty}(T - \kappa _{\iota},T)},\ \forall s\in[T - \kappa _{\iota},T]$. Hence, Assumption \ref{as4.1}-(ii) and Assumption \ref{as4.1}-(iv) with $M=2\|Y^{1}\|_{\mathbb{S}^{\infty}(T - \kappa _{\iota},T)}$ implies
$$|\tilde{\Lambda} _s|\leq \phi(0)(1+|Z^{1}_s|+2\|Z^1_s\|_{L^2(\Omega)}),\ \ \ \
|\tilde{\Gamma} _s\rho_s(x)|\leq \big(1-\delta_{2\|Y^{1}\|_{\mathbb{S}^{\infty}(T - \kappa _{\iota},T)}}\big)(1\wedge |x|).$$
Furthermore, since $Z^{1}\cdot B$ is a BMO-martingale and the process $\tilde{\Gamma} \rho$ is bounded as noted above, it follows that $\tilde{\Lambda}\textbf{1}_{[T-\kappa _{\iota}, T]}\cdot B+ \tilde{\Gamma}\rho\textbf{1}_{[T-\kappa _{\iota}, T]}\cdot\tilde{N}$ is also a BMO-martingale. Then, by Lemma \ref{le2.2}, the Dol$\mathrm{\acute{e}}$ans-Dade exponential
$$(\mathcal{E}_t)_{t\in[0,T]}:=\mathcal{E}\big(\tilde{\Lambda}\textbf{1}_{[T-\kappa _{\iota}, T]}\cdot B+ \tilde{\Gamma}\rho\textbf{1}_{[T-\kappa _{\iota}, T]}\cdot\tilde{N}\big)$$
is a uniformly integrable martingale. Consequently, $d\widetilde{\mathbb{P}} = \mathcal{E}_T d\mathbb{P}$ defines a probability measure equivalent to $\mathbb{P}$. Under $\widetilde{\mathbb{P}}$, the processes
$$dB^1_s:=dB_s-\tilde{\Lambda} _s\textbf{1}_{[T-\kappa _{\iota}, T]}(s)ds,\ \ \tilde{N}^{1}(ds,dx):=\tilde{N}(ds,dx)-\tilde{\Gamma} _s\rho_s(x)\textbf{1}_{[T-\kappa _{\iota}, T]}(s)\lambda(dx)ds$$
are respectively a Brownian motion and a compensated random measure with compensator $\big(1+\tilde{\Gamma} _s\rho_s(x)\textbf{1}_{[T-\kappa _{\iota}, T]}(s)\big)\lambda(dx)ds$. Hence,  equation \eqref{eq4.2} can be rewritten as
\begin{equation*}
\begin{aligned}
Y^{1}_t=&\ \eta+\int_t^T\!\! \Big[g\big(s,Y^{1}_s,Z^{1}_s,(U^{1}_s)^{\rho}(\cdot),\mathbb{P}_{(Y^{1}_s,Z^{1}_s,(U^{1}_s)^{\rho}(\cdot))}\big)
-g\big(s,0,Z^{1}_s,(U^{1}_s)^{\rho}(\cdot),\mathbb{P}_{(0,Z^{1}_s,(U^{1}_s)^{\rho}(\cdot))}\big)\\
&+g\big(s,0,0,0,\mathbb{P}_{(0,Z^{1}_s,(U^{1}_s)^{\rho}(\cdot))}\big)\Big]ds
-\int_t^T Z^{1}_sdB^{1}_s
-\int_t^T \int_E U^{1}_s(x)\tilde{N}^{1}(ds,dx),\ \ \ t\in[T - \kappa _{\iota}, T].
\end{aligned}
\end{equation*}
 Applying It$\mathrm{\hat{o}}$'s formula to $|Y^{1}_t|^2$ and using Assumption \ref{as4.1}-(i),(ii),  for $r \in [T - \kappa_{\iota}, t],$ we obtain
\begin{equation*}
\begin{aligned}
&\tilde{\mathbb{E}}_r\Big[|Y^{1}_t|^2+\int_t^T |Z^{1}_s|^2ds+\int_t^T \int_E|U^{1}_s(x)|^2(1+\tilde{\Gamma} _s\rho_s(x))\lambda(dx)ds\Big]\\
=\ &\tilde{\mathbb{E}}_r\Big[|\eta|^2+2\int_t^T Y^{1}_s \Big(g\big(s,Y^{1}_s,Z^{1}_s,(U^{1}_s)^{\rho}(\cdot),\mathbb{P}_{(Y^{1}_s,Z^{1}_s,(U^{1}_s)^{\rho}(\cdot))}\big)
-g\big(s,0,Z^{1}_s,(U^{1}_s)^{\rho}(\cdot),\mathbb{P}_{(0,Z^{1}_s,(U^{1}_s)^{\rho}(\cdot))}\big)\\
&+g\big(s,0,0,0,\mathbb{P}_{(0,Z^{1}_s,(U^{1}_s)^{\rho}(\cdot))}\big)\Big)ds\Big]\\
\leq\ & \tilde{\mathbb{E}}_r\Big[|\eta|^2+\int_t^T 2|Y^{1}_s|\big[|\theta_s|+L\big(|Y^{1}_s|+\|Y^{1}_s\|_{L^2(\Omega)}\big)\big]ds\Big]\\
\leq\ & \tilde{\mathbb{E}}_r\Big[|\eta|^2+\int_t^T \big[|\theta_s|^2+(L^2+2L+1)|Y^{1}_s|^2+\mathbb{E}[|Y^{1}_s|^2]\big]ds\Big]\\
\leq\ & K^2_1+K_3+\int_t^T(L^2+2L+2)\|Y^{1}\|_{\mathbb{S}^{\infty}(s,T)}^2 ds,
\end{aligned}
\end{equation*}
where the last inequality follows from Assumption  \ref{as4.1}-(v). From the definition of $\bar{C}$, we get
\begin{equation}\label{eq4.3}
\begin{aligned}
\|Y^{1}\|_{\mathbb{S}^{\infty}(t,T)}^2\leq \bar{C}+\int_t^T \bar{C} \|Y^{1}\|_{\mathbb{S}^{\infty}(s,T)}^2ds,\ \ \ \  t\in[T - \kappa _{\iota}, T].
\end{aligned}
\end{equation}
Recall  equation \eqref{eq4.1}
\begin{equation*}
\begin{aligned}
\Theta(t)= \bar{C}+\int_t^T \bar{C} \Theta (s)ds,\ \ \ t\in[0,T].
\end{aligned}
\end{equation*}
By using  the comparison theorem, we have
\begin{equation*}
\begin{aligned}
\|Y^{1}\|_{\mathbb{S}^{\infty}(t,T)}^2\leq \Theta(t),\ \ \ t\in[T - \kappa _{\iota}, T].
\end{aligned}
\end{equation*}
By taking $T - \kappa _{\iota}$ as the terminal time and $Y_{T - \kappa _{\iota}}$
as the terminal condition, we apply Theorem \ref{th3.1} to obtain a local solution  $(Y^2, Z^2,U^2)$ of mean-field BSDEJ \eqref{eq3.1} on  $[T-2\kappa _{\iota},T- \kappa _{\iota}]$. We further define, for $t\in [T-2\kappa _{\iota},T],$
$$(\bar{Y}_t,\bar{Z}_t,\bar{U}_t):=\left(Y_{t}^{1}, Z_{t}^{1},U_t^1\right) {\textbf{1}}_{\left[T-\kappa_{\iota}, T\right]}(t)+\left(Y_{t}^{2}, Z_{t}^{2},U_t^2\right) {\textbf{1}}_{\left[T-2 \kappa_{\iota}, T-\kappa_{\iota}\right)}(t).$$
Therefore, the triple $(\bar{Y},\bar{Z},\bar{U})$ constitutes a solution to  mean-field BSDEJs \eqref{eq3.1} on the interval $[T-2\kappa _{\iota},T]$.
Moreover, repeating the above argument yields
\begin{equation*}
\begin{aligned}
\|\bar{Y}\|_{\mathbb{S}^{\infty}(t,T)}^2\leq \Theta(t),\ \ \ t\in[T - 2\kappa _{\iota}, T].
\end{aligned}
\end{equation*}
By iterating this procedure a finite number of times, we extend the local solution successively backward in time and obtain a solution $(Y, Z,U)$ to mean-field BSDEJ \eqref{eq3.1} on the entire interval $[0, T]$, with $Y$  uniformly bounded by $\iota= \Theta(0)$. In other words, there exists a positive constant $M_1$ depending only on $L,\ K_1,\ K_3$ and $ T$ such that
$$\|Y\|_{\mathbb{S}^{\infty}(0,T)}\leq M_{1}.$$
Next, we show that $Z\cdot B$ and $ U\cdot \tilde{N}$ are BMO-martingales. Applying It$\mathrm{\hat{o}}$'s formula to $e^{2\gamma Y_{\tau}}$ with $\tau\in \mathcal{T}_0^T$, we have\\
\begin{equation*}
\begin{aligned}
e^{2\gamma Y_{\tau}}=\ &\mathbb{E}_\tau \Big[e^{2\gamma \eta}+\int _{\tau}^{T}2\gamma e^{2\gamma Y_{s}}g\big(s,Y_s,Z_s,U_s^{\rho}(\cdot), \mathbb{P}_{(Y_s,Z_s,U_s^{\rho}(\cdot))}\big)ds-\int _{\tau }^{T}2\gamma ^{2}e^{2\gamma Y_{s}}|Z_{s}|^{2}d{s}\\
&
 -\int _{\tau}^{T}\int _{E}e^{2\gamma Y_{s}}\big(e^{2\gamma U_{s}(x)}-1-2\gamma U_{s}(x)\big)\lambda(dx)ds\Big]\\
\leq\ & \mathbb{E}_\tau \Big[e^{2\gamma \eta}+\int _{\tau}^{T}2\gamma e^{2\gamma Y_{s}} \big(|\theta_{s}|+L|Y_{s}|+L\|Y_s\|_{L^2(\Omega)}+[U_{s}]_\gamma\big)ds-\int _{\tau}^{T}\gamma ^{2}e^{2\gamma Y_{s}}|Z_{s}|^{2}d{s}\\
& -\int _{\tau}^{T}\int _{E}e^{2\gamma Y_{s}}\big(e^{2\gamma U_{s}(x)}-1-2\gamma U_{s}(x)\big)\lambda(dx)ds\Big].
\end{aligned}
\end{equation*}
Using the following equation
$$2\gamma[U_{s}]_\gamma-\int _{E}\big(e^{2\gamma U_{s}(x)}-1-2\gamma U_{s}(x)\big)\lambda(dx)=-\int _{E}\big(e^{\gamma U_{s}(x)}-1\big)^{2}\lambda(dx),$$
we obtain
\begin{equation}\label{eq4.4}
\begin{aligned}
&\mathbb{E}_\tau \Big[\int _{\tau}^{T}\Big(\gamma ^{2}e^{2\gamma Y_{s}}|Z_{s}|^{2}+\int _{E}e^{2\gamma Y_{s}}\big(e^{\gamma U_{s}(x)}-1\big)^2\lambda(dx)\Big)ds\Big]\\
\leq\ & \mathbb{E}_\tau \Big[e^{2\gamma \eta}+\!\int _{\tau}^{T}\!\!2\gamma e^{2\gamma Y_{s}} \big(|\theta_{s}|+L|Y_{s}|+L\|Y_s\|_{L^2(\Omega)}\big)ds\Big]
\leq e^{2\gamma K_1}+2\gamma e^{2\gamma M_1}\big(\sqrt{K_3T}+2LM_1T\big).
\end{aligned}
\end{equation}
Applying the same argument to $e^{-2\gamma Y_{\tau}}$ with $\tau\in \mathcal{T}_0^T$, yields the following estimate
\begin{equation}\label{eq4.5}
\begin{aligned}
&\mathbb{E}_\tau \Big[\int _{\tau}^{T}\!\!\!\Big(\gamma ^{2}e^{-2\gamma Y_{s}}|Z_{s}|^{2}+\int _{E}\!\!e^{-2\gamma Y_{s}}\big(1-e^{-\gamma U_{s}(x)}\big)^2\lambda(dx)\Big)ds\Big]
\leq e^{2\gamma K_1}+2\gamma e^{2\gamma M_1}\big(\sqrt{K_3T}+2LM_1T\big).
\end{aligned}
\end{equation}
Noticing the fact that $e^{m}+e^{-m}\geq 2$ and $a(e^{m}-1)^{2}+\frac{1}{a}(1-e^{-m})^{2}\geq 2(e^{m}-1)(1-e^{-m})=2(e^{m}+e^{-m}-2)\geq 2m^{2}$ one obtains
\begin{equation*}
\begin{aligned}
&\mathbb{E}_{\tau}\Big[\int _{\tau} ^{T}\Big(|Z_{s}|^{2}+\int _{E}|U_{s}(x)|^{2}\lambda(dx)\Big)ds\Big]
\leq  \frac{1}{\gamma ^2}e^{2\gamma K_1}+\frac{2}{\gamma } e^{2\gamma M_1}\big(\sqrt{K_3T}+2LM_1T\big).
\end{aligned}
\end{equation*}
Therefore, by Lemma \ref{le2.1}, we have
$$\|Z\|^2_{\mathbb{H}^{2}_{BMO}(0,T)}\leq \frac{1}{\gamma ^2}e^{2\gamma K_1}+\frac{2}{\gamma } e^{2\gamma M_1}\big(\sqrt{K_3T}+2LM_1T\big)=:M_2,$$
and
\begin{equation*}
\begin{aligned}
\|U\|^2_{\mathbb{J}^{2}_{BMO}(0,T)}\leq\ &\|U\|^2_{\mathbb{J}^{2}_{B}(0,T)}+4\|Y\|^2_{\mathbb{S}^{\infty}(0,T)}\\
\leq\ &4M^2_1+\frac{1}{\gamma ^2}e^{2\gamma K_1}+\frac{2}{\gamma } e^{2\gamma M_1}\big(\sqrt{K_3T}+2LM_1T\big)=:M_3.
\end{aligned}
\end{equation*}
Finally, we prove the uniqueness. Let $(Y, Z,U)$ and $(\bar{Y}, \bar{Z},\bar{U})$ be two  solutions of  mean-field BSDEJs \eqref{eq3.1}.
For all $s\in[0,T]$ and $ x\in E$, we denote
$$\Delta Y_s:=Y_s-\bar{Y}_s,\ \ \Delta Z_s:=Z_s-\bar{Z}_s,\ \ \Delta U_s(x):=U_s(x)-\bar{U}_s(x).$$
Then for arbitrary $\varepsilon > 0$ and  $t \in [T -\varepsilon, T]$, we have
\begin{equation*}
\begin{aligned}
\Delta Y_t=\ &\int_t^T\big[g\big(s,Y_s,\bar{Z}_s,\bar{U}_s^{\rho}(\cdot),\mathbb{P}_{(Y_s,Z_s,U_s^{\rho}(\cdot))}\big)
-g\big(s,\bar{Y}_s,\bar{Z}_s,\bar{U}_s^{\rho}(\cdot),\mathbb{P}_{(\bar{Y}_s,\bar{Z}_s,\bar{U}_s^{\rho}(\cdot))}\big)\big]ds\\
&-\int_t^T \Delta Z_sdB^{\mathbb{Q}}_s
-\int_t^T \int_E \Delta U_s(x)\tilde{N}^{\mathbb{Q}}(ds,dx),
\end{aligned}
\end{equation*}
where $dB^{\mathbb{Q}}_s$ and $\tilde{N}^{\mathbb{Q}}(ds,dx)$ are defined as in \eqref{eq3.15} of Theorem \ref{th3.1}, except that $P,\ Q,\ H,\ \bar{P},\ \bar{Q}$ and $\bar{H}$ in the definitions of $\Lambda$ and $\Gamma$ are replaced by  $Y,\ Z,\ U,\ \bar{Y},\ \bar{Z}$ and $\bar{U}$, respectively.
Arguing as in the proof of \eqref{eq3.17}-\eqref{eq3.20}, we obtain, for all $t \in [T -\varepsilon, T]$,
\begin{equation*}
\begin{aligned}
&\|\Delta Y\|^{2}_{\mathbb{S}^{\infty}(t,T)}+\|\Delta Z\|^2_{\mathbb{H}^{2}_{BMO}(t,T)}+\|\Delta U\|^2_{\mathbb{J}^{2}_{BMO}(t,T)}\\
\leq\ & 12\Big(1+\frac{6}{c^2_1}\Big)\varepsilon\Big[4\varepsilon L^2 \|\Delta Y\|^2_{\mathbb{S}^{\infty}(t,T)}+\phi(M_1)^2\|\Delta Z\|^2_{\mathbb{H}^{2}_{BMO}(t,T)}
+ L^2 C_\rho^2\|\Delta U\|^2_{\mathbb{J}^{2}_{BMO}(t,T)}\Big].
\end{aligned}
\end{equation*}
Clearly, if $\varepsilon$ is small enough, such that
$$12\Big(1+\frac{6}{c^2_1}\Big)\varepsilon\max\Big\{4\varepsilon L^2,\ \phi(M_1)^2,\ L^2 C_\rho^2\Big\}<1, $$
then  we can obtain $Y=\bar{Y},\ Z = \bar{Z},\ U = \bar{U}$ on the interval $[T -\varepsilon, T]$. Repeating this argument within a finite steps, the uniqueness can be obtained. The proof is then complete.
\end{proof}
\begin{remark}
\normalfont
Observe that $Y_t$ is uniformly bounded for all $t\in[0,T]$. This property is essential for passing from a local solution to a global one. By inequality \eqref{eq4.3} and Gronwall's inequality, Assumption \ref{as4.1} guarantees that $Y_t$ admits a universal  bound $\iota$ valid for all $t\in[0,T]$.
\end{remark}
\begin{example}\label{ex4.1}
 \normalfont
Suppose that Assumption \ref{as4.1}-(iii) holds and the terminal  $\eta$ is bounded. For simplicity, let $\rho$ be deterministic and independent of time $t$. For every $(t,y,z,u,\mu)\in [0,T]\times \mathbb{R}^{1+d}\times \mathcal{L}^{2}(\lambda)\times \mathcal{P}_{2}(\mathbb{R}^{1+d+1})$, we consider the following generator:
$$g\big(t,y,z,u^{\rho}(\cdot),\mu\big)=|y|+\frac{\gamma}{2}|z|^2+A\frac{u^{\rho}(\cdot)|u^{\rho}(\cdot)|}{1+2K_0|u^{\rho}(\cdot)|}
+W_2(\mu_1,\delta_0),$$
 where $\gamma>0,\ 0<A\leq \min\Big\{\frac{1}{2},\inf\limits_{u\in \mathcal{L}^{2}(\lambda),u^{\rho}(\cdot)\neq 0}\frac{[u]_{\gamma}(1+2K_0|u^{\rho}(\cdot)|)}{|u^{\rho}(\cdot)|^2}\Big\}$. \\
Based on Example \ref{ex3.1}, it is easy to verify that the above coefficient $g$ satisfies Assumption \ref{as4.1} with $\theta=0,\ L=1$. Hence, the mean-field BSDEJs with terminal condition $\eta$ and generator $g$ has a unique solution.
\end{example}

Next, we investigate the well-posedness of global solutions for equation \eqref{eq3.1} under another assumptions.
\begin{assumption}\label{as4.2}
\normalfont
$\mathrm{(i)}$\
For all $t\in [0,T],\ y\in \mathbb{R},\ z\in \mathbb{R}^d,\ u\in \mathcal{L}^{2}(\lambda),$ and $ \mu\in \mathcal{P}_2(\mathbb{R}^{1+d+1})$,
\begin{equation*}
\begin{aligned}
&-\theta _{t}-\frac{\gamma}{2}|z|^{2}-[-u]_{\gamma}-\gamma _{0}W_2(\mu_2,\delta_0)^{1+\alpha}-\beta _{0}W_2(\mu_3,\delta_0)\\
\leq\ & g\big(t,y,z,u^{\rho_t}(\cdot),\mu\big)
\leq  \theta _{t}+\frac{\gamma}{2}|z|^{2}+[u]_{\gamma}+\gamma _{0}W_2(\mu_2,\delta_0)^{1+\alpha}+\beta _{0}W_2(\mu_3,\delta_0),\ d\mathbb{P}\times dt\mbox{-a.e.}
\end{aligned}
\end{equation*}
$\mathrm{(ii)}$\ There exists  a constant $ C\geq 0$ such that, $ d\mathbb{P}\times dt\mbox{-a.e.}$, for all $t\in [0,T],\ y,\ \bar{y}\in \mathbb{R},\ z,\ \bar{z}\in \mathbb{R}^d,\ u\in \mathcal{L}^{2}(\lambda),$ and $ \mu,\ \bar{\mu}\in \mathcal{P}_2(\mathbb{R}^{1+d+1})$,
\begin{equation*}
\begin{aligned}
&\big|g\big(t,y,z,u^{\rho_t}(\cdot),\mu\big)-g\big(t,\bar{y},\bar{z},u^{\rho_t}(\cdot),\bar{\mu}\big)\big|\\
\leq &C\Big[|y-\bar{y}|+\big(1+|z|+|\bar{z}|+W_2(\mu _{2},\delta _{0})+W_2(\bar{\mu} _{2},\delta _{0})\big)|z-\bar{z}|+W_2(\mu_1,\bar{\mu}_1)+W_2(\mu _{3},\bar{\mu} _{3})\\
&+\big(1+|z|^\alpha+|\bar{z}|^\alpha+W_2(\mu _{2},\delta _{0})^{\alpha}+W_2(\bar{\mu} _{2},\delta _{0})^{\alpha}\big)W_2(\mu _{2},\bar{\mu} _{2})\Big].
\end{aligned}
\end{equation*}
$\mathrm{(iii)}$\ For all $t\in[0,T], $ and $ x\in E$,
$$|\rho_t(x)|\leq K_0(1\wedge |x|),\ \ d\mathbb{P}\times dt\mbox{-a.e.}$$
$\mathrm{(iv)}$
For each $M>0$, there exist constants $K_M\geq0$ and  $\delta_M\in(0,1)$ such that, $d\mathbb{P}\times dt\mbox{-a.e.}$, for all $t\in [0,T],\ y\in \mathbb{R},\ z\in \mathbb{R}^d,\ u,\ \bar{u}\in \mathcal{L}^{2}(\lambda)\cap \mathcal{L}^\infty(\lambda)$ with
$\|u\|_{\mathcal{L}^\infty(\lambda)}\vee\|\bar u\|_{\mathcal{L}^\infty(\lambda)}\le M,$ and $ \mu\in \mathcal{P}_{2}(\mathbb{R}^{1+d+1}),$
\begin{equation*}
\begin{aligned}
\big|g\big(t,y,z,u^{\rho_t}(\cdot),\mu\big)-g\big(t,y,z,\bar{u}^{\rho_t}(\cdot),\mu\big)\big|\leq K_M|u^{\rho_t}(\cdot)-\bar{u}^{\rho_t}(\cdot)|.
\end{aligned}
\end{equation*}
Moreover, for all $t\in[0,T]$ and $ x\in E$,
 $$K_M|\rho_t(x)|\leq (1-\delta_M)(1\wedge |x|),\ \ \ d\mathbb{P}\times dt\mbox{-a.e.}$$
 $\mathrm{(v)}$ There exist  constants $K_{1}> 0$ and $K_{2}> 0$ such that
 $$\| \eta\|_{\infty}\leq K_{1}\ \ \ \text{and}\ \ \   \Big\| \int _{0}^{T}\theta _{t} dt\Big\|_{\infty}\leq K_{2}.$$
\end{assumption}
 \begin{theorem} \label{th4.2}
Under Assumption \ref{as4.2}, if the parameters $\gamma_0$ and $ \beta_0$ are sufficiently small,  mean-field BSDEJs (\ref{eq3.1}) has a unique global solution $(Y,Z,U)\in \mathbb{S}^{\infty}(0,T;\mathbb{R})\times \mathbb{H}^{2}_{BMO}(0,T;\mathbb{R}^d)\times \mathbb{J}^{2}_{BMO}(0,T;\mathbb{R})$. Moreover, there exist positive constants $\widetilde{M}_{1},\ \widetilde{M}_{2}$ and $ \widetilde{M}_{3}$ depending only on $ K_{1},\  K_{2},$ and  $\gamma$ such that
\begin{equation}\label{eq4.6}
\begin{aligned}
\|Y\|_{\mathbb{S}^{\infty}(0,T)}\leq \widetilde{M}_{1},\ \ \ \|Z\|^2_{\mathbb{H}^{2}_{BMO}(0,T)}\leq  \widetilde{M}_{2},\ \ \
  \|U\|^2_{\mathbb{J}^{2}_{BMO}(0,T)}\leq \widetilde{M}_{3}.
\end{aligned}
\end{equation}
\end{theorem}
\begin{proof}
We construct a Cauchy sequence $\{(Y^n,Z^n,U^n)\}_{n\geq 1}$ to prove the result. For this purpose, we let $Y^1=0,\ Z^1=0,\ U^1=0$, and for $n\geq 1$,
\begin{equation}\label{eq4.7}
\begin{aligned}
Y^{n+1}_t=\ &\eta+\int_t^Tg\big(s,Y_s^n,Z_s^{n+1},(U_s^{n+1})^{\rho}(\cdot),\mathbb{P}_{(Y_s^n,Z_s^{n},(U_s^{n})^{\rho}(\cdot))}\big)ds-\int_t^T Z_s^{n+1}dB_s\\
&-\int_t^T \int_E U_s^{n+1}(x)\tilde{N}(ds,dx),\ \ \ \ \ t\in [0, T].
\end{aligned}
\end{equation}
Similar to the proof in Theorem \ref{th3.1}, under  Assumptions \ref{as4.2}, the generators of equation \eqref{eq4.7} satisfy the conditions of Proposition \ref{prop3.1} with $\beta_1:=\beta_0C_\rho$.  Therefore, it follows from Proposition \ref{prop3.1}  that  BSDEJs \eqref{eq4.7}  admits, for each $n \in \mathbb{N}$, a unique solution $(Y^{n+1}, Z^{n+1}, U^{n+1})\in\mathbb{S}^{\infty}(0,T;\mathbb{R})\times \mathbb{H}^{2}_{BMO}(0,T;\mathbb{R}^d)\times \mathbb{J}^{2}_{BMO}(0,T;\mathbb{R})$.  Following a derivation process similar to that of  \eqref{eq3.3} and \eqref{eq3.4}, and applying Lemma \ref{le2.1}, we obtain the following estimates, for all $t \in [0, T]$:
\begin{equation}\label{eq4.8}
\begin{aligned}
\|Y^{n+1}\|_{\mathbb{S}^{\infty}(t,T)}\leq\ & K_1+K_2
+\gamma _{0}\|Z^n\|_{\mathbb{H}^{2}_{BMO}(t,T)}^{1+\alpha}T^{\frac{1-\alpha}{2}}+\beta _{1}\|U^n\|_{\mathbb{J}^{2}_{BMO}(t,T)}T^{\frac{1}{2}},
\end{aligned}
\end{equation}
\begin{equation}\label{eq4.9}
\begin{aligned}
\|Z^{n+1}\|^2_{\mathbb{H}^{2}_{BMO}(t,T)} \leq & \frac{1}{\gamma ^{2}}e^{2\gamma K_1}+\frac{2}{\gamma} e^{2\gamma \| Y^{n+1}\|_{\mathbb{S}^{\infty}(t,T)}}\Big[K_2+\gamma _{0}\|Z^n\|_{\mathbb{H}^{2}_{BMO}(t,T)}^{1+\alpha}T^{\frac{1-\alpha}{2}}
+\beta _{1}\|U^n\|_{\mathbb{J}^{2}_{BMO}(t,T)}T^{\frac{1}{2}}\Big]
\end{aligned}
\end{equation}
and
\begin{equation}\label{eq4.10}
\begin{aligned}
\|U^{n+1}\|^2_{\mathbb{J}^{2}_{BMO}(t,T)}\leq\ & \|U^{n+1}\|^2_{\mathbb{J}^{2}_{B}(t,T)} +4\|Y^{n+1}\|^2_{\mathbb{S}^{\infty}(t,T)}\\
\leq & \frac{1}{\gamma ^{2}}e^{2\gamma K_1}+\frac{2}{\gamma} e^{2\gamma \| Y^{n+1}\|_{\mathbb{S}^{\infty}(t,T)}}\Big[K_2+\gamma _{0}\|Z^n\|_{\mathbb{H}^{2}_{BMO}(t,T)}^{1+\alpha}T^{\frac{1-\alpha}{2}}
 +\beta _{1}\|U^n\|_{\mathbb{J}^{2}_{BMO}(t,T)}T^{\frac{1}{2}}\Big]\\
 &+4\|Y^{n+1}\|^2_{\mathbb{S}^{\infty}(t,T)}.
\end{aligned}
\end{equation}
We now prove by induction that for each $n\geq1,$  the following result holds: for all $t\in[0,T]$,
\begin{equation}\label{eq4.11}
\begin{aligned}
\|Y^n\|_{\mathbb{S}^{\infty}(t,T)}\leq\widetilde{M}_{1}, \ \ \|Z^n\|^2_{\mathbb{H}^{2}_{BMO}(t,T)}\leq  \widetilde{M}_{2},\ \ \
  \|U^n\|^2_{\mathbb{J}^{2}_{BMO}(t,T)}\leq \widetilde{M}_{3},
  \end{aligned}
\end{equation}
where $$\widetilde{M}_{1}:=2(K_{1}+K_2),\ \widetilde{M}_{2}:=\frac{2}{\gamma ^2}e^{2\gamma K_1}+\frac{4K_2}{\gamma}e^{2\gamma \widetilde{M}_1},\ \widetilde{M}_{3}:=\widetilde{M}_{2}+8\widetilde{M}_{1}^2.$$
In fact, it is obvious that triple of processes $(Y^1,Z^1,U^1)$ satisfy estimate \eqref{eq4.11}, and assume inequality \eqref{eq4.11} holds for $n$. Then we just need to prove that inequality \eqref{eq4.11} holds for $n+1$. If we  choose sufficiently small parameters $\gamma_0$ and $\beta_0$ such that  the following inequality hold: 
\begin{equation}\label{eq4-12}
\begin{aligned}
\gamma _{0}\widetilde{M}_{2}^{\frac{1+\alpha}{2}}T^{\frac{1-\alpha}{2}}+\beta _{1}\widetilde{M}_{3}^{\frac{1}{2}}T^{\frac{1}{2}}\leq \min\Big\{\frac{\widetilde{M}_{1}}{2},
\frac{\gamma \widetilde{M}_{2}}{4}e^{-2\gamma \widetilde{M}_1}\Big\}.
\end{aligned}
\end{equation}
Then, from \eqref{eq4.8}, we deduce that
$$\|Y^{n+1}\|_{\mathbb{S}^{\infty}(t,T)}\leq\widetilde{M}_{1},\ \ \ t\in[0,T].$$
Substituting this bound into \eqref{eq4.9}-\eqref{eq4.10}, we obtain
$$\|Z^{n+1}\|^2_{\mathbb{H}^{2}_{BMO}(t,T)}\leq  \widetilde{M}_{2},\ \ \ \ \
  \|U^{n+1}\|^2_{\mathbb{J}^{2}_{BMO}(t,T)}\leq \widetilde{M}_{3},\ \ \ t\in[0,T].$$
Hence, the inequality \eqref{eq4.11} also holds for $n+1$. Therefore, by induction,  the inequality \eqref{eq4.11} holds for all $n\geq 1$.

Next, we prove that the sequence $\{(Y^n,Z^n,U^n)\}_{n\geq 1}$ converges on $[0,T]$. For an interval $[s,t]$, we define
$$D_{n,m}(s,t):=\|Y^{n}-Y^m\|^2_{\mathbb{S}^{\infty}(s,t)}+\|Z^{n}-Z^m\|^2_{\mathbb{H}^{2}_{BMO}(s,t)}+\|U^{n}-U^m\|^2_{\mathbb{J}^{2}_{BMO}(s,t)}.$$
Similar to the argument in the proof Theorem \ref{th3.1} that the mapping $\Theta$ is a contraction on $\mathcal{B}_{\varepsilon}(L_{1}, L_{2}, L_{3})$, for $\kappa:=t-s<1$, we obtain
\begin{equation}\label{eq4.12a}
\begin{aligned}
D_{n+1,m+1}(s,t)\le C\|Y^{n+1}_t-Y^{m+1}_t\|_{\infty}^2+C(\kappa+\kappa^{1-\alpha})D_{n,m}(s,t),
  \end{aligned}
\end{equation}
where  $C$ independent of $n,\ m,\ s$ and $t$. We can choose a sufficiently small $\kappa>0$ such that $C(\kappa+\kappa^{1-\alpha})<1$. By inequality \eqref{eq4.12a}, we get the sequence $\{(Y^n,Z^n,U^n)\}_{n\geq 1}$ is a Cauchy sequence on interval $[T-\kappa,T]$.  Furthermore, using inequality \eqref{eq4.12a}, we have
$$\limsup_{n,m\rightarrow\infty}D_{n,m}(T-\kappa,T)\le C(\kappa+\kappa^{1-\alpha})\limsup_{n,m\rightarrow\infty}D_{n,m}(T-\kappa,T)<\infty,$$
which implies that $\lim_{n,m\rightarrow\infty}D_{n,m}(T-\kappa,T)=0.$ Repeating this argument backwards  a finite number of times on $[0,T]$, we conclude that the sequence $\{(Y^n,Z^n,U^n)\}_{n\geq 1}$ converges in $\mathbb{S}^{\infty}(0,T;\mathbb{R})\times \mathbb{H}^{2}_{BMO}(0,T;\mathbb{R}^d)\times \mathbb{J}^{2}_{BMO}(0,T;\mathbb{R})$. Let $(Y,Z,U)$ be its limit, it can be easily confirmed that   $(Y,Z,U)$ is a solution of equation \eqref{eq3.1}. Moreover, $(Y,Z,U)$ satisfies  estimate \eqref{eq4.11}.

 Finally, let $(Y,Z,U)$ and $ (\bar{Y},\bar{Z},\bar{U})$ be two solutions of equation \eqref{eq3.1}. Apply the same estimate \eqref{eq4.12a} to their difference, with $C$ depending on the common BMO norm bounds. We choose a sufficiently small $\kappa>0$  such that $C(\kappa+\kappa^{1-\alpha})<1$. Since the terminal values coincide, we obtain the two solutions  coincide on the last interval $[T-\kappa,T]$.  Repeating the argument backwards proves uniqueness on $[0,T]$. This completes the proof.
\end{proof}
\begin{example}\label{ex4.2}
 \normalfont
Suppose that Assumption \ref{as4.2}-(iii) holds with $C_\rho=1$ and the terminal  $\|\eta\|_{\infty}\le 1$. For simplicity, let $\rho$ be deterministic and independent of time $t,\ T=1$ and $\gamma=1$. For every $(t,y,z,u,\mu)\in [0,1]\times \mathbb{R}^{1+d}\times \mathcal{L}^{2}(\lambda)\times \mathcal{P}_{2}(\mathbb{R}^{1+d+1})$, we consider the following generator:
$$g\big(t,y,z,u^{\rho}(\cdot),\mu\big)=\frac{1}{2}|z|^2+A\frac{u^{\rho}(\cdot)|u^{\rho}(\cdot)|}{1+2K_0|u^{\rho}(\cdot)|}
+10^{-3} W_2(\mu_2,\delta_0)^{\frac{4}{3}}+10^{-3} W_2(\mu_3,\delta_0),$$
 where $ 0<A\leq \min\Big\{\frac{1}{2},\inf\limits_{u\in \mathcal{L}^{2}(\lambda),u^{\rho}(\cdot)\neq 0}\frac{[u]_{1}(1+2K_0|u^{\rho}(\cdot)|)}{|u^{\rho}(\cdot)|^2}\Big\}$. \\
Based on Example \ref{ex3.1}, it is easy to verify that the above coefficient $g$ satisfies Assumption \ref{as4.2} with $\theta=0,\ \alpha=\frac{1}{3}$ and $K_1=K_2=1$. Moreover, from the proof of Theorem \ref{th4.2}, we know that $\widetilde{M}_1=4,\ \widetilde{M}_2=2e^2+4e^8,\ \widetilde{M}_3=128+2e^2+4e^8$ in \eqref{eq4.11} and the small condition \eqref{eq4-12} is satisfied. 
Hence, the mean-field BSDEJs with terminal condition $\eta$ and generator $g$ has a unique solution.
\end{example}

\subsection{A special case.}
Motivated by Proposition \ref{prop3.1}, we next investigate the solvability of a class of mean-field BSDEJs whose generator admits the additive decomposition $f(t,z,u,\nu):= f_1(t,z,u) + f_2(t,\nu)$, for $ (t,z,u,\nu)\in[0,T]\times \mathbb{R}^d\times \mathcal{L}^{2}(\lambda)\times\mathcal{P}_2(\mathbb{R}^{1+d}\times \mathcal{L}^{2}(\lambda))$. Here  $f_1$ has quadratic growth in $z$ and $f_2$ has quadratic growth in $\nu_2$. More precisely, we consider the following mean-field BSDEJ:
\begin{equation}\label{eq4.13}
\begin{aligned}
Y_t=\ &\eta+\int_t^T \!\big[f_1(s,Z_s,U_s)+f_2\big(s, \mathbb{P}_{(Y_s,Z_s,U_s)}\big)\big]ds-\int_t^T\!\! Z_sdB_s
-\int_t^T\!\!\int_E\! U_s(x)\tilde{N}(ds,dx),\ t\in [0, T].
\end{aligned}
\end{equation}
Here,  $f_1:\Omega\times [0,T]\times\mathbb{R}^d\times \mathcal{L}^2(\lambda)\rightarrow \mathbb{R}$ is $\mathcal{P}^0$-measurable for every fixed $(z,u)\in \mathbb{R}^d\times \mathcal{L}^{2}(\lambda)$ and $ f_2: [0,T]\times  \mathcal{P}_2(\mathbb{R}^{1+d}\times \mathcal{L}^2(\lambda))\rightarrow \mathbb{R}$ is  Borel measurable. $f_1,\ f_2$ and terminal value $\eta$ satisfy the following assumption.

\begin{assumption}\label{as4.3}
\normalfont
$\mathrm{(i)}$\
There exists a constant $C\geq 0$ such that, $d\mathbb{P}\times dt\mbox{-a.e.}$, for all $t\in[0,T],\ z\in \mathbb{R}^d,\ u\in \mathcal{L}^{2}(\lambda), $ and $ \nu\in \mathcal{P}_2(\mathbb{R}^{1+d}\times \mathcal{L}^2(\lambda))$,
\begin{equation*}
\begin{aligned}
&-\theta _{t}-\frac{\gamma}{2}|z|^{2}-[-u]_{\gamma}\leq  f_1(t,z,u)\leq \theta _{t}+\frac{\gamma}{2}|z|^{2}+[u]_{\gamma},\\
&|f_2(t,\nu)|\leq C\big(1+W_2(\nu_1,\delta_0)+W_2(\nu_2,\delta_0)^2+W_2(\nu_3,\delta_0)\big).
\end{aligned}
\end{equation*}
$\mathrm{(ii)}$\ There exists a constant $C\ge 0$ such that, $d\mathbb{P}\times dt\mbox{-a.e.}$, for all $t\in [0,T],\ z,\ \bar{z}\in \mathbb{R}^d,\ u\in \mathcal{L}^{2}(\lambda)$ and $ \nu,\ \bar{\nu}\in\mathcal{P}_2(\mathbb{R}^{1+d}\times \mathcal{L}^2(\lambda))$ with $\nu(\mathbb{R}\times \cdot\times \cdot)=\bar{\nu}(\mathbb{R}\times \cdot\times \cdot)$,
\begin{equation*}
\begin{aligned}
&|f_1(t,z,u)-f_1(t,\bar{z},u)|\leq C(1+|z|+|\bar{z}|)|z-\bar{z}|,\\
&|f_2(t,\nu)-f_2(t,\bar{\nu})|\leq CW_2(\nu_1,\bar{\nu}_1).
\end{aligned}
\end{equation*}
Moreover, for each $M>0$,\ there exists a constant $C_M\ge 0$ such that,   for all $t\in [0,T],\ z\in \mathbb{R}^d,$ and $ u,\ \bar{u}\in \mathcal{L}^{2}(\lambda)\cap \mathcal{L}^{\infty}(\lambda)$ with $\|u\|_{\mathcal L^\infty(\lambda)}\vee\|\bar u\|_{\mathcal L^\infty(\lambda)}\le M$,
\begin{equation*}
\begin{aligned}
|f_1(t,z,u)-f_1(t,z,\bar{u})|\leq\ &C_M\|u-\bar{u}\|_{\mathcal{L}^{2}(\lambda)},\ d\mathbb{P}\times dt\mbox{-a.e.}
\end{aligned}
\end{equation*}
$\mathrm{(iii)}$
For each $M>0$,\ $t\in [0,T],\ z\in \mathbb{R}^d$ and $ u,\ \bar{u}\in \mathcal{L}^{2}(\lambda)\cap \mathcal{L}^{\infty}(\lambda)$ with $\|u\|_{\mathcal L^\infty(\lambda)}\vee\|\bar u\|_{\mathcal L^\infty(\lambda)}\le M$, there exists a $\mathcal{P}^0\otimes \mathcal{B}(E)$-measurable process $\Pi^{z,u,\bar{u}}$ satisfying
\begin{equation*}
\begin{aligned}
f_1(t,z,u)-f_1(t,z,\bar{u})\leq \int_E\Pi_t^{z,u,\bar{u}}(x)(u(x)-\bar{u}(x))\lambda(dx),\ \  \ d\mathbb{P}\times dt\mbox{-a.e.},
\end{aligned}
\end{equation*}
and $(-1+\delta_M)(1\wedge |x|)\leq\Pi_t^{z,u,\bar{u}}(x)\leq D_M(1\wedge |x|)$. Here, $D_M>0$ and $\delta_M>0$ depend only on $M$.\\
 $\mathrm{(iv)}$ There exist constants $K_{1}> 0$ and $K_{2}> 0$ such that
 $$\| \eta\|_{\infty}\leq K_{1}\ \ \ \text{and}\ \ \  \Big\| \int _{0}^{T}\theta _{t} dt\Big\|_{\infty}\leq K_{2}.$$
 \end{assumption}
\begin{theorem}\label{th4.3}
Under Assumption \ref{as4.3},   mean-field BSDEJs \eqref{eq4.13} admits a unique global solution $(Y,Z,U)\in \mathbb{S}^{\infty}(0,T;\mathbb{R})\times \mathbb{H}^{2}_{BMO}(0,T;\mathbb{R}^d)\times \mathbb{J}^{2}_{BMO}(0,T;\mathbb{R})$. Moreover, there exist positive constants $\widehat{M}_{1},\ \widehat{M}_{2}$ and $ \widehat{M}_{3}$ depending only on $ K_{1},\  K_{2},\ \gamma,\ C$ and $ T$ such that
\begin{equation}\label{eq4.14}
\begin{aligned}
\|Y\|_{\mathbb{S}^{\infty}(0,T)}\leq \widehat{M}_{1},\ \ \ \|Z\|^2_{\mathbb{H}^{2}_{BMO}(0,T)}\leq  \widehat{M}_{2},\ \ \
  \|U\|^2_{\mathbb{J}^{2}_{BMO}(0,T)}\leq \widehat{M}_{3}.
\end{aligned}
\end{equation}
\end{theorem}
\begin{proof}
First, Proposition \ref{prop3.1} implies that the following equation:
\begin{equation}\label{eq4.15}
\begin{aligned}
\bar{Y}_t=\ &\eta+\int_t^T f_1(s,Z_s,U_s)ds-\int_t^T Z_sdB_s
-\int_t^T \int_E U_s(x)\tilde{N}(ds,dx),\  t\in [0, T]
\end{aligned}
\end{equation}
admits a unique solution $(\bar{Y},Z,U)\in \mathbb{S}^{\infty}(0,T;\mathbb{R})\times \mathbb{H}^{2}_{BMO}(0,T;\mathbb{R}^d)\times \mathbb{J}^{2}_{BMO}(0,T;\mathbb{R})$ such that
\begin{equation*}
\begin{aligned}
\|\bar{Y}\|_{\mathbb{S}^{\infty}(0,T)}\leq \breve{M}_{1},\ \ \ \|Z\|^2_{\mathbb{H}^{2}_{BMO}(0,T)}\leq  \breve{M}_{2},\ \ \
  \|U\|^2_{\mathbb{J}^{2}_{BMO}(0,T)}\leq \breve{M}_{3},
\end{aligned}
\end{equation*}
where $\breve{M}_{1}:=K_1+K_2,\ \breve{M}_{2}:=\gamma^{-2}e^{2\gamma K_1}+2\gamma^{-1}K_2e^{2\gamma \breve{M}_{1}},\ \breve{M}_{3}:=\breve{M}_{2}+4\breve{M}^2_{1}$. Now, we consider the following equation:
\begin{equation}\label{eq4.16}
\begin{aligned}
 Y_t=\bar{Y}_t+\int_t^T f_2\big(s,\mathbb{P}_{(Y_s,Z_s,U_s)}\big)ds,\ \ \ t\in[0,T].
 \end{aligned}
\end{equation}
To obtain the existence of the solution, we now apply Picard iteration. Let $Y_{t}^{0}:= 0,\ t \in[0, T]$, and for $n \geq 0$, we define
\begin{equation}\label{eq4.17}
\begin{aligned}
Y_{t}^{n+1}:=\bar{Y}_{t}+\int_{t}^{T} f_2\big(s,\mathbb{P}_{(Y^n_s,Z_s,U_s)}\big)d s,\ \ \ t \in[0, T].
\end{aligned}
\end{equation}
By induction, we conclude that $Y^{n+1} \in \mathbb{S}^{\infty}(0,T;\mathbb{R}),\ n \geq 0$. Indeed, for $n=0$, by Assumption \ref{as4.3}-(i) we know that for any $t\in[0,T]$,
\begin{equation*}
\begin{aligned}
|Y^1_t|\leq\ &|\bar{Y}_t|+\int_t^T\big|f_2\big(s,\mathbb{P}_{(0,Z_s,U_s)}\big)\big|ds\\
\leq\ &\|\bar{Y}\|_{\mathbb{S}^{\infty}(0,T)}+\int_t^T C\Big(1+\mathbb{E}[|Z_s|^2]+\mathbb{E}\Big[\int_E|U_s(x)|^2\lambda(dx)\Big]^{\frac{1}{2}}\Big)ds\\
\leq\ &\breve{M}_{1}+CT+C\breve{M}_2+C\breve{M}^{\frac{1}{2}}_3T^{\frac{1}{2}}<+\infty,
\end{aligned}
\end{equation*}
this implies that $Y^1\in \mathbb{S}^{\infty}(0,T;\mathbb{R})$. Next, we assume $Y^n\in\mathbb{S}^{\infty}(0,T;\mathbb{R})$ and then prove $Y^{n+1}\in\mathbb{S}^{\infty}(0,T;\mathbb{R}).$ From Assumption \ref{as4.3}-(i), for any $t\in[0,T]$, we have
\begin{equation*}
\begin{aligned}
|Y^{n+1}_t|\leq\ &|\bar{Y}_t|+\int_t^T\big|f_2\big(s,\mathbb{P}_{(Y^n_s,Z_s,U_s)}\big)\big|ds\\
\leq\ &\breve{M}_{1}+CT+CT\|Y^n\|_{\mathbb{S}^{\infty}(0,T)}+C\breve{M}_2+C\breve{M}^{\frac{1}{2}}_3T^{\frac{1}{2}}<+\infty,
\end{aligned}
\end{equation*}
so $Y^{n+1}\in\mathbb{S}^{\infty}(0,T;\mathbb{R}).$ Therefore, by induction we know that $Y^{n+1} \in \mathbb{S}^{\infty}(0,T;\mathbb{R}),\ n \geq 0$.

Let us put $\Delta Y^{n}:=Y^{n}-Y^{n-1}$. Then we have
\begin{equation}\label{eq4.18}
\begin{aligned}
\Delta Y_{t}^{n+1}=\int_{t}^{T} \big[f_2\big(s,\mathbb{P}_{(Y^n_s,Z_s,U_s)}\big)-f_2\big(s,\mathbb{P}_{(Y^{n-1}_s,Z_s,U_s)}\big)\big]d s,\ \ \ t\in[0,T] .
\end{aligned}
\end{equation}
For the sake of convenience in notation, we denote
$$\Delta f_2^{n}(s):=f_2\big(s,\mathbb{P}_{(Y^n_s,Z_s,U_s)}\big)-f_2\big(s,\mathbb{P}_{(Y^{n-1}_s,Z_s,U_s)}\big),\ \ \ s \in[0, T].$$
Assumption \ref{as4.3}-(ii) implies that
\begin{equation}\label{eq4.19}
\begin{aligned}
\left|\Delta f_2^{n}(s)\right| \leq C\big(\mathbb{E}\big[|\Delta Y_{s}^{n}|^{2}\big]\big)^{\frac{1}{2}}\leq C\|\Delta Y^{n}\|_{\mathbb{S}^{\infty}(s,T)},\ \ \ s \in[0, T].
\end{aligned}
\end{equation}
Hence, from equation \eqref{eq4.18} and inequality \eqref{eq4.19} we get
\begin{equation*}
\begin{aligned}
\left|\Delta Y_{t}^{n+1}\right| \leq \int_t^T\left|\Delta f_2^{n}(s)\right| d s \leq C\int_t^T\|\Delta Y^{n}\|_{\mathbb{S}^{\infty}(s,T)} d s,\ \ \ t \in[0, T],
\end{aligned}
\end{equation*}
so we have
\begin{equation}\label{eq4.20}
\begin{aligned}
\|\Delta Y^{n+1}\|_{\mathbb{S}^{\infty}(t,T)}\leq  C\int_t^T\|\Delta Y^{n}\|_{\mathbb{S}^{\infty}(s,T)} d s,\ \ \ t\in[0,T].
\end{aligned}
\end{equation}
Applying \eqref{eq4.20} repeatedly we get for all $n \geq 1$,
$$
\begin{aligned}
\|\Delta Y^{n+1}\|_{\mathbb{S}^{\infty}(t,T)} & \leq\ C^n \int_{t}^{T} \int_{t_{1}}^{T} \cdots \int_{t_{n-1}}^{T} \|\Delta Y^{1}\|_{\mathbb{S}^{\infty}(t_n,T)} d t_{n} \cdots d t_{2} d t_{1} \\
& \leq\ \frac{C^{n}(T-t)^n}{n!} \|\Delta Y^{1}\|_{\mathbb{S}^{\infty}(t,T)}, \ \ \ t \in[0, T].
\end{aligned}
$$
Define $\widetilde{C}:=\|\Delta Y^{1}\|_{\mathbb{S}^{\infty}(0,T)}$. Thus, for $m>n$, we have
$$\left\|Y^{m}-Y^{n}\right\|_{\mathbb{S}^{\infty}(0,T)} \leq \sum_{k=n+1}^{m}\|\Delta Y^{k}\|_{\mathbb{S}^{\infty}(0,T)} \leq \widetilde{C} \sum_{k \geq n} \frac{C^{k}T^k}{k!} \rightarrow 0, \quad \text { as } n \rightarrow+\infty.$$
Therefore, we conclude that there exists a process $Y \in \mathbb{S}^{\infty}(0,T; \mathbb{R})$ such that
$$\|Y^{n}-Y\|_{\mathbb{S}^{\infty}(0,T)} \rightarrow 0, \ \ \  \text { as } n \rightarrow+\infty.$$
By taking the limit as $n \rightarrow+\infty$ in \eqref{eq4.17}, it can be easily verified that $Y$ satisfies \eqref{eq4.16}.

In conclusion, from equations \eqref{eq4.15} and \eqref{eq4.16} we get $(Y, Z,U) \in \mathbb{S}^{\infty}(0,T;\mathbb{R})\times \mathbb{H}^{2}_{BMO}(0,T;\mathbb{R}^d)\times \mathbb{J}^{2}_{BMO}(0,T;\mathbb{R})$ is a solution of  mean-field BSDEJs \eqref{eq4.13}. Moreover, by Assumptions \ref{as4.3}-(i) and equation \eqref{eq4.16}, we obtain
$$\|Y\|_{\mathbb{S}^{\infty}(t,T)}\le\breve{M}_{1}+CT+C\breve{M}_2+C\breve{M}^{\frac{1}{2}}_3T^{\frac{1}{2}}+\int_t^TC\|Y\|_{\mathbb{S}^{\infty}(s,T)}ds,\ \ \ t\in[0,T].$$
 Using Gronwall's inequality, we have $\|Y\|_{\mathbb{S}^{\infty}(t,T)}\le \big(\breve{M}_{1}+CT+C\breve{M}_2+C\breve{M}^{\frac{1}{2}}_3T^{\frac{1}{2}}\big)e^{C(T-t)}$ for $t\in[0,T]$. Hence, estimate \eqref{eq4.14} holds.

Next, we will prove uniqueness. Assume that $(Y^1,Z^1,U^1)$ and $(Y^2,Z^2,U^2)$ are two solutions of  equation \eqref{eq4.13}. We denote
$$\bar{Y}^i_t:=Y^i_t-\int_t^T f_2\big(s,\mathbb{P}_{(Y^i_s,Z^i_s,U^i_s)}\big)ds,\ \ t\in[0,T],\ \ i=1,2.$$
 Then, the triple of processes $(\bar{Y}^i,Z^i,U^i)_{i=1,2}$ satisfies the following equation:
 \begin{equation*}
\begin{aligned}
\bar{Y}^i_t=\ &\eta+\int_t^T f_1(s,Z^i_s,U^i_s)ds-\int_t^T Z^i_sdB_s
-\int_t^T \int_E U^i_s(x)\tilde{N}(ds,dx),\  \ t\in [0, T].
\end{aligned}
\end{equation*}
Finally, note that the above equation admits a unique solution, which entails that $\bar{Y}^1=\bar{Y}^2,\ Z^1=Z^2$ and $U^1=U^2$ on  $[0,T]$. By equation \eqref{eq4.16} and Assumption \ref{as4.3}-(ii), for $t\in[0,T]$, we obtain
$$|Y^1_t-Y^2_t|\le \int_t^T|f_2(s,\mathbb{P}_{(Y^1_s,Z_s,U_s)})-f_2(s,\mathbb{P}_{(Y^2_s,Z_s,U_s)})|ds\le C\int_t^T\|Y^1-Y^2\|_{\mathbb{S}^{\infty}(s,T)}ds.$$
By Gronwall's inequality, we obtain $Y^1=Y^2$ on $[0,T]$. This completes the proof.
\end{proof}
\begin{example}\label{ex4.3}
\normalfont
For all $t\in[0,T],\ z\in\mathbb{R}^d,\ u\in \mathcal{L}^2(\lambda)$ and $\nu\in \mathcal{P}_2(\mathbb{R}^{1+d}\times \mathcal{L}^2(\lambda))$, we consider the following functions:
\begin{equation*}
\begin{aligned}
 f_1(t,z,u)=\ &\frac{\gamma}{2}|z|^2+\int_E(1\wedge|x|)\frac{\gamma |u(x)|^2}{2(1+\gamma|u(x)|)}\lambda(dx),\\
f_2(t,\nu)=\ &W_2(\nu_1,\delta_0)+W_2(\nu_2,\delta_0)^2+W_2(\nu_3,\delta_0).
\end{aligned}
\end{equation*}
 If let terminal condition $\eta$ is bounded, then the mean-field  BSDEJs with terminal condition $\eta$ and generator $f:=f_1+f_2$ has a unique solution. By Theorem \ref{th4.3}, we know that we only verify that Assumption \ref{as4.3}(i)-(iii) hold. Indeed,  it suffices to verify that $f_1$ satisfies the requirements. We define $q_\gamma(r):=\frac{\gamma |r|^2}{2(1+\gamma|r|)}$. By the elementary inequality, for $x\in E$, we obtain
$$0\le q_\gamma(u(x))\le \frac{e^{\gamma u(x)}-\gamma u(x)-1}{\gamma}.$$
Therefore, $0\le \int_E(1\wedge|x|)q_\gamma(u(x))\lambda(dx)\le [u]_\gamma$. This implies that $-\frac{\gamma}{2}|z|^2-[-u]_\gamma\le0\le f_1(t,z,u)\le \frac{\gamma}{2}|z|^2+[u]_\gamma$. Moreover, by direct calculation, we have
$$q'_\gamma(r) = \frac{\operatorname{sgn}(r)}{2} \Big( 1 - \frac{1}{(1 + \gamma |r|)^2} \Big), \ \ \  |q'_\gamma(r)| \le \frac{1}{2}.$$
It's obvious that $f_1$ is Lipschitz in $u$ and
$$f_1(t,z,u)-f_1(t,z,\bar{u})=\int_E \Pi_t(x)(u(x)-\bar{u}(x))\lambda(dx),$$
where $\Pi_t(x):=(1\wedge|x|)\int_0^1q'_\gamma(\bar{u}(x)+\sigma(u(x)-\bar{u}(x)))d\sigma$ and $-\frac{1}{2}(1\wedge|x|)\le\Pi_t(x)\le \frac{1}{2}(1\wedge|x|)$. Hence, $f_1$ satisfies Assumption \ref{as4.3}(i)-(iii).
\end{example}

\section{Comparison Theorem}
In this section, we give a comparison theorem for the solutions of mean-field BSDEJs under quadratic-exponential growth. For simplicity, we will consider the following equation:
\begin{equation}\label{eq5.1}
\begin{aligned}
Y_t=\ &\eta+\int_t^T g\big(s,Y_s,Z_s,U_s^{\rho}(\cdot), \mathbb{P}_{Y_s}\big)ds-\int_t^T Z_sdB_s
-\int_t^T \int_E U_s(x)\tilde{N}(ds,dx),\  t\in [0,T].
\end{aligned}
\end{equation}
We now state the comparison theorem. For notational convenience, a mean-field BSDEJ of the form \eqref{eq5.1} will be referred to as the equation associated with the pair $(\eta,g)$, where $g$ denotes the generator and $\eta$ the terminal condition.
\begin{theorem}\label{th5.1}
Let $g^i(s,y,z,u^{\rho_s}(\cdot),\mu),\ i=1,2,$ satisfy Assumption \ref{as4.1}-(i),(iii),(v). Moreover, we suppose

(1) One of the both generators satisfies Assumption \ref{as4.1}-(ii),(iv).

(2) One of the both generators satisfies: for all $s\in  [0,T],\ (y,z,u)\in \mathbb{R}^{1+d}\times \mathcal{L}^2(\lambda),\ \theta_1,\ \theta_2\in L^2_{\mathcal{F}}(\Omega;\mathbb{R})$, there exists a positive constant $L$ such that
$$g^i\big(s,y,z,u^{\rho_s}(\cdot),\mathbb{P}_{\theta_1}\big)-g^i\big(s,y,z,u^{\rho_s}(\cdot),\mathbb{P}_{\theta_2}\big)\leq L\|(\theta_1-\theta_2)^+\|_{L^2(\Omega)}.$$
Let $\eta^1,\ \eta^2$ satisfy Assumption \ref{as4.1}-(v) and denote by $(Y^1, Z^1,U^1)$ and $(Y^2, Z^2,U^2)$ the global  solutions of  mean-field BSDEJs \eqref{eq5.1} with parameters $(\eta^1, g^1)$ and $(\eta^2, g^2)$, respectively.
Then, if $\eta^1\leq \eta^2,\ \mathbb{P}\mbox{-a.s.}$, and $g^1\leq g^2,\ d\mathbb{P}\times ds\mbox{-a.e.}$,  it holds that also
$$Y^1_s\leq Y^2_s,\ \ \ s \in [0, T],\ \ \ \mathbb{P}\mbox{-a.s.}$$
\end{theorem}
\begin{proof}
Without loss of generality, we assume that conditions (1) and (2) hold for $g^1$. For all $s \in [t, T]$ and $x\in E$, we denote
$$\delta Y_s:=Y^1_s-Y^2_s,\ \ \delta Z_s:=Z^1_s-Z^2_s,\ \ \delta U_s(x):=U^1_s(x)-U^2_s(x),\ \ \delta \eta:=\eta^1-\eta^2,$$
$$\delta g_s:= g^1\big(s,Y^2_s,Z^2_s,(U^2_s)^{\rho}(\cdot), \mathbb{P}_{Y^2_s}\big)-g^2\big(s,Y^2_s,Z^2_s,(U^2_s)^{\rho}(\cdot), \mathbb{P}_{Y^2_s}\big).$$
Then we have
\begin{equation}\label{eq5.2}
\begin{aligned}
\delta Y_t=\ &\delta \eta+\int_t^T \big[g^1\big(s,Y^1_s,Z^1_s,(U^1_s)^{\rho}(\cdot), \mathbb{P}_{Y^1_s}\big)-g^2\big(s,Y^2_s,Z^2_s,(U^2_s)^{\rho}(\cdot), \mathbb{P}_{Y^2_s}\big)\big]ds\\
& -\int_t^T \delta Z_sdB_s-\int_t^T \int_E \delta U_s(x)\tilde{N}(ds,dx),\ \ \ t\in[0,T].
\end{aligned}
\end{equation}
For the generator of  mean-field BSDEJs \eqref{eq5.2}, we note that for all $s \in [t, T]$,
\begin{equation*}
\begin{aligned}
&g^1\big(s,Y^1_s,Z^1_s,(U^1_s)^{\rho}(\cdot), \mathbb{P}_{Y^1_s}\big)-g^2\big(s,Y^2_s,Z^2_s,(U^2_s)^{\rho}(\cdot), \mathbb{P}_{Y^2_s}\big)\\
\le\ &\delta _y g^1_s \delta Y_s+\delta _z g^1_s \delta Z_s+\delta _u g^1_s\int_E \delta U_s(x)\rho_s(x)\lambda(dx)+\delta _\mu g^1_s \|(\delta Y_s)^+\|_{L^2(\Omega)}+\delta g_s,
\end{aligned}
\end{equation*}
where
\begin{equation*}
\begin{aligned}
\delta _y g^1_s=\ &\frac{g^1\left(s,Y^1_s,Z^1_s,(U^1_s)^{\rho}(\cdot), \mathbb{P}_{Y^1_s}\right)-g^1\left(s,Y^2_s,Z^1_s,(U^1_s)^{\rho}(\cdot), \mathbb{P}_{Y^1_s}\right)}{\delta Y_s}{\textbf{1}}_{\{\delta Y_s\neq0\}},\\
\end{aligned}
\end{equation*}
\begin{equation*}
\begin{aligned}
\delta _z g^1_s=\ &\frac{g^1\left(s,Y^2_s,Z^1_s,(U^1_s)^{\rho}(\cdot), \mathbb{P}_{Y^1_s}\right)-g^1\left(s,Y^2_s,Z^2_s,(U^1_s)^{\rho}(\cdot), \mathbb{P}_{Y^1_s}\right)}{|\delta Z_s|^2}{(\delta Z_s)^\top\textbf{1}}_{\{\delta Z_s\neq0\}},\\
\delta _u g^1_s=\ &\frac{g^1\left(s,Y^2_s,Z^2_s,(U^1_s)^{\rho}(\cdot), \mathbb{P}_{Y^1_s}\right)-g^1\left(s,Y^2_s,Z^2_s,(U^2_s)^{\rho}(\cdot), \mathbb{P}_{Y^1_s}\right)}{\int_E\delta U_s(x)\rho_s(x)\lambda(dx)}{\textbf{1}}_{\{\int_E\delta U_s(x)\rho_s(x)\lambda(dx)\neq0\}},\\
\delta _\mu g^1_s=\ &\frac{g^1\left(s,Y^2_s,Z^2_s,(U^2_s)^{\rho}(\cdot), \mathbb{P}_{Y^1_s}\right)-g^1\left(s,Y^2_s,Z^2_s,(U^2_s)^{\rho}(\cdot), \mathbb{P}_{Y^2_s}\right)}{\|(\delta Y_s)^+\|_{L^2(\Omega)}}{\textbf{1}}_{\{\|(\delta Y_s)^+\|_{L^2(\Omega)}\neq0\}}.
\end{aligned}
\end{equation*}
It follows from  Assumption \ref{as4.1} and  conditions of Theorem \ref{th5.1} that
\begin{equation*}
\begin{aligned}
&|\delta _y g^1_s|\vee \delta _\mu g^1_s\leq L, \ \ \ \ |\delta _z g^1_s|\leq \phi(M_1)(1+|Z^1_s|+|Z^2_s|),\ \ \ \ |\delta _u g^1_s\rho_s(x)|\leq (1-\delta_{2M_1})(1\wedge |x|).
\end{aligned}
\end{equation*}
where $M_1$ is given in \eqref{4.1}. Furthermore, since $Z^1\cdot B$ and $Z^2\cdot B$ are BMO-martingales and the process $\delta _u g^1\rho$ is bounded as noted above, it follows that $\delta _z g^1\cdot B+\delta _u g^1\rho\cdot \tilde{N}$ is a BMO martingale.
 Moreover, by Lemma \ref{le2.2}, the Dol$\mathrm{\acute{e}}$ans-Dade exponential
$$(\mathcal{E}_t)_{t\in[0,T]}:=\mathcal{E}\big(\delta _z g^1\cdot B+\delta _u g^1\rho\cdot\tilde{N}\big)$$
is a uniformly integrable martingale. Consequently, $d\mathbb{Q} = \mathcal{E}_T d\mathbb{P}$ defines a probability measure equivalent to $\mathbb{P}$. Under $\mathbb{Q}$, the processes
$$dB^{\mathbb{Q}}_s:=dB_s-\delta _z g^1_sds,\ \ \tilde{N}^{\mathbb{Q}}(ds,dx):=\tilde{N}(ds,dx)-\delta _u g^1_s\rho_s(x)\lambda(dx)ds$$
are respectively a Brownian motion and a compensated random measure with compensator $(1+\delta_u g^1_s \rho_s(x))\lambda(dx)ds$. Therefore, under $\mathbb{Q}$, equation \eqref{eq5.2} can be rewritten as follows: for all $t\in[0,T]$,
\begin{equation*}
\begin{aligned}
\delta Y_t\le \delta \eta+\int_t^T\!\big(\delta _y g^1_s \delta Y_s+\delta _\mu g^1_s \|(\delta Y_s)^+\|_{L^2(\Omega)}+\delta g_s\big)ds
-\int_t^T \!\delta Z_sdB^{\mathbb{Q}}_s-\int_t^T\! \int_E\! \delta U_s(x)\tilde{N}^{\mathbb{Q}}(ds,dx).
\end{aligned}
\end{equation*}
For all $s\in [t,T]$, we define $R_s := \exp\{\int_t^s \delta_y g^1_r dr\}$. Then, applying It$\mathrm{\hat{o}}$'s formula to $R_s \delta Y_s$ together with  inequalities $\delta \eta \le 0$ and $\delta g_s \le 0$ yields
\begin{equation*}
\begin{aligned}
\delta Y_t\le\ &R_T\delta \eta+\int_t^T\! R_s\big[\delta _\mu g^1_s \|(\delta Y_s)^+\|_{L^2(\Omega)}+\delta g_s\big]ds
 -\int_t^T\! R_s\delta Z_sdB^{\mathbb{Q}}_s-\int_t^T\! \int_E\! R_s\delta U_s(x)\tilde{N}^{\mathbb{Q}}(ds,dx)\\
\leq\ & \int_t^T R_s\delta _\mu  g^1_s \|(\delta Y_s)^+\|_{L^2(\Omega)}ds -\int_t^T R_s\delta Z_sdB^{\mathbb{Q}}_s-\int_t^T \int_E R_s\delta U_s(x)\tilde{N}^{\mathbb{Q}}(ds,dx).
\end{aligned}
\end{equation*}
Taking the conditional expectation $\mathbb{E}^{\mathbb{Q}}_t[\cdot]$  on both sides of the above inequality, we obtain
$$\delta Y_t\leq \mathbb{E}^{\mathbb{Q}}_t\Big[\int_t^T R_s\delta _\mu  g^1_s \|(\delta Y_s)^+\|_{L^2(\Omega)}ds\Big]\leq Le^{TL}\int_t^T\big(\mathbb{E}[|(\delta Y_s)^+|^2]\big)^{\frac{1}{2}}ds.$$
Using H$\mathrm{\ddot{o}}$lder's inequality, we get
$$|(\delta Y_t)^+|^2\leq L^2e^{2TL}\Big(\int_t^T\big(\mathbb{E}[|(\delta Y_s)^+|^2]\big)^{\frac{1}{2}}ds\Big)^2\leq TL^2e^{2TL}\int_t^T\mathbb{E}[|(\delta Y_s)^+|^2]ds.$$
Taking the expectation $\mathbb{E}[\cdot]$ on both sides of the above inequality and applying Gronwall's lemma yields
 $Y^1_t\leq Y^2_t,\ \forall t\in[0,T],\ \mathbb{P}\mbox{-a.s.}$ This completes the proof.
\end{proof}

\begin{appendices}
\section{}
The following lemma allows us to establish the existence and the uniqueness of solutions to equation \eqref{eq3.2}.
\begin{lemma}\label{leA.1}
Assume that the generator $f:=f(\omega,t, z,u):\Omega \times[0,T] \times\mathbb{R}^d \times \mathcal{L}^{2}(\lambda)\rightarrow \mathbb{R}$ is  $\mathcal{P}^0$-measurable  for each fixed $(z,u)\in\mathbb{R}^d \times \mathcal{L}^{2}(\lambda)$ and that the terminal $\eta:\Omega\rightarrow \mathbb{R}$ is $\mathcal{F}_T$-measurable. $\eta$ and $f$ satisfy the following conditions:\\
$\mathrm{(1)}$  There exists a predictable  process $A:\Omega\times [0,T]\rightarrow [0,+\infty)$ such that, $d\mathbb{P}\times dt\mbox{-a.e.}$, for all $t\in [0,T],\ z\in \mathbb{R}^d,$ and $ u\in \mathcal{L}^{2}(\lambda)$,
\begin{equation*}
\begin{aligned}
-A _{t}-\frac{\gamma}{2}|z|^{2}-[-u]_{\gamma}\leq f(t,z,u)\leq A _{t}+\frac{\gamma}{2}|z|^{2}+[u]_{\gamma}.
\end{aligned}
\end{equation*}
$\mathrm{(2)}$\
There exists a constant $C>0$ and a function $d:[0,T]\rightarrow [0,+\infty)$ with $D:=\int_0^Td_s^2ds<+\infty$ such that, for each $M>0$,\ there exists  constant $ C_M\ge 0$ satisfying,  for all $t\in [0,T],\ z,\ \bar{z}\in \mathbb{R}^d,$ and $ u,\ \bar{u}\in \mathcal{L}^{2}(\lambda)\cap \mathcal{L}^{\infty}(\lambda)$ with $\|u\|_{\mathcal L^\infty(\lambda)}\vee\|\bar u\|_{\mathcal L^\infty(\lambda)}\le M$,
\begin{equation*}
\begin{aligned}
|f(t,z,u)-f(t,\bar{z},\bar{u})|\leq C(1+|z|+|\bar{z}|+d_t)|z-\bar{z}|+C_M\|u-\bar{u}\|_{\mathcal{L}^{2}(\lambda)},\  \ d\mathbb{P}\times dt\mbox{-a.e.}
\end{aligned}
\end{equation*}
$\mathrm{(3)}$\ For each $M>0$,  for all $t\in [0,T],\ z\in \mathbb{R}^d$ and $ u,\ \bar{u}\in \mathcal{L}^{2}(\lambda)\cap \mathcal{L}^{\infty}(\lambda)$ with $\|u\|_{\mathcal L^\infty(\lambda)}\vee\|\bar u\|_{\mathcal L^\infty(\lambda)}\le M$, there exists a $\mathcal{P}^0\otimes \mathcal{B}(E)$-measurable process $\Pi^{z,u,\bar{u}}$ satisfying
\begin{equation*}
\begin{aligned}
f(t,z,u)-f(t,z,\bar{u})
\leq\ &\int_E\Pi_t^{z,u,\bar{u}}(x)(u(x)-\bar{u}(x))\lambda(dx),\ \ \ d\mathbb{P}\times dt\mbox{-a.e.},
\end{aligned}
\end{equation*}
and $(-1+\delta_M)(1\wedge |x|)\leq\Pi_t^{z,u,\bar{u}}(x)\leq D_M(1\wedge |x|)$. Here, $D_M>0$ and $\delta_M>0$ depend only on $M$.\\
 $\mathrm{(4)}$ There exist constants $K> 0$ and $K_{A}> 0$ such that
 $$\| \eta\|_{\infty}\leq K\ \ \ \text{and}\ \ \ \Big\| \int _{0}^{T}A _{t} dt\Big\|_{\infty}\leq K_{A}.$$
Then the following BSDEJs
\begin{equation}\label{eqA.1}
\begin{aligned}
Y_t=\eta+\int_t^T f(s,Z_s,U_s)ds-\int_t^T Z_sdB_s-\int_t^T \int_E U_s(x)\tilde{N}(ds,dx),\ t\in [0,T]
\end{aligned}
\end{equation}
admits a unique  solution $(Y,Z,U)\in \mathbb{S}^{\infty}(0,T;\mathbb{R})\times \mathbb{H}^{2}_{BMO}(0,T;\mathbb{R}^d)\times \mathbb{J}^{2}_{BMO}(0,T;\mathbb{R})$.
\end{lemma}
\begin{proof}
We prove this lemma by using the random compensator version of \cite[Theorem 4.1]{FT18}, as stated in \cite[Remark 2.1]{FT18}. This extension  requires the noises have the weak predictable representation property. Below, we verify the conditions specified in \cite{FT18}.

Step 1. Construction a  measure with  finite second moment. For $x\in E$, we define
\begin{equation*}
H(x):=\left\{
\begin{aligned}
&x,\ \ \ \ \ \ \ \ \ \ \ \ \ \ \ \ \ \ \ \  \ \ |x|\le 1,\\
&\mathrm{sgn}(x)(2-|x|^{-1}),\ \ |x|>1.
\end{aligned}
\right.
\end{equation*}
Then $H$ is a bijection from $E$ onto $E_H:=(-2,2)\setminus\{0\}$.  For all $t\in[0,T]$,\ and $F\in \mathcal{B}(E_H)$,  we denote $N_H((0,t]\times F):=N((0,t]\times H^{-1}(F))$ and   $\lambda_H(F):=\lambda(H^{-1}(F))$, then measure $\lambda_H$ has a finite second moment and the compensator of $N_H$ is $\lambda_H(dy)dt$. For all $u \in \mathcal{L}^{2}(\lambda)$ and $y\in E_H$, we define $u_H(y):=u(H^{-1}(y))$. Moreover, for $v \in \mathcal{L}^{2}(\lambda_H)$, we introduce $f_H(t,z,v):=f(t,z,v\circ H)$. Then $\|u_H\|_{\mathcal{L}^{2}(\lambda_H)}=\|u\|_{\mathcal{L}^{2}(\lambda)}$ and $f_H(t,z,u_H)=f(t,z,u)$. For simplicity, we use $E,\ \lambda,\ N$ and $f$ to denote $E_H,\ \lambda_H,\ N_H$ and $f_H$.

Step 2. Introduce a new compensator.
For $t\in[0,T]$, we define $$a_t:=1+A_t+d^2_t,\ \ \ L_t:=\int_0^ta_sds,\ \ \ \sigma:=L_T,\ \  \ S:=T+K_A+D+1.$$
Then, $L$ is a continuous and increasing function, and $\sigma<S$. For $r\in[0,S]$, we define
$$\tau_r:=\inf\{t\in[0,T]: L_t\ge r\}\wedge T,$$
with $\inf \varnothing :=+\infty$. Therefore, $\tau_r=T$ for $r\ge \sigma$. Moreover, if let $c_r:=a_{\tau_r}^{-1}$ for $r<\sigma$, and $c_r=0$ for $r\ge \sigma$, then it's clearly that
\begin{equation}\label{eqA.2}
\begin{aligned}
d\tau_r=c_rdr,\ \ \ 0<c_r\le 1,\ \ \ c_rA_{\tau_r}\le 1,\ \ \ \sqrt{c_r}d_{\tau_r}\le 1
\end{aligned}
\end{equation}
 for $r<\sigma$.

In the filtration $\mathbb{G}:=\{\mathcal{G}_r\}_{r\in[0,S]}$, where $\mathcal{G}_r:=\mathcal{F}_{\tau_r}$, for $r\in[0,S]$, we denote
$$\widehat{B}_r:=\int_0^{\tau_r}\sqrt{a_s}dB_s,\ \ \ \ \widehat{N}((0,r]\times F):={N}((0,\tau_r]\times F).$$
Then the process $\widehat{B}$ is a continuous local martingale with $\langle\widehat{B}\rangle_r=r\wedge\sigma$, and for all  $F\in \mathcal{B}(E)$ with $\lambda(F)<\infty$, the process
$$\Big\{\widehat{N}((0,r]\times F)-\lambda(F)\int_0^rc_vdv\Big\}_{r\in[ 0,S]}=\big\{{N}((0,\tau_r]\times F)-\tau_r\lambda(F)\big\}_{r\in[ 0,S]}$$
is a $\mathbb{G}$-martingale. Hence $\widehat{B}$ is a Brownian motion stopped at $\sigma$, and $\widehat{N}$ is an integer-valued random measure with compensator $\widehat{\nu}_r(dx)dr$, where $\widehat{\nu}_r(dx):=c_r\lambda(dx)$. The compensated martingale measure of $\widehat{N}$ is
\begin{equation}\label{eqA.3}
\begin{aligned}
 \widetilde{\widehat{N}}(dr,dx):=\widehat{N}(dr,dx)-\widehat{\nu}_r(dx)dr.
\end{aligned}
\end{equation}
Moreover, for $[u]_{\gamma,\widehat{\nu}_r}:=\int_E\frac{e^{\gamma u(x)}-\gamma u(x)-1}{\gamma}\widehat{\nu}_r(dx)$, it's obvious that $[u]_{\gamma,\widehat{\nu}_r}=c_r[u]_\gamma$.

 In order to apply \cite[Theorem 4.1]{FT18} on the interval $[0,S]$, extend $\widehat{B}$ after $\sigma$ by independent Brownian increments to obtain a new Brownian motion $W$, and corresponding enlarge the filtration $\mathbb{G}$ (added only the information after time $\sigma$). 
 Because the time change preserves the weak predictable representation property of the original noises before $\sigma$ and we only  add independent Brownian increments  after $\sigma$,  it follows that $(W,\widetilde{\widehat{N}})$ has the required weak predictable representation property. Moreover, the independent extension does not change the compensator in \eqref{eqA.3}.

Step 3. Construction of a generator satisfying the requirements. For $(r,z,u)\in[0,S]\times \mathbb{R}^d\times \mathcal{L}^2(\widehat{\nu}_r)$, we define
\begin{equation}\label{eqA.4}
\widehat{f}(r,z,u):=\left\{
\begin{aligned}
&c_rf\Big(\tau_r,\frac{z}{\sqrt{c_r}},u\Big), \ \ \ \ \ \ r< \sigma,\\
&0,\ \ \ \ \ \ \ \ \ \ \ \ \ \ \ \ \ \ \ \ \ \  \  \ \  r\ge\sigma.
\end{aligned}
\right.
\end{equation}
 Then we claim that the generator $\widehat{f}$ satisfies Assumptions 3.1,\ 3.2 and 4.1 in \cite{FT18}. For Assumption 3.1, by Assumption (1) and $c_rA_{\tau_r}\le 1$ in \eqref{eqA.2}, we obtain
$$-1-\frac{\gamma}{2}|z|^{2}-[-u]_{\gamma,\widehat{\nu}_r}\leq\ \widehat{f}(r,z,u)\leq 1+\frac{\gamma}{2}|z|^{2}+[u]_{\gamma,\widehat{\nu}_r}.$$
For Assumption 3.2,  the Assumption (2) and $ 0<c_r\le 1, \ \sqrt{c_r}d_{\tau_r}\le 1$ in \eqref{eqA.2} imply that
\begin{equation*}
\begin{aligned}
|\widehat{f}(r,z,u)-\widehat{f}(r,\bar{z},\bar{u})|\le \ &C(\sqrt{c_r}+|z|+|\bar{z}|+\sqrt{c_r}d_{\tau_r})|z-\bar{z}|+c_rC_M\|u-\bar{u}\|_{\mathcal{L}^{2}(\lambda)}\\
\le\ &C(2+|z|+|\bar{z}|)|z-\bar{z}|+C_M\|u-\bar{u}\|_{\mathcal{L}^{2}(\widehat{\nu}_r)},
\end{aligned}
\end{equation*}
where we used $c_r\|v\|_{\mathcal{L}^{2}(\lambda)}=\sqrt{c_r}\|v\|_{\mathcal{L}^{2}(\widehat{\nu}_r)}$.\\
For Assumption 4.1,  we define $\widehat{\Pi}_r^{z,u,\bar{u}}:=\Pi_{\tau_r}^{z/\sqrt{c_r},u,\bar{u}}$ for $r<\sigma$ and $\widehat{\Pi}_r^{z,u,\bar{u}}:=0$ for $r\ge \sigma$, using  Assumption (3), we have
 $$\widehat{f}(r,z,u)-\widehat{f}(r,z,\bar{u})
\leq\ \int_E\widehat{\Pi}_r^{z,u,\bar{u}}(x)(u(x)-\bar{u}(x))\widehat{\nu}_r(dx),$$
and $(-1+\delta_M)(1\wedge |x|)\leq\widehat{\Pi}_t^{z,u,\bar{u}}(x)\leq D_M(1\wedge |x|)$. Here, $D_M>0$ and $\delta_M>0$ depend only on $M$.
Therefore, according to \cite[Remark 2.1]{FT18}, applying the random compensator version of \cite[Theorem 4.1]{FT18}, we obtain the following equation
\begin{equation}\label{eqA.5}
\begin{aligned}
\widehat{Y}_r=\eta+\int_r^S \widehat{f}(t,\widehat{Z}_t,\widehat{U}_t)dt-\int_r^S \widehat{Z}_tdW_t-\int_r^S \int_E \widehat{U}_t(x)\widetilde{\widehat{N}}(dt,dx),\ r\in [0,S]
\end{aligned}
\end{equation}
admits a unique  solution $(\widehat{Y},\widehat{Z},\widehat{U})\in \mathbb{S}^{\infty}(0,S;\mathbb{R})\times \mathbb{H}^{2}_{BMO}(0,S;\mathbb{R}^d)\times \mathbb{J}^{2}_{BMO}(0,S;\mathbb{R})$.

Step 4. The existence and the uniqueness of equation \eqref{eqA.1}. Since $\mathcal{G}_\sigma=\mathcal{F}_{\tau_\sigma}=\mathcal{F}_T$, so we know that $\eta$ is $\mathcal{G}_\sigma$-measurable. Taking conditional expectation $\mathbb{E}_{\mathcal{G}_r}[\cdot]$ in \eqref{eqA.5} and using $\widehat{f}(r,z,u)=0$ for $r\ge \sigma$, we get $\widehat{Y}_r=\eta$ and $\widehat{Z}_r=0$ for $r\ge \sigma$. We also choose   $\widehat{U}_r(x)=0$ for $r\ge \sigma$. For $t\in[0,T]$ and $x\in E$, we define
\begin{equation}\label{eqA.6}
\begin{aligned}
Y_t:=\widehat{Y}_{L_t},\ \ \ Z_t:=\sqrt{a_t}\widehat{Z}_{L_t},\ \ \ U_t(x):=\widehat{U}_{L_t}(x).
\end{aligned}
\end{equation}
Due to the filtration  before $\sigma$  unchanged, the processes in \eqref{eqA.6} are still adapted or predictable in the original filtration. By change of variables, the definition of $L$ and $c_{L_s}a_s=1$, for $t\in[0,T]$, we obtain
$$\int_{L_t}^\sigma\widehat{f}(r,\widehat{Z}_r,\widehat{U}_r)dr=\int_{t}^T c_{L_s}a_sf\Big(\tau_{L_s},\frac{\widehat{Z}_{L_s}}{\sqrt{c_{L_s}}},\widehat{U}_{L_s}\Big)ds=\int_t^Tf(s,\sqrt{a_s}\widehat{Z}_{L_s},\widehat{U}_{L_s})ds
=\int_t^Tf(s,Z_s,U_s)ds.$$
Using the definition of $W$, for $t\in[0,T]$, we have
$$\int_{L_t}^\sigma\widehat{Z}_rdW_r= \int_t^T\widehat{Z}_{L_s}\sqrt{a_s}dB_s=\int_t^TZ_sdB_s.$$
Similarly, by equation \eqref{eqA.3}, for $t\in[0,T]$, we get
\begin{equation*}
\begin{aligned}
\int_{L_t}^\sigma\int_E\widehat{U}_r(x)\widetilde{\widehat{N}}(dr,dx)=&\int_{L_t}^\sigma\int_E\widehat{U}_r(x)\big(\widehat{N}(dr,dx)-c_r\lambda(dx)dr\big)\\
=&\int_t^T\int_E\widehat{U}_{L_s}(x)\big(N(ds,dx)-c_{L_s}a_s\lambda(dx)ds\big)\\
=&\int_t^T\int_E\widehat{U}_{L_s}(x)\tilde{N}(ds,dx)
=\int_t^T\int_EU_s(x)\tilde{N}(ds,dx).
\end{aligned}
\end{equation*}
Combining this with equation \eqref{eqA.5}, we get equation \eqref{eqA.1} holds.  For each $\mathbb{F}$-stopping time $\theta\in[0,T]$ , we know that $L_\theta$ is a $\mathbb{G}$-stopping time. Due to
 $$\int_\theta^T|Z_s|^2ds=\int_{L_\theta}^\sigma|\widehat{Z}_r|^2dr,\ \ \
 \int_\theta^T\int_E|U_s(x)|^2\lambda(dx)ds=\int_{L_\theta}^\sigma\int_E|\widehat{U}_r(x)|^2\widehat{\nu}_r(dx)dr,$$
 and noting that the jumps size remains unchanged, it's easy to verify that $(Y,Z,U)\in \mathbb{S}^{\infty}(0,T;\mathbb{R})\times \mathbb{H}^{2}_{BMO}(0,T;\mathbb{R}^d)\times \mathbb{J}^{2}_{BMO}(0,T;\mathbb{R})$. Finally, since equation \eqref{eqA.5} has a unique solution, the solution of equation \eqref{eqA.1} is also unique. This completes the proof.

\end{proof}

\end{appendices}


\begin{thebibliography}{99}

\bibitem{BBP97}G. Barles, R. Buckdahn, E. Pardoux. Backward stochastic differential equations and integral-partial differential
 equations. \emph{Stochastics and Stochastics Reports}. 60, 57-83, 1997.

\bibitem{B06} D. Becherer. Bounded solutions to backward SDE's with jumps for utility optimization and indifference hedging. \emph{Annals of Applied Probability}. 16 (4), 2027-2054, 2006.

\bibitem{BH26} D. Becherer, S. Hesse. Common noise by random measures: constructing mean-field equilibria for competitive investment and hedging. \emph{arXiv:2408.01175v4}, 3, 2026.

\bibitem{B73} J.-M. Bismut. Conjugate convex functions in optimal stochastic control. \emph{Journal of Mathematical Analysis and Applications}. 44 (2), 384-404, 1973.


\bibitem{BH06} P. Briand, Y. Hu. BSDE with quadratic growth and unbounded terminal value.  \emph{Probability Theory and Related Fields}. 136 (4), 604-618, 2006.

\bibitem{BH08} P. Briand, Y. Hu. Quadratic BSDEs with convex generators and unbounded terminal condition.  \emph{Probability
Theory and Related Fields}. 141 (3-4), 543-567, 2008.

\bibitem{BDLP09} R. Buckdahn, B. Djehiche, J. Li, S. Peng. Mean-field backward stochastic differential equations. A limit approach. \emph{Annals of Probability}. 37 (4), 1524-1565, 2009.

\bibitem{BLP09} R. Buckdahn, J. Li, S. Peng. Mean-field backward stochastic differential equations and related partial differential equations. \emph{Stochastic Processes and their Applications}. 119 (10), 3133-3154, 2009.

\bibitem{BLPR17} R. Buckdahn, J. Li, S. Peng, C. Rainer. Mean-field stochastic differential equations and associated PDEs. \emph{Annals of Probability}. 45 (2), 824-878, 2017.


\bibitem{FHT20}S. Fan, Y. Hu, S. Tang. On the uniqueness of solutions to quadratic BSDEs with non-convex generators and unbounded terminal conditions. \emph{Comptes Rendus Mathematique}. 358 (2), 227-235, 2020.



\bibitem{FJ13} S. Fan, L. Jiang. Multidimensional BSDEs with weak monotonicity and general growth generators.
\emph{Acta Mathematica Sinica, English Series}.  29 (10), 1885-1906, 2013.

\bibitem{FT18}M. Fujii, A. Takahashi.  Quadratic-exponential growth BSDEs with jumps and their Malliavin's differentiability. \emph{Stochastic Processes and their Applications}. 128 (6), 2083-2130, 2018.

\bibitem{FT19}M. Fujii, A. Takahashi. Anticipated Backward SDEs with jumps and quadratic-exponential growth drivers. \emph{Stochastics and Dynamics}. 19 (3), 1950020, 2019.


\bibitem{HHTW25}T. Hao, Y. Hu, S. Tang, J. Wen. Mean-field backward stochastic differential equations and nonlocal PDEs with quadratic growth. \emph{Annals of Applied Probability}. 35 (3), 2128-2174, 2025.

\bibitem{HIM05} Y. Hu, P. Imkeller, M. Muller. Utility maximization in incomplete markets. \emph{Annals of Applied Probability}. 15 (3), 1691-1712, 2005.


\bibitem{HLW21}Y. Hu, X. Li, J. Wen. Anticipated backward stochastic differential equations with quadratic growth. \emph{Journal of Differential Equations}. 270, 1298-1331, 2021.

\bibitem{HT16} Y. Hu, S. Tang. Multi-dimensional backward stochastic differential equations of diagonally quadratic generators. \emph{Stochastic Processes and their Applications}. 126 (4), 1066-1086, 2016.

\bibitem{HWX24}Y. Hu,  J. Wen, J. Xiong. Backward doubly stochastic differential equations and SPDEs with quadratic growth. \emph{Stochastic Processes and their Applications}. 175, 104405, 2024.

\bibitem{HMC06}M. Huang, R.P. Malham$\mathrm{\acute{e}}$, P.E. Caines. Large population stochastic dynamic games:
closed-loop McKean-Vlasov systems and the Nash certainty equivalence principle. \emph{Communications in Information and Systems}. 6 (3), 221-252, 2006.


\bibitem{JLW24} W. Jiang, J. Li, Q. Wei.  General mean-field BSDEs with diagonally quadratic generator in multi-dimension. \emph{Discrete and Continuous Dynamical Systems}. 44 (10), 2957-2984, 2024.

\bibitem{K56}M. Kac. Foundations of kinetic theory. \emph{Proceedings of the Third Berkeley Symposium on Mathematical Statistics and Probability}.  171-197, 1956.

\bibitem{KT03} M. Kohlmann, S. Tang. Minimization of risk and linear quadratic optimal control theory. \emph{SIAM Journal on Control and Optimization}. 42 (3), 1118-1142, 2003.


\bibitem{KPQ97}N. El Karoui, S. Peng, M.C. Quenez. Backward stochastic differential equations in finance. \emph{Mathematical Finance}. 7 (1), 1-71, 1997.

\bibitem{K79} N. Kazamaki. Transformation of $H^p$-martingales by a change of law. \emph{Zeitschrift f$\ddot{u}$r Wahrscheinlichkeitstheorie und verwandte Gebiete}. 46, 343-349, 1979.

\bibitem{K94}N. Kazamaki. \emph{Continuous exponential martingales and BMO}. Springer-Verlag, Berlin. 1994.

\bibitem{KPZ15} N. Kazi-Tani, D. Possama, C. Zhou. Quadratic BSDEs with jumps: a fixed-point approach. \emph{Electron Journal of Probability}. 20 (66), 1-28, 2015.

\bibitem{K00} M. Kobylanski. Backward stochastic differential equations and partial differential equations with
quadratic growth. \emph{Annals of Probability}. 28 (2), 558-602, 2000.

\bibitem{LL07} J.M. Lasry, P.L. Lions. Mean field games. \emph{Japanese Journal of Mathematics}. 2, 229-260, 2007.

\bibitem{LM97}J.-P. Lepeltier, J. San Martin. Backward stochastic differential equations with continuous coefficient. \emph{Statistics \emph{\&} Probability Letters}. 32 (4), 425-430, 1997.


\bibitem{L18} J. Li. Mean-field forward and backward SDEs with jumps and associated nonlocal quasi-linear integral-PDEs. \emph{Stochastic Processes and their Applications}. 128 (9), 3118-3180, 2018.

\bibitem{LLZ18} J. Li, H. Liang, X. Zhang. General mean-field BSDEs with continuous coefficients. \emph{Journal of Mathematical Analysis and Applications}. 466 (1), 264-280, 2018.
\bibitem{LQ11} A. Lim, M.C. Quenez. Exponential utility maximization and indifference pricing in an incomplete market with defaults. \emph{Electronic Journal of Probability}. 16, 1434-1464, 2011.

\bibitem{M09}M.A. Morlais. Utility maximization in a jump market model. \emph{Stochastics: An International Journal of Probability and Stochastics Processes}. 81 (1), 1-27, 2009.

\bibitem{NO25}  K. Nasroallah, Y. Ouknine. Backward stochastic differential equations with central value reflection. \emph{Probability, Uncertainty and Quantitative Risk}. 10 (2), 293-318, 2025.

\bibitem{PP90}  E. Pardoux, S. Peng. Adapted solution of a backward stochastic differential equation. \emph{Systems \emph{\&} Control Letters}. 14 (1), 55-61, 1990.

\bibitem{PP92}  E. Pardoux, S. Peng.  Backward stochastic differential equations and quasilinear parabolic partial differential equations. \emph{Stochastic Partial Differential Equations and their Applications}. 200-217, 1992.

\bibitem{P92} S. Peng. A generalized dynamic programming principle and Hamilton-Jacobi-Bellman equation. \emph{Stochastics and Stochastic Reports}. 38, 119-134, 1992.



\bibitem{TL94}S. Tang, X. Li. Necessary conditions for optimal control of stochastic systems with random jumps. \emph{SIAM Journal on Control and Optimization}. 32 (5), 1447-1475, 1994.

\bibitem{S97}R. Situ. On solutions of backward stochastic differential equations with jumps and applications. \emph{Stochastic Processes and their Applications}. 66 (2), 209-236, 1997.




\end{thebibliography}
\end{document}